%Paper: alg-geom/9511004
%From: Gerard.LAUMON@math.u-psud.fr ( Gerard LAUMON ORSAY )
%Date: Tue, 7 Nov 1995 15:37:06 +0100

\magnification=1200
% =======================================================
% ================== Reglages generaux ==================
\hsize=13.50cm
\vsize=18cm
\parindent=12pt   \parskip=0pt
\pageno=1

\ifnum\mag=\magstep1
\hoffset=-2mm   % offset horizontal en \magnification=1200
\voffset=.8cm   % offset horizontal en \magnification=1200
\fi

% -------------------  Reglages --------------------------
% -------- auquels il ne vaut mieux pas toucher ----------

\pretolerance=500 \tolerance=1000  \brokenpenalty=5000

% -------------- Debut des macros privees ----------------
\catcode`\@=11
% ---------------------Les fontes ------------------------

\font\eightrm=cmr8         \font\eighti=cmmi8
\font\eightsy=cmsy8        \font\eightbf=cmbx8
\font\eighttt=cmtt8        \font\eightit=cmti8
\font\eightsl=cmsl8        \font\sixrm=cmr6
\font\sixi=cmmi6           \font\sixsy=cmsy6
\font\sixbf=cmbx6

% Fontes AMS

\font\tengoth=eufm10
\font\tenbboard=msbm10
\font\eightgoth=eufm10 at 8pt
\font\eightbboard=msbm10 at 8 pt
\font\sevengoth=eufm7
\font\sevenbboard=msbm7
\font\sixgoth=eufm7 at 6 pt
\font\fivegoth=eufm5

 \font\tencyr=wncyr10

\font\eightcyr=wncyr10 at 8 pt

\font\sevencyr=wncyr10 at 7 pt

\font\sixcyr=wncyr10 at 6 pt

% Pour que les accents se placent correctement
%en mode math en corps 8 et 6

\skewchar\eighti='177 \skewchar\sixi='177
\skewchar\eightsy='60 \skewchar\sixsy='60

% Nouvelles familles pour les maths

\newfam\gothfam           \newfam\bboardfam
\newfam\cyrfam

\def\tenpoint{%
  \textfont0=\tenrm \scriptfont0=\sevenrm
  \scriptscriptfont0=\fiverm
  \def\rm{\fam\z@\tenrm}%
  \textfont1=\teni  \scriptfont1=\seveni
  \scriptscriptfont1=\fivei
  \def\oldstyle{\fam\@ne\teni}\let\old=\oldstyle
  \textfont2=\tensy \scriptfont2=\sevensy
  \scriptscriptfont2=\fivesy
  \textfont\gothfam=\tengoth \scriptfont\gothfam=\sevengoth
  \scriptscriptfont\gothfam=\fivegoth
  \def\goth{\fam\gothfam\tengoth}%
  \textfont\bboardfam=\tenbboard
  \scriptfont\bboardfam=\sevenbboard
  \scriptscriptfont\bboardfam=\sevenbboard
  \def\bb{\fam\bboardfam\tenbboard}%
  \textfont\cyrfam=\tencyr \scriptfont\cyrfam=\sevencyr
  \scriptscriptfont\cyrfam=\sixcyr
  \def\cyr{\fam\cyrfam\tencyr}%
  \textfont\itfam=\tenit
  \def\it{\fam\itfam\tenit}%
  \textfont\slfam=\tensl
  \def\sl{\fam\slfam\tensl}%
  \textfont\bffam=\tenbf \scriptfont\bffam=\sevenbf
  \scriptscriptfont\bffam=\fivebf
  \def\bf{\fam\bffam\tenbf}%
  \textfont\ttfam=\tentt
  \def\tt{\fam\ttfam\tentt}%
  \abovedisplayskip=12pt plus 3pt minus 9pt
  \belowdisplayskip=\abovedisplayskip
  \abovedisplayshortskip=0pt plus 3pt
  \belowdisplayshortskip=4pt plus 3pt
  \smallskipamount=3pt plus 1pt minus 1pt
  \medskipamount=6pt plus 2pt minus 2pt
  \bigskipamount=12pt plus 4pt minus 4pt
  \normalbaselineskip=12pt
  \setbox\strutbox=\hbox{\vrule height8.5pt depth3.5pt width0pt}%
  \let\bigf@nt=\tenrm       \let\smallf@nt=\sevenrm
  \normalbaselines\rm}

\def\eightpoint{%
  \textfont0=\eightrm \scriptfont0=\sixrm
  \scriptscriptfont0=\fiverm
  \def\rm{\fam\z@\eightrm}%
  \textfont1=\eighti  \scriptfont1=\sixi
  \scriptscriptfont1=\fivei
  \def\oldstyle{\fam\@ne\eighti}\let\old=\oldstyle
  \textfont2=\eightsy \scriptfont2=\sixsy
  \scriptscriptfont2=\fivesy
  \textfont\gothfam=\eightgoth \scriptfont\gothfam=\sixgoth
  \scriptscriptfont\gothfam=\fivegoth
  \def\goth{\fam\gothfam\eightgoth}%
  \textfont\cyrfam=\eightcyr \scriptfont\cyrfam=\sixcyr
  \scriptscriptfont\cyrfam=\sixcyr
  \def\cyr{\fam\cyrfam\eightcyr}%
  \textfont\bboardfam=\eightbboard
  \scriptfont\bboardfam=\sevenbboard
  \scriptscriptfont\bboardfam=\sevenbboard
  \def\bb{\fam\bboardfam}%
  \textfont\itfam=\eightit
  \def\it{\fam\itfam\eightit}%
  \textfont\slfam=\eightsl
  \def\sl{\fam\slfam\eightsl}%
  \textfont\bffam=\eightbf \scriptfont\bffam=\sixbf
  \scriptscriptfont\bffam=\fivebf
  \def\bf{\fam\bffam\eightbf}%
  \textfont\ttfam=\eighttt
  \def\tt{\fam\ttfam\eighttt}%
  \abovedisplayskip=9pt plus 3pt minus 9pt
  \belowdisplayskip=\abovedisplayskip
  \abovedisplayshortskip=0pt plus 3pt
  \belowdisplayshortskip=3pt plus 3pt
  \smallskipamount=2pt plus 1pt minus 1pt
  \medskipamount=4pt plus 2pt minus 1pt
  \bigskipamount=9pt plus 3pt minus 3pt
  \normalbaselineskip=9pt
  \setbox\strutbox=\hbox{\vrule height7pt depth2pt width0pt}%
  \let\bigf@nt=\eightrm     \let\smallf@nt=\sixrm
  \normalbaselines\rm}

\tenpoint

% Definition des petites capitales qui reagissent au
%\tenpoint et \eightpoint
% La syntaxe est celle d'un changement de fonte :
% {\pc FERMAT}, {\pc EUCLIDE} et {\pc G\"ODEL}.

\def\pc#1{\bigf@nt#1\smallf@nt}         \def\pd#1 {{\pc#1} }

% ----------------dactylographie francaise ------------------

\catcode`\;=\active
\def;{\relax\ifhmode\ifdim\lastskip>\z@\unskip\fi
\kern\fontdimen2  -1.2 \fontdimen3 \string;}

\catcode`\:=\active
\def:{\relax\ifhmode\ifdim\lastskip>\z@
\unskip\fi\penalty\@M\ \fi\string:}

\catcode`\!=\active
\def!{\relax\ifhmode\ifdim\lastskip>\z@
\unskip\fi\kern\fontdimen2  -1.1
\fontdimen3 \string!}

\catcode`\?=\active
\def?{\relax\ifhmode\ifdim\lastskip>\z@
\unskip\fi\kern\fontdimen2  -1.1
\fontdimen3 \string?}

\def\^#1{\if#1i{\accent"5E\i}\else{\accent"5E #1}\fi}
\def\"#1{\if#1i{\accent"7F\i}\else{\accent"7F #1}\fi}

\frenchspacing

% ------------------ Le format de sortie ---------------------
% Haut et bas de page

\newtoks\auteurcourant      \auteurcourant={\hfil}
\newtoks\titrecourant       \titrecourant={\hfil}

\newtoks\hautpagetitre      \hautpagetitre={\hfil}
\newtoks\baspagetitre       \baspagetitre={\hfil}

\newtoks\hautpagegauche
\hautpagegauche={\eightpoint\rlap{\folio}\hfil\the\auteurcourant
\hfil}
\newtoks\hautpagedroite
\hautpagedroite={\eightpoint\hfil\the\titrecourant\hfil
\llap{\folio}}

\newtoks\baspagegauche      \baspagegauche={\hfil}
\newtoks\baspagedroite      \baspagedroite={\hfil}

\newif\ifpagetitre          \pagetitretrue

% \nopagenumbers : c'est un peu violent, mais a marche. Alors ...

\headline={\ifpagetitre\the\hautpagetitre
\else\ifodd\pageno\the\hautpagedroite\else\the\hautpagegauche
\fi\fi}

\footline={\ifpagetitre\the\baspagetitre\else
\ifodd\pageno\the\baspagedroite\else\the\baspagegauche\fi\fi
\global\pagetitrefalse}

% Redefinition de \raggedbottom pour avoir plus de mou
%en bas de page
% (necesssaire quand il y a beaucoup de grumeaux, des grosses
% formules centrees et pas beaucoup de texte entre)

\def\raggedbottom{\topskip 10pt plus 36pt\r@ggedbottomtrue}

% --------------------- Macros de mise en page --------------

% Un point-tiret

\def\pointir{\unskip . --- \ignorespaces}

% Macros Bigbreak et \Medbreak pour que les blancs verticaux
%ne s'ajoutent pas

\def\Medbreak{\vskip-\lastskip\medbreak}

\def\rem#1\endrem{%
\Medbreak
{\it#1\unskip} : }

% ------------------------------------------------------------
% enonces de theoremes avec numerotation apres
% #1 = THEOREME, COROLLAIRE, etc.
% #2 = numero (par exemple 3, 3.1, etc.)
% #3 = l'enonce du th proprememnt dit.

\long\def\th#1 #2\enonce#3\endth{%
   \Medbreak
   {\pc#1} {#2\unskip}\pointir{\it #3}\medskip}

\long\def\tha#1 #2\enonce#3\endth{%
   \Medbreak
   {\pc#1} {#2\unskip}\par\nobreak{\it #3}\medskip}

% ------------------------------------------------------------

% les differents retraits, voir aussi \item

\def\decale#1{\smallbreak\hskip 28pt\llap{#1}\kern 5pt}
\def\decaledecale#1{\smallbreak\hskip 34pt\llap{#1}\kern 5pt}

% ------------------------------------------------------------
% ----------------- ce que Knuth n'a pas fait ----------------
% ------------------------------------------------------------

% pour avoir des messages raisonnables avec les
%lettres accentu\'ees

\let\@ldmessage=\message

\def\message#1{{\def\pc{\string\pc\space}%
                \def\'{\string'}\def\`{\string`}%
                \def\^{\string^}\def\"{\string"}%
                \@ldmessage{#1}}}

% --------------------- Divers gadgets -----------------------

\def\qed{\raise -2pt\hbox{\vrule\vbox to 10pt{\hrule width 4pt
                 \vfill\hrule}\vrule}}

\def\cqfd{\unskip\penalty 500\quad\qed\medbreak}

% ------------------------------------------------------------
% fin des macros privees
% ------------------------------------------------------------
\catcode`\@=12

% pour qu'il la ferme
\showboxbreadth=-1  \showboxdepth=-1

\expandafter\ifx\csname amssym.def\endcsname\relax \else
\endinput\fi
%
%  Store the catcode of the @ in the csname so that it can be restored later.
\expandafter\edef\csname amssym.def\endcsname{%
       \catcode`\noexpand\@=\the\catcode`\@\space}
%  Set the catcode to 11 for use in private control sequence names.
\catcode`\@=11
%
%  Include all definitions related to the fonts msam, msbm and eufm, so that
%  when this file is used by itself, the results with respect to those fonts
%  are equivalent to what they would have been using AMS-TeX.
%  Most symbols in fonts msam and msbm are defined using \newsymbol;
%  however, a few symbols that replace composites defined in plain must be
%  defined with \mathchardef.

\def\undefine#1{\let#1\undefined}
\def\newsymbol#1#2#3#4#5{\let\next@\relax
 \ifnum#2=\@ne\let\next@\msafam@\else
 \ifnum#2=\tw@\let\next@\msbfam@\fi\fi
 \mathchardef#1="#3\next@#4#5}
\def\mathhexbox@#1#2#3{\relax
 \ifmmode\mathpalette{}{\m@th\mathchar"#1#2#3}%
 \else\leavevmode\hbox{$\m@th\mathchar"#1#2#3$}\fi}
\def\hexnumber@#1{\ifcase#1 0\or 1\or 2\or 3\or 4\or 5\or 6\or 7\or 8\or
 9\or A\or B\or C\or D\or E\or F\fi}

\font\tenmsa=msam10
\font\sevenmsa=msam7
\font\fivemsa=msam5
\newfam\msafam
\textfont\msafam=\tenmsa
\scriptfont\msafam=\sevenmsa
\scriptscriptfont\msafam=\fivemsa
\edef\msafam@{\hexnumber@\msafam}
\mathchardef\dabar@"0\msafam@39
\def\dashrightarrow{\mathrel{\dabar@\dabar@\mathchar"0\msafam@4B}}
\def\dashleftarrow{\mathrel{\mathchar"0\msafam@4C\dabar@\dabar@}}

\def\ulcorner{\delimiter"4\msafam@70\msafam@70 }
\def\urcorner{\delimiter"5\msafam@71\msafam@71 }
\def\llcorner{\delimiter"4\msafam@78\msafam@78 }
\def\lrcorner{\delimiter"5\msafam@79\msafam@79 }
%    Note that there should not be a final space after the digits for a
%    \mathhexbox@.
\def\yen{{\mathhexbox@\msafam@55}}
\def\checkmark{{\mathhexbox@\msafam@58}}
\def\circledR{{\mathhexbox@\msafam@72}}
\def\maltese{{\mathhexbox@\msafam@7A}}

\font\tenmsb=msbm10
\font\sevenmsb=msbm7
\font\fivemsb=msbm5
\newfam\msbfam
\textfont\msbfam=\tenmsb
\scriptfont\msbfam=\sevenmsb
\scriptscriptfont\msbfam=\fivemsb
\edef\msbfam@{\hexnumber@\msbfam}

\def\widehat#1{\setbox\z@\hbox{$\m@th#1$}%
 \ifdim\wd\z@>\tw@ em\mathaccent"0\msbfam@5B{#1}%
 \else\mathaccent"0362{#1}\fi}
\def\widetilde#1{\setbox\z@\hbox{$\m@th#1$}%
 \ifdim\wd\z@>\tw@ em\mathaccent"0\msbfam@5D{#1}%
 \else\mathaccent"0365{#1}\fi}
\font\teneufm=eufm10
\font\seveneufm=eufm7
\font\fiveeufm=eufm5
\newfam\eufmfam
\textfont\eufmfam=\teneufm
\scriptfont\eufmfam=\seveneufm
\scriptscriptfont\eufmfam=\fiveeufm

\let\goth\frak

%  Restore the catcode value for @ that was previously saved.
\csname amssym.def\endcsname
\expandafter\ifx\csname pre amssym.tex at\endcsname\relax \else
\endinput\fi
%%  Otherwise we store the catcode of the @ in the csname.
\expandafter\chardef\csname pre amssym.tex at\endcsname=\the
\catcode`\@
%%  Set the catcode to 11 for use in private control sequence names.
\catcode`\@=11
%%  Load amssym.def if necessary: If \newsymbol is undefined, do nothing
%%  and the following \input statement will be executed; otherwise
%%  change \input to a temporary no-op.
\begingroup\ifx\undefined\newsymbol \else\def\input#1 {\endgroup}\fi
\input amssym.def \relax
%%  Most symbols in fonts msam and msbm are defined using \newsymbol.
%% A few
%%  that are delimiters or otherwise require special treatment have already
%%  been defined as soon as the fonts were loaded.  Finally, a few symbols
%%  that replace composites defined in plain must be undefined first.
\newsymbol\boxdot 1200
\newsymbol\boxplus 1201
\newsymbol\boxtimes 1202
\newsymbol\square 1003
\newsymbol\blacksquare 1004
\newsymbol\centerdot 1205
\newsymbol\lozenge 1006
\newsymbol\blacklozenge 1007
\newsymbol\circlearrowright 1308
\newsymbol\circlearrowleft 1309
\undefine\rightleftharpoons
\newsymbol\rightleftharpoons 130A
\newsymbol\leftrightharpoons 130B
\newsymbol\boxminus 120C
\newsymbol\Vdash 130D
\newsymbol\Vvdash 130E
\newsymbol\vDash 130F
\newsymbol\twoheadrightarrow 1310
\newsymbol\twoheadleftarrow 1311
\newsymbol\leftleftarrows 1312
\newsymbol\rightrightarrows 1313
\newsymbol\upuparrows 1314
\newsymbol\downdownarrows 1315
\newsymbol\upharpoonright 1316
 
\newsymbol\downharpoonright 1317
\newsymbol\upharpoonleft 1318
\newsymbol\downharpoonleft 1319
\newsymbol\rightarrowtail 131A
\newsymbol\leftarrowtail 131B
\newsymbol\leftrightarrows 131C
\newsymbol\rightleftarrows 131D
\newsymbol\Lsh 131E
\newsymbol\Rsh 131F
\newsymbol\rightsquigarrow 1320
\newsymbol\leftrightsquigarrow 1321
\newsymbol\looparrowleft 1322
\newsymbol\looparrowright 1323
\newsymbol\circeq 1324
\newsymbol\succsim 1325
\newsymbol\gtrsim 1326
\newsymbol\gtrapprox 1327
\newsymbol\multimap 1328
\newsymbol\therefore 1329
\newsymbol\because 132A
\newsymbol\doteqdot 132B
 
\newsymbol\triangleq 132C
\newsymbol\precsim 132D
\newsymbol\lesssim 132E
\newsymbol\lessapprox 132F
\newsymbol\eqslantless 1330
\newsymbol\eqslantgtr 1331
\newsymbol\curlyeqprec 1332
\newsymbol\curlyeqsucc 1333
\newsymbol\preccurlyeq 1334
\newsymbol\leqq 1335
\newsymbol\leqslant 1336
\newsymbol\lessgtr 1337
\newsymbol\backprime 1038
\newsymbol\risingdotseq 133A
\newsymbol\fallingdotseq 133B
\newsymbol\succcurlyeq 133C
\newsymbol\geqq 133D
\newsymbol\geqslant 133E
\newsymbol\gtrless 133F
\newsymbol\sqsubset 1340
\newsymbol\sqsupset 1341
\newsymbol\vartriangleright 1342
\newsymbol\vartriangleleft 1343
\newsymbol\trianglerighteq 1344
\newsymbol\trianglelefteq 1345
\newsymbol\bigstar 1046
\newsymbol\between 1347
\newsymbol\blacktriangledown 1048
\newsymbol\blacktriangleright 1349
\newsymbol\blacktriangleleft 134A
\newsymbol\vartriangle 134D
\newsymbol\blacktriangle 104E
\newsymbol\triangledown 104F
\newsymbol\eqcirc 1350
\newsymbol\lesseqgtr 1351
\newsymbol\gtreqless 1352
\newsymbol\lesseqqgtr 1353
\newsymbol\gtreqqless 1354
\newsymbol\Rrightarrow 1356
\newsymbol\Lleftarrow 1357
\newsymbol\veebar 1259
\newsymbol\barwedge 125A
\newsymbol\doublebarwedge 125B
\undefine\angle
\newsymbol\angle 105C
\newsymbol\measuredangle 105D
\newsymbol\sphericalangle 105E
\newsymbol\varpropto 135F
\newsymbol\smallsmile 1360
\newsymbol\smallfrown 1361
\newsymbol\Subset 1362
\newsymbol\Supset 1363
\newsymbol\Cup 1264
 
\newsymbol\Cap 1265
 
\newsymbol\curlywedge 1266
\newsymbol\curlyvee 1267
\newsymbol\leftthreetimes 1268
\newsymbol\rightthreetimes 1269
\newsymbol\subseteqq 136A
\newsymbol\supseteqq 136B
\newsymbol\bumpeq 136C
\newsymbol\Bumpeq 136D
\newsymbol\lll 136E
 
\newsymbol\ggg 136F
 
\newsymbol\circledS 1073
\newsymbol\pitchfork 1374
\newsymbol\dotplus 1275
\newsymbol\backsim 1376
\newsymbol\backsimeq 1377
\newsymbol\complement 107B
\newsymbol\intercal 127C
\newsymbol\circledcirc 127D
\newsymbol\circledast 127E
\newsymbol\circleddash 127F
\newsymbol\lvertneqq 2300
\newsymbol\gvertneqq 2301
\newsymbol\nleq 2302
\newsymbol\ngeq 2303
\newsymbol\nless 2304
\newsymbol\ngtr 2305
\newsymbol\nprec 2306
\newsymbol\nsucc 2307
\newsymbol\lneqq 2308
\newsymbol\gneqq 2309
\newsymbol\nleqslant 230A
\newsymbol\ngeqslant 230B
\newsymbol\lneq 230C
\newsymbol\gneq 230D
\newsymbol\npreceq 230E
\newsymbol\nsucceq 230F
\newsymbol\precnsim 2310
\newsymbol\succnsim 2311
\newsymbol\lnsim 2312
\newsymbol\gnsim 2313
\newsymbol\nleqq 2314
\newsymbol\ngeqq 2315
\newsymbol\precneqq 2316
\newsymbol\succneqq 2317
\newsymbol\precnapprox 2318
\newsymbol\succnapprox 2319
\newsymbol\lnapprox 231A
\newsymbol\gnapprox 231B
\newsymbol\nsim 231C
\newsymbol\ncong 231D
\newsymbol\diagup 201E
\newsymbol\diagdown 201F
\newsymbol\varsubsetneq 2320
\newsymbol\varsupsetneq 2321
\newsymbol\nsubseteqq 2322
\newsymbol\nsupseteqq 2323
\newsymbol\subsetneqq 2324
\newsymbol\supsetneqq 2325
\newsymbol\varsubsetneqq 2326
\newsymbol\varsupsetneqq 2327
\newsymbol\subsetneq 2328
\newsymbol\supsetneq 2329
\newsymbol\nsubseteq 232A
\newsymbol\nsupseteq 232B
\newsymbol\nparallel 232C
\newsymbol\nmid 232D
\newsymbol\nshortmid 232E
\newsymbol\nshortparallel 232F
\newsymbol\nvdash 2330
\newsymbol\nVdash 2331
\newsymbol\nvDash 2332
\newsymbol\nVDash 2333
\newsymbol\ntrianglerighteq 2334
\newsymbol\ntrianglelefteq 2335
\newsymbol\ntriangleleft 2336
\newsymbol\ntriangleright 2337
\newsymbol\nleftarrow 2338
\newsymbol\nrightarrow 2339
\newsymbol\nLeftarrow 233A
\newsymbol\nRightarrow 233B
\newsymbol\nLeftrightarrow 233C
\newsymbol\nleftrightarrow 233D
\newsymbol\divideontimes 223E
\newsymbol\varnothing 203F
\newsymbol\nexists 2040
\newsymbol\Finv 2060
\newsymbol\Game 2061
\newsymbol\mho 2066
\newsymbol\eth 2067
\newsymbol\eqsim 2368
\newsymbol\beth 2069
\newsymbol\gimel 206A
\newsymbol\daleth 206B
\newsymbol\lessdot 236C
\newsymbol\gtrdot 236D
\newsymbol\ltimes 226E
\newsymbol\rtimes 226F
\newsymbol\shortmid 2370
\newsymbol\shortparallel 2371
\newsymbol\smallsetminus 2272
\newsymbol\thicksim 2373
\newsymbol\thickapprox 2374
\newsymbol\approxeq 2375
\newsymbol\succapprox 2376
\newsymbol\precapprox 2377
\newsymbol\curvearrowleft 2378
\newsymbol\curvearrowright 2379
\newsymbol\digamma 207A
\newsymbol\varkappa 207B
\newsymbol\Bbbk 207C
\newsymbol\hslash 207D
\undefine\hbar
\newsymbol\hbar 207E
\newsymbol\backepsilon 237F
%  Restore the catcode value for @ that was previously saved.
\catcode`\@=\csname pre amssym.tex at\endcsname
%---------------------------------------------------
\font\quatorzebf=cmbx10 at 14pt
\font\douzebf=cmbx10 at 12pt
\font\quatorzemi=cmmi10 at 14pt
\font\dixmi=cmmi10 at 10pt
\def\({{\rm (}}
\def\){{\rm )}}
%---------------------------------------------------
\def\maprightover#1{\smash{\mathop{\longrightarrow}
\limits^{#1}}}
\def\mapleftover#1{\smash{\mathop{\longleftarrow}
\limits^{#1}}}
\def\maprightunder#1{\smash{\mathop{\longrightarrow}
\limits_{#1}}}
\def\mapleftunder#1{\smash{\mathop{\longleftarrow}
\limits_{#1}}}
\def\mapdownleft#1{\llap{$\vcenter
{\hbox{$\scriptstyle#1$}}$}\Big\downarrow}
\def\mapdownright#1{\Big\downarrow
\rlap{$\vcenter{\hbox{$\scriptstyle#1$}}$}}
%------------------------------------------
\auteurcourant={G. LAUMON}
\titrecourant={FAISCEAUX AUTOMORPHES POUR $GL_n$}

%================================================
%================================================
%================================================

\centerline{\quatorzebf Faisceaux automorphes pour
$\hbox{\quatorzemi GL}_{\hbox{\dixmi n}}$:}
\vskip 2mm
\centerline{\quatorzebf la premi\`ere construction
de Drinfeld}
\vskip 15mm
\centerline{G\'erard Laumon}
\vskip 25mm

Drinfeld a associ\'e \`a tout syst\`eme
local irr\'eductible $L$ de rang $2$ sur une surface
de Riemann compacte $X$ un ``faisceau automorphe''
$W_{L,2}$ sur l'espace de modules ${\cal F}{\it ib}_X^2$
des fibr\'es vectoriels de rang $2$ sur $X$.

En fait, Drinfeld
a donn\'e deux constructions de $W_{L,2}$.
La premi\`ere est bas\'ee sur l'existence d'un mod\`ele
de Whittaker pour toute repr\'esentation automorphe
cuspidale pour $GL_2$. Elle garde un sens en caract\'eristique
positive et a \'et\'e expos\'ee dans [Dr].
Si $n$ est un entier $\geq 2$ arbitraire,
toute repr\'esentation automorphe
cuspidale pour $GL_n$ admet un mod\`ele de Whittaker
et il est possible d'\'etendre conjecturalement
cette premi\`ere construction en rang $n$ (cf. [La]).

La seconde construction est bas\'ee sur un
calcul explicite de l'anneau des sections globales
d'un certain faisceau d'op\'erateurs diff\'erentiels sur
l'espace de modules ${\cal F}{\it ib}_X^2$. Elle n'a pour l'instant de sens
qu'en caract\'eristique nulle mais, par contre, elle est de port\'ee
plus g\'en\'erale que la premi\`ere construction car elle s'\'etend
conjecturalement aux $G$-syst\`emes
locaux pour une groupe r\'eductif $G$
arbitraire.

Dans ces deux expos\'es, donn\'es \`a Luminy
lors du Colloque ``Fibr\'es vectoriels sur les
courbes et th\'eorie de Langlands'',
en juin 1995, je reprends la premi\`ere construction
de Drinfeld en rang $n$ arbitraire, d'un point de vue l\'eg\`erement
diff\'erent de celui de [La] (ce point de vue a aussi
\'et\'e consid\'er\'e par Bernstein et Vilonen).
\vskip 10mm

\centerline{\douzebf Table des mati\`eres}
\vskip 5mm

\noindent {\bf Expos\'e I}
\vskip 3mm

\noindent \line{1. Le diagramme de champs fondamental\dotfill
$\,\,\,3$}
\vskip 1mm

\noindent \line{2. Faisceaux de Whittaker\dotfill $\,\,\,5$}
\vskip 1mm

\noindent \line{3. La construction fondamentale\dotfill $\,\,\,7$}
\vskip 1mm

\noindent \line{4. Op\'erateurs de Hecke\dotfill $\,\,\,9$}
\vskip 3mm

\noindent {\bf Expos\'e II}
\vskip 3mm

\noindent \line{1. Passage \`a la caract\'eristique $p$\dotfill
$\,\,19$}
\vskip 1mm

\noindent \line{2. Fonction trace de Frobenius de $W_L$\dotfill
$\,\,20$}
\vskip 1mm

\noindent \line{3. La fonction trace de Frobenius de $W_{L,n}^\circ$
\dotfill $\,\,21$}
\vskip 1mm

\noindent \line{4. Op\'erateurs de Hecke et fonctions trace de
Frobenius\dotfill $\,\,28$}
\vskip 3mm

\noindent {\bf Compl\'ements}
\vskip 3mm

\noindent \line{1. Variante sch\'ematique du diagramme
fondamental\dotfill $\,\,36$}
\vskip 1mm

\noindent \line{2. Vari\'et\'es caract\'eristiques\dotfill $\,\,39$}
\vskip 1mm

\noindent \line{3. Une variante de la construction du complexe
$W_{L,n}^\circ$\dotfill $\,\,48$}
\vskip 3mm

\noindent {\bf Bibliographie}

\vfill\eject
\centerline{\douzebf Expos\'e I}
\vskip 20mm
{\bf 1. Le diagramme de champs fondamental.}
\vskip 5mm
On fixe une courbe $X$, connexe, projective, lisse et de genre
$g\geq 2$ sur le corps des nombres complexes.
On fixe de plus un entier $n\geq 1$ et un entier $m_0>
n(3n-1)(g-1)$.

Pour chaque $i=1,\ldots ,n$, on pose
$$
m_i=m_0+i(i-1)(g-1)
$$
et on note ${\cal C}{\it oh}_{X}^{i,m_i}$ le champ des
${\cal O}_X$-Modules coh\'erents ${\cal M}_i$ de rang
g\'en\'erique $i$ et de degr\'e $m_i$. C'est un champ
alg\'ebrique connexe et lisse de dimension $i^2(g-1)$ sur ${\bb C}$.

On note
$$
{\cal C}_i\subset {\cal C}{\it oh}_{X}^{i,m_i}
$$
l'ouvert form\'e des ${\cal M}_i$ tels que
$$
{\rm Hom}_{{\cal O}_X}({\cal M}_i,(\Omega_X^1)^{\otimes (2n-1)})
=(0).
$$
Pour chaque ${\cal M}_i$ dans ${\cal C}_i$ on a
$$
{\rm Hom}_{{\cal O}_X}({\cal M}_i,(\Omega_X^1)^{\otimes j})=(0)
\qquad (\forall j\leq 2n-1),
$$
puisque $H^0(X,\Omega_X^1)\not= (0)$ par hypoth\`ese sur $g$, et
donc
$$
{\rm Ext}_{{\cal O}_X}^1((\Omega_X^1)^{\otimes j},{\cal M}_i)=(0)
\qquad (\forall j\leq 2n-2)
$$
par dualit\'e de Serre.

\th LEMME 1.1
\enonce
L'ouvert ${\cal C}_i$ est non vide et de type fini sur ${\bb C}$.

Plus pr\'ecis\'ement, pour $i=0$,
on a ${\cal C}_0={\cal C}{\it oh}_{X}^{0,m_0}$.
Pour $i=1,\ldots ,n$,
l'ouvert ${\cal C}_i$ de ${\cal C}{\it oh}_{X}^{i,m_i}$ contient l'ouvert
${\cal F}{\it ib}_{X}^{i,m_i,{\rm ss}}$ des ${\cal O}_X$-Modules
localement libres semi-stables
de rang $i$ et de degr\'e $m_i$ et, pour tout point ${\cal M}_i$
de ${\cal C}_i$, on a
$$
{m_i-{\rm deg}\,{\cal M}_i^{\rm tors}\over i}\geq
\mu^-({\cal M}_i/{\cal M}_i^{\rm tors})\geq 4(n-1)(g-1),
$$
o\`u ${\cal M}_i^{\rm tors}$ est le plus grand
sous-${\cal O}_X$-Module de
torsion de ${\cal M}_i$ et o\`u, pour tout ${\cal O}_X$-Module
localement libre de rang fini, on a pos\'e
$$
\mu^-({\cal E})={\rm Inf}\big\{{{\rm deg}({\cal E}/{\cal F})\over
{\rm rg}({\cal E}/{\cal F})}\mid (0)\subset {\cal F}\subsetneqq
{\cal E}\hbox{ et }{\cal E}/{\cal F}\hbox{ est sans torsion}\big\}.
$$
\endth

\rem Preuve
\endrem
Pour $i=0$, chaque ${\cal M}_0$ est engendr\'e par ses sections
globales, de sorte que ${\cal C}_0={\cal C}{\it oh}_{X}^{0,m_0}$ admet une
pr\'esentation globale
$$
{\rm Quot}_{{\cal O}_X^{m_0}/X}^{0,m_0}\twoheadrightarrow
{\cal C}{\it oh}_{X}^{0,m_0}
$$
et est donc de type fini.

Supposons donc $i\geq 1$.
Si ${\cal L}_i$ un point de ${\cal F}{\it ib}_{X}^{i,m_i,{\rm ss}}$,
on a ${\rm Hom}_{{\cal O}_X}({\cal L}_i,
(\Omega_X^1)^{\otimes (2n-1)})=(0)$ car $(2n-1)(2g-2)<m_i/i$
par hypoth\`ese sur $m_0$, d'o\`u l'inclusion
${\cal F}{\it ib}_{X}^{i,m_i,{\rm ss}}\subset {\cal C}_i$.
De plus, pour tout ${\cal O}_X$-Module localement libre de rang
fini ${\cal G}$ qui est semi-stable, on a
$$
{\rm Hom}_{{\cal O}_X}({\cal G},(\Omega_X^1)^{\otimes (2n-1)})
\not=(0)
$$
d\`es que
$$
\mu ({\cal G})={{\rm deg}({\cal G})\over {\rm rg}({\cal G})}
<(2n-2)(2g-2)=4(n-1)(g-1).
$$
Par suite, si un tel ${\cal G}$ est quotient d'un ${\cal M}_i\in
{\cal C}_i$, on a n\'ecessairement $\mu ({\cal G})\geq 4(n-1)(g-1)$.
D'o\`u la conclusion

\hfill\hfill\cqfd

Pour chaque $i=0,\ldots ,n$, notons
$$
\pi_i:{\cal E}_i\rightarrow {\cal C}_i
$$
le fibr\'e vectoriel
de fibre ${\rm Hom}_{{\cal O}_X}((\Omega_X^1)^{\otimes (i-1)},
{\cal M}_i)$ en ${\cal M}_i\in {\cal C}_i$ et
$$
\pi_i^\vee :{\cal E}_i^\vee\rightarrow {\cal C}_i
$$
le fibr\'e vectoriel dual, de fibre
${\rm Ext}_{{\cal O}_X}^1({\cal M}_i,(\Omega_X^1)^{\otimes i})$
en ${\cal M}_i$. Les fibr\'es $\pi_i$ et $\pi_i^\vee$ sont
tous les deux  de rang $m_0-i^2(g-1)$
et ${\cal E}_i$ et ${\cal E}_i^\vee$ sont donc tous les
deux des champs alg\'ebriques lisses
et connexes, de dimension $m_0$, sur ${\bb C}$. Un point de
${\cal E}_i$ (resp. ${\cal E}_i^\vee$)
est une fl\`eche
$$
(\Omega_X^1)^{\otimes (i-1)}\rightarrow {\cal M}_i
$$
(resp. une extension
$$
(\Omega_X^1)^{\otimes i}\hookrightarrow {\cal M}_{i+1}
\twoheadrightarrow  {\cal M}_i\,)
$$
de ${\cal O}_X$-Modules avec ${\cal M}_i\in {\cal C}_i$.
\vskip 2mm

Pour chaque $i=1,\ldots ,n$, on notera
$$
{\cal E}_i^\circ\subset {\cal E}_i
$$
l'ouvert des fl\`eches de ${\cal O}_X$-Modules
$$
(\Omega_X^1)^{\otimes (i-1)}\rightarrow {\cal M}_i
$$
qui sont injectives (cet ouvert est non vide d'apr\`es le lemme 1.1)
et on notera
$$
{\cal E}_{i-1}^{\vee\circ}\subset {\cal E}_{i-1}^\vee
$$
l'ouvert des extensions de ${\cal O}_X$-Modules
$$
(\Omega_X^1)^{\otimes (i-1)}\hookrightarrow {\cal M}_i
\twoheadrightarrow  {\cal M}_{i-1}
$$
telles que l'on ait ${\rm Hom}_{{\cal O}_X}({\cal M}_i,
(\Omega_X^1)^{\otimes (2n-1)})=(0)$.
On d\'efinit un isomorphisme de ${\bb C}$-champs alg\'ebriques
$$
{\cal E}_i^\circ\buildrel\sim\over\longrightarrow
{\cal E}_{i-1}^{\vee\circ}
$$
en envoyant
$$
(\Omega_X^1)^{\otimes (i-1)}\hookrightarrow {\cal M}_i
$$
sur
$$
(\Omega_X^1)^{\otimes (i-1)}\hookrightarrow {\cal M}_i
\twoheadrightarrow  {\cal M}_{i-1},
$$
o\`u ${\cal M}_{i-1}$ est le quotient de ${\cal M}_i$ par l'image de
$(\Omega_X^1)^{\otimes (i-1)}$. On a bien
$$
{\rm Hom}_{{\cal O}_X}({\cal M}_{i-1},
(\Omega_X^1)^{\otimes (2n-1)})=(0)
$$
car on a un plongement
$$
{\rm Hom}_{{\cal O}_X}({\cal M}_{i-1},
(\Omega_X^1)^{\otimes (2n-1)})\hookrightarrow
{\rm Hom}_{{\cal O}_X}({\cal M}_i,(\Omega_X^1)^{\otimes (2n-1)})
=(0).
$$
\vskip 2mm

Le diagramme fondamental est alors le diagramme de morphismes
de champs alg\'ebriques sur ${\bb C}$
$$
\matrix{&&\kern -3mm{\cal E}_n\supset {\cal E}_n^\circ\cong
{\cal E}_{n-1}^{\vee\circ}\subset
{\cal E}_{n-1}^\vee\kern -3mm &\cr
&\swarrow&&\searrow \cr
{\cal C}_n\kern -3mm&&&\cr}\quad\cdots\quad
\matrix{&&&\kern -3mm{\cal E}_1\supset {\cal E}_1^\circ\cong
{\cal E}_0^{\vee\circ}
\subset {\cal E}_0^\vee\kern -3mm&&\cr
\searrow &&\swarrow &&\searrow &\cr
&\kern -3mm{\cal C}_1\kern -3mm &&&&\kern -3mm{\cal C}_0.\cr}
$$

On a des actions naturelles du groupe
multiplicatif ${\bb G}_m$ (sur ${\bb C}$) sur les fibr\'es vectoriels
du diagramme ci-dessus, par multiplication dans les fibres.
Il est facile de voir que, pour tout $i=1,\ldots,n$, les ouverts
${\cal E}_i^\circ\subset {\cal E}_i$ et ${\cal E}_{i-1}^{\vee\circ}
\subset {\cal E}_{i-1}^\vee$ sont ${\bb G}_m$-invariants et
l'isomorphisme ${\cal E}_i^\circ\cong {\cal E}_{i-1}^{\vee\circ}$
est ${\bb G}_m$-\'equivariant.
\vskip 5mm

{\bf 2. Faisceaux de Whittaker.}
\vskip 5mm

Soit $L$ un syst\`eme local
irr\'eductible de ${\bb C}$-espaces vectoriels de rang $n$
sur $X$. Une construction, due \`a Deligne,
associe \`a $L$ un faisceau constructible de
${\bb C}$-espaces vectoriels
$L^{(m_0)}$ sur le produit sym\'etrique
$X^{(m_0)}=X^{m_0}/{\goth S}_{m_0}$, espace de modules
des diviseurs effectifs de degr\'e $m_0$ sur $X$: on consid\`ere
le rev\^etement
$$
r:X^{m_0}\rightarrow X^{(m_0)},~(x_1,\ldots ,x_{m_0})\mapsto
x_1+\cdots +x_{m_0},
$$
fini \'etale galoisien (ramifi\'e d\`es
que $m_0\geq 2$), de groupe de Galois ${\goth S}_{m_0}$,
et on pose
$$
L^{(m_0)}=\bigl(r_*(L^{\boxtimes m_0})\bigr)^{{\goth S}_{m_0}}.
$$
On v\'erifie facilement que $L^{(m_0)}[m_0]$ est un
faisceau pervers irr\'eductible sur le
${\bb C}$-sch\'ema lisse et
connexe, de dimension $m_0$, $X^{(m_0)}$.

Dans [La] \S 3 nous avons prolong\'e cette construction
en associant \`a $L$ un faisceau pervers irr\'eductible $W_L$,
dit de Whittaker, sur le champ ${\cal C}_0=
{\cal C}{\it oh}_{X}^{0,m_0}$, tel que l'on ait
$$
\iota_{(m_0)}^*W_L[m_0]\cong L^{(m_0)}[m_0],
$$
o\`u $\iota_{(m_0)}:X^{(m_0)}\rightarrow
{\cal C}_0$ est le morphisme lisse,
purement de dimension relative $m_0$, qui envoie
$D$ sur ${\cal O}_D$. Rappelons cette construction.

On note $\widetilde{\cal C}_0=\widetilde{{\cal C}oh}_{X}^{0,m_0}$
le champ alg\' ebrique sur ${\bb C}$ param\`etrant les drapeaux
$$
{\cal M}_0^{\bullet}=\bigl({\cal M}_0^{m_0}\twoheadrightarrow
{\cal M}_0^{m_0-1}\twoheadrightarrow\cdots\twoheadrightarrow
{\cal M}_0^{1}\twoheadrightarrow{\cal M}_0^{0}=(0)\bigr),
$$
o\`u chaque ${\cal M}_0^i$ est un ${\cal O}_X$-Module coh\'erent de
rang g\'en\'erique $0$ et de degr\'e $i$. On a des projections
naturelles
$$
\def\normalbaselines{\baselineskip20pt\lineskip3pt
\lineskiplimit3pt}
\matrix{\widetilde{\cal C}_0&\maprightover{\rho}&X^{m_0}\cr
\mapdownright{\pi}&&\cr
{\cal C}_0&&\cr}
$$
avec
$$
\pi ({\cal M}_0^{\bullet})={\cal M}_0^{m_0}
$$
et
$$
\rho ({\cal M}_0^{\bullet})=(x_1,\ldots ,x_{m_0})
$$
o\`u $\{x_i\}$ est le support de ${\rm Ker}({\cal M}_0^i
\twoheadrightarrow{\cal M}_0^{i-1})$.
Le morphisme $\pi$ est repr\'esentable et projectif
et on a un carr\'e cart\'esien
$$
\def\normalbaselines{\baselineskip20pt\lineskip3pt
\lineskiplimit3pt}
\matrix{X^{m_0}&\maprightover{\iota_{m_0}}&
\widetilde{\cal C}_0\cr
\mapdownleft{r}&\square &\mapdownright{\pi}\cr
X^{(m_0)}&\maprightunder{\iota_{(m_0)}}&{\cal C}_0\cr}
$$
avec
$$
\rho\circ \iota_{m_0}={\rm id}_{X^{m_0}}.
$$

\th TH\'EOR\`EME 2.1 ([La] Thm. (3.3.1))
\enonce
Le complexe $R\pi_*\rho^*(L^{\boxtimes m_0})$ est en fait un
faisceau pervers, plac\'e en degr\'e $0$,
sur le champ alg\'ebrique ${\cal C}_0$
\(lisse purement de dimension $0$ sur ${\bb C}$\).

De plus, ce faisceau pervers est le prolongement
interm\'ediaire de sa restriction \`a l'ouvert dense
$\iota_{(m_0)}(X^{(m_0)})\subset {\cal C}_0$.
\endth

On a
$$
\iota_{(m_0)}^*R\pi_*\rho^*(L^{\boxtimes m_0})\cong
r_*(L^{\boxtimes m_0})
$$
et donc, le groupe sym\'etrique ${\goth S}_{m_0}$ agit sur
$\iota_{(m_0)}^*R\pi_*\rho^*(L^{\boxtimes m_0})$.
D'apr\`es le th\'eor\`eme 2.1, cette action se prolonge
en une action de ${\goth S}_{m_0}$ sur
$R\pi_*\rho^*(L^{\boxtimes m_0})$.
Alors, par d\'efinition, $W_L$ est le faisceau pervers
$$
\bigl(R\pi_*\rho^*(L^{\boxtimes m_0})\bigr)^{{\goth S}_{m_0}}.
\leqno (2.2)
$$

Il r\'esulte du th\'eor\`eme 2.1 que $W_L$ est aussi le
prolongement interm\'ediaire de sa restriction \`a l'ouvert dense
$\iota_{(m_0)}(X^{(m_0)})\subset {\cal C}_0$.
En particulier, $W_L$ est un faisceau
pervers irr\'eductible (cf. [La] Thm. (3.3.10)).
\vskip 5mm

{\bf 3. La construction fondamentale.}
\vskip 5mm

Rappelons qu'un faisceau constructible de
${\bb C}$-espaces vectoriels
sur un fibr\'e vectoriel est dit monodromique (au sens de Verdier)
si sa restriction \`a chaque droite \'epoint\'ee (i.e. \`a chaque
orbite de ${\bb G}_m$ dans le compl\'ementaire de la section nulle)
est localement constante.

Pour tout $i=1,\ldots ,n-1$,
notons $D_{\rm mon}^{\rm b}({\cal E}_i,{\bb C})$ et
$D_{\rm mon}^{\rm b}({\cal E}_i^\vee,{\bb C})$
les sous-cat\'egories pleines de $D_{\rm c}^{\rm b}({\cal E}_i,
{\bb C})$ et $D_{\rm c}^{\rm b}({\cal E}_i^\vee,{\bb C})$
respectivement dont les objets sont les complexes
de faisceaux de ${\bb C}$-espaces vectoriels \`a
cohomologie constructible et monodromique.
On dispose alors de la transformation de Fourier
$$
{\cal F}_{i} :D_{\rm mon}^{\rm b}({\cal E}_i,{\bb C})\rightarrow
D_{\rm mon}^{\rm b}({\cal E}_i^\vee ,{\bb C})
$$
pour le fibr\'e ${\cal E}_i\rightarrow {\cal C}_i$. C'est une
\'equivalence de
cat\'egories qui commute \`a la dualit\'e et qui transforme
les faisceaux pervers monodromiques sur ${\cal E}_i$ en des
faisceaux pervers
monodromiques sur ${\cal E}_i^\vee$ (cf. [Br] Ch. II, \S VI).

Si $K_i\in {\rm ob}\,D_{\rm c}^{\rm b}({\cal E}_i,{\bb C})$ est
muni d'un isomorphisme
$$
\mu_i^*K_i\cong {\bb C}\boxtimes K_i
$$
dans $D_{\rm c}^{\rm b}({\bb G}_m\times {\cal E}_i,{\bb C})$,
o\`u $\mu_i:{\bb G}_m\times {\cal E}_i\rightarrow {\cal E}_i$
est l'action naturelle
de ${\bb G}_m$ par multiplication dans les fibres de ${\cal E}_i$,
$K_i$ est trivialement monodromique et ${\cal F}_i(K_i)$ est lui
aussi
muni d'un isomorphisme
$$
\mu_i^{\vee *}{\cal F}_i(K_i)\cong {\bb C}\boxtimes
{\cal F}_i(K_i)
$$
dans $D_{\rm c}^{\rm b}({\bb G}_m\times {\cal E}_i^\vee ,{\bb C})$,
o\`u $\mu_i^\vee :{\bb G}_m\times {\cal E}_i^\vee
\rightarrow {\cal E}_i^\vee$ est
l'action naturelle de ${\bb G}_m$ par multiplication dans les fibres
de ${\cal E}_i^\vee$.
\vskip 2mm

Par r\'ecurrence
sur $i$, on d\'efinit un faisceau pervers irr\'eductible $W_{L,i}$
sur ${\cal E}_i$, muni d'un isomorphisme
$$
\mu_i^*W_{L,i}\cong {\bb C}\boxtimes W_{L,i}
$$
dans $D_{\rm c}^{\rm b}({\bb G}_m\times {\cal E}_i,{\bb C})$,
en posant
$$
W_{L,1}=j_{1!*}j_0^{\vee *}\pi_0^{\vee *}W_L[m_0]
$$
et
$$
W_{L,i}=j_{i!*}j_{i-1}^{\vee *}{\cal F}_{i-1} (W_{L,i-1})
$$
pour $i=2,\ldots ,n$, o\`u $j_i:{\cal E}_i^\circ
\hookrightarrow {\cal E}_i$ et
$j_{i-1}^\vee :{\cal E}_i^\circ\cong {\cal E}_{i-1}^{\vee\circ}
\hookrightarrow {\cal E}_{i-1}^\vee$ sont les inclusions.

\th CONJECTURE 3.1
\enonce
Il existe un faisceau pervers irr\'eductible $K_L$ sur
${\cal C}_n$, unique \`a isomorphisme pr\`es, tel que
$$
W_{L,n}\cong\pi_n^*K_L[m_0-n^2(g-1)].
$$

\endth

On d\'efinit aussi par r\'ecurrence sur $i$ l'objet
$W_{L,i}^\circ$ de $D_{\rm c}^{\rm b}({\cal E}_i^\circ,{\bb C})$,
muni d'un isomorphisme
$$
\mu_i^{\circ *}W_{L,i}^\circ\cong {\bb C}\boxtimes W_{L,i}^\circ
$$
dans $D_{\rm c}^{\rm b}({\bb G}_m\times {\cal E}_i^\circ,{\bb C})$, o\`u
$\mu_i^\circ$ est la restriction \`a l'ouvert
${\cal E}_i^\circ\subset {\cal E}_i$ de l'action $\mu_i$,
en posant
$$
W_{L,1}^\circ =j_0^{\vee *}\pi_0^{\vee *}W_L[m_0]
$$
et
$$
W_{L,i}^\circ =j_{i-1}^{\vee *}{\cal F}_{i-1} (j_{i-1,!}W_{L,i-1}^\circ)
$$
pour $i=2,\ldots ,n$. On peut alors renforcer la conjecture 3.1 en
demandant de plus~:

\th CONJECTURE 3.2
\enonce
Pour chaque $i=2,\ldots ,n$, $W_{L,i}^\circ$ est pervers
irreductible
{\rm (}pour $i=1$ c'est automatique{\rm )} et
on a des isomorphismes
$$
j_{i-1,!}W_{L,i-1}^\circ\cong j_{i-1,!*}W_{L,i-1}^\circ\cong
W_{L,i-1}.
$$
\endth

Pour tout entier $\ell$, soit ${\rm Pic}_X^{\ell}$ la
composante connexe du sch\'ema de Picard
de $X$ qui param\`etre les
classes d'isomorphie de fibr\'es en droites de degr\'e $\ell$ sur
$X$ (${\rm Pic}_X^{\ell}$ est une vari\'et\'e connexe,
projective et lisse de dimension $g$).
Le quotient ab\'elien maximal du groupe fondamental de $X$ et
le groupe fondamental de ${\rm Pic}_X^{\ell}$ sont
canoniquement isomorphes (ces deux groupes sont canoniquement
isomorphes \`a $H_1(X,{\bb Z})$).
Le syt\`eme local $\bigwedge^nL$ de rang $1$ sur $X$ induit donc un
syst\`eme de rang $1$ sur ${\rm Pic}_X^{\ell}$ que l'on notera
$\langle L\rangle$.

Soient $m_0'$ un entier $\geq m_0+g$, $m_n'=m_0'+n(n-1)(g-1)$
et ${\cal C}_n'$ l'ouvert du champ ${\cal C}{\it oh}_X^{n,m_n'}$
form\'e des ${\cal M}_n'$ tels que
${\rm Hom}_{{\cal O}_X}({\cal M}_n',
(\Omega_X^1)^{\otimes (2n-1)})=(0)$. On a un morphisme de champs
$$
\mu_n :{\rm Pic}_X^{m_n'-m_n}\times
{\cal C}_n\rightarrow {\cal C}_n',~
({\cal L},{\cal M}_n)\mapsto
{\cal L}\otimes_{{\cal O}_X}{\cal M}_n
$$
qui est lisse, \`a fibres connexes de dimension $g$ (tout
${\cal L}\in {\rm Pic}_X^{m_n'-m_n}$ admet des
sections globales non nulles puisque $m_n'-m_n\geq g$ et donc
${\cal L}\otimes_{{\cal O}_X}{\cal M}_n$ est bien dans
${\cal C}_n'$ pour tout ${\cal M}_n\in {\cal C}_n$).

\th CONJECTURE 3.3
\enonce
Si $K_L$ et $K_L'$ sont les objets de $D_{\rm c}^{\rm b}({\cal C}_n,
{\bb C})$ et $D_{\rm c}^{\rm b}({\cal C}_n',{\bb C})$ respectivement
qui sont conjecturalement associ\'es en $3.1$ au syst\`eme local
irr\'eductible $L$ de rang $n$ sur $X$, on a
$$
\mu_n^*K_L'\cong \langle L\rangle\boxtimes K_L.
$$
\endth

Supposons la conjecture 3.1 v\'erifi\'ee
pour tout entier $m_0>n(3n-1)(g-1)$.
Alors, la conjecture 3.3 permet de calculer le
prolongement interm\'ediaire de $K_L$ \`a
${\cal C}{\it oh}_X^{n,m_n}$ tout entier
puisque, pour tout ${\cal O}_X$-Module coh\'erent ${\cal M}_n$
de rang g\'en\'erique $n$ sur $X$, il existe $\ell\geq g$ et
${\cal L}$ dans ${\rm Pic}_X^{\ell}$ tels que
$$
{\rm Hom}_{{\cal O}_X}({\cal L}\otimes_{{\cal O}_X}{\cal M}_n,
(\Omega_X^1)^{\otimes (2n-1)})=(0).
$$
De plus, si on note encore $K_L$ ce prolongement interm\'ediaire,
cette conjecture permet d'\'etendre le faisceau pervers
$K_L$ sur
$\coprod_{m_n>2n(2n-1)(g-1)}{\cal C}{\it oh}_X^{n,m_n}$
au champ ${\cal C}{\it oh}_X^n$ (des ${\cal O}_X$-Modules
coh\'erents de rang g\'en\'erique $n$ sur $X$) tout entier.
\vskip 2mm

{\pc REMARQUE} 3.4\pointir Pour tout entier $m_0\geq 0$,
on a un morphisme de champs
$$
{\rm det}_0:{\cal C}{\it oh}_X^{0,m_0}\rightarrow
{\rm Pic}_X^{m_0}
$$
qui envoie ${\cal O}_{D_1+D_2+\cdots +D_n}\oplus
{\cal O}_{D_2+\cdots +D_n}\oplus\cdots\oplus
{\cal O}_{D_n}$ sur ${\cal O}_X(D_1+2D_2+\cdots +nD_n)$.
Pour tout $i=1,\ldots ,n$ et tout entier $m_i$, on a
un morphisme de champs
$$
{\rm det}_i:{\cal C}{\it oh}_X^{i,m_i}\rightarrow
{\rm Pic}_X^{m_i}
$$
qui envoie ${\cal L}_i\oplus {\cal T}_i$, o\`u ${\cal L}_i$ est
localement libre de rang $i$ et ${\cal T}_i$ est de torsion, sur
$\bigr(\bigwedge^i{\cal L}_i\bigl)
\otimes_{{\cal O}_X}{\rm det}_0({\cal T}_i)$.
Il est clair que
$$
{\rm det}_1\circ\pi_1\circ j_1={\rm det}_0\circ\pi_0^\vee\circ
j_0^\vee :{\cal E}_1^\circ\rightarrow {\rm Pic}_X^{m_0}
$$
(on a $m_0=m_1$) et que, pour $L$ de rang $1$, on a
$W_L={\rm det}_0^*\langle L\rangle$.

Il s'ensuit que, pour $L$ de rang $1$, on a simplement
$$
K_L={\rm det}_1^*\langle L\rangle .
$$
\vskip 5mm

{\bf 4. Op\'erateurs de Hecke.}
\vskip 5mm

Pour chaque entier $i=1,\ldots ,n$ et chaque entier
$m_i'> m_i$, consid\'erons le champ
$$
{\cal H}{\it ecke}_{X}^{i,m_i,m_i'}
$$
des triplets $({\cal M}_i,{\cal M}_i',\alpha_i)$, o\`u
${\cal M}_i\in {\rm ob}\,{\cal C}{\it oh}_{X}^{i,m_i}$,
${\cal M}_i'\in {\rm ob}\,{\cal C}{\it oh}_{X}^{i,m_i'}$ et
$\alpha_i :{\cal M}_i\hookrightarrow {\cal M}_i'$ est
un monomorphisme de ${\cal O}_X$-Modules. Alors, si on pose
$n_0=m_i'-m_i$, ${\rm Coker}(\alpha_i)$ est un point du champ
${\cal C}{\it oh}_X^{0,n_0}$.

On a des morphismes de champs
$$
p_i:{\cal H}{\it ecke}_{X}^{i,m_i,m_i'}\rightarrow
{\cal C}{\it oh}_X^{0,n_0}
\times {\cal C}{\it oh}_{X}^{i,m_i},~
({\cal M}_i,{\cal M}_i',\alpha_i)\mapsto
({\rm Coker}(\alpha_i),{\cal M}_i),
$$
et
$$
q_i:{\cal H}{\it ecke}_{X}^{i,m_i,m_i'}\rightarrow
{\cal C}{\it oh}_X^{i,m_i'},~
({\cal M}_i,{\cal M}_i',\alpha_i)\mapsto {\cal M}_i'.
$$
Le morphisme $q_i$ est repr\'esentable~: sa fibre en ${\cal M}_i'$
est le sch\'ema ${\rm Quot}_{{\cal M}_i'/X}^{0,n_0}$
des quotients ${\cal M}_i'\twoheadrightarrow {\cal N}_0$
avec ${\cal N}_0$ de rang g\'en\'erique $0$ et de degr\'e $n_0$.
Par contre, le morphisme $p_i$ n'est pas
repr\'esentable en g\'en\'eral~: sa fibre en
$({\cal N}_0,{\cal M}_i)$ est le champ de Picard associ\'e
au complexe $R{\rm Hom}_{{\cal O}_X}({\cal N}_0,{\cal M}_i)[1]$
(l'ensemble des classes d'isomorphie d'objets est en bijection avec
${\rm Ext}_{{\cal O}_X}^1({\cal N}_0,{\cal M}_i)$ et le
groupe des automorphismes de chaque objet
est isomorphe au groupe additif
${\rm Hom}_{{\cal O}_X}({\cal N}_0,{\cal M}_i)$).
Cependant $p_i$ est lisse, purement de dimension relative
$in_0$.

Dans la section 1 on a introduit l'ouvert ${\cal C}_i$ de
${\cal C}{\it oh}_X^{i,m_i}$. Introduisons de m\^eme l'ouvert
${\cal C}'_i$ de ${\cal C}{\it oh}_{X}^{i,m_i'}$  form\'e des
${\cal M}_i'$ qui v\'erifient la condition
$$
{\rm Hom}_{{\cal O}_X}({\cal M}_i',
(\Omega_X^1)^{\otimes (2n-1)})=(0).
$$
Il est clair que
$$
p_i^{-1}({\cal C}{\it oh}_X^{0,n_0}\times
{\cal C}_i)\subset q_i^{-1}({\cal C}'_i).
$$
On notera simplement ${\cal H}_i$ l'ouvert
$p_i^{-1}({\cal C}{\it oh}_X^{0,n_0}\times {\cal C}_i)$
de ${\cal H}{\it ecke}_{X}^{i,m_i,m_i'}$
et encore
$$
\matrix{&&{\cal H}_i&&\cr
\noalign{\smallskip}
&{}^{p_i}\kern -1.5mm\swarrow&&\searrow^{q_i}\cr
\noalign{\smallskip}
{\cal C}{\it oh}_X^{0,n_0}\times {\cal C}_i\kern -3mm&&&&
\kern -3mm{\cal C}_i'\cr}
$$
les restrictions \`a cet ouvert des projections $p_i$ et $q_i$
d\'efinies ci-dessus.

On d\'efinit l'op\'erateur de Hecke
$$
T_i:D_{\rm c}^{\rm b}({\cal C}_i',{\bb C})\rightarrow
D_{\rm c}^{\rm b}({\cal C}{\it oh}_X^{0,n_0}\times {\cal C}_i,{\bb C})
$$
par
$$
T_i(K')=Rp_{i!}q_i^*K'[in_0].
$$

\th TH\'EOR\`EME 4.1
\enonce
Si on note $V_L$, $W_L$ et $W_L'$
les faisceaux  pervers irr\'eductibles de Whittaker associ\'es
\`a $L$ sur ${\cal C}{\it oh}_X^{0,n_0}$,
${\cal C}_0$ et ${\cal C}_0'$ respectivement, on a
$$
T_0(W'_L)\cong V_L\boxtimes W_L
$$
dans $D_{\rm c}^{\rm b}({\cal C}{\it oh}_X^{0,n_0}\times{\cal C}_0,{\bb C})$.
\endth

\rem Preuve
\endrem
Notons ${\cal U}'$ l'ouvert de ${\cal C}_0'$ form\'e des
${\cal M}_0'$ qui sont isomorphes \`a
${\cal O}_{x_1'+\cdots +x_{m_0'}'}$ pour des points
$x_1',\ldots ,x_{m_0'}'$ de $X$ deux \`a deux distincts et
${\cal V}$ l'ouvert de ${\cal C}{\it oh}_X^{0,n_0}\times{\cal C}_0$
form\'e des couples $({\cal N}_0,{\cal M}_0)$ tels que ${\cal N}_0$
soit isomorphe \`a ${\cal O}_{y_1+\cdots +y_{n_0}}$ et
${\cal M}_0$ soit isomorphe \`a ${\cal O}_{x_1+\cdots +x_{m_0}}$
pour des points $y_1,\ldots ,y_{n_0},
x_1,\ldots ,x_{m_0}$ de $X$ deux \`a deux distincts.

Consid\'erons dans un premier temps
la restriction de $T_0(W_L')$ \`a l'ouvert ${\cal V}$ de
${\cal C}{\it oh}_X^{0,n_0}\times{\cal C}_0$.
Il est clair que $p_0$ est un isomorphisme au dessus de cet ouvert,
que
$$
p_0^{-1}({\cal V})=q_0^{-1}({\cal U}')
$$
et que $q_0$  est un rev\^etement fini \'etale de degr\'e
$\bigl({m_0'\atop n_0}\bigr)$ (i.e. le nombre des
sous-ensembles $\{y_1,\ldots ,y_{n_0}\}$ \`a $n_0$ \'el\'ements
de $\{x_1',\ldots ,x_{m_0'}'\}$). Or la fibre de $V_L\boxtimes
W_L$ en $({\cal O}_{y_1+\cdots +y_{n_0}},
{\cal O}_{x_1+\cdots +x_{m_0}})\in {\cal V}$
est $L_{y_1}\otimes\cdots\otimes L_{y_{n_0}}
\otimes L_{x_1}\otimes\cdots\otimes L_{x_{m_0}}$ et la fibre de
$W_L'$ en ${\cal O}_{x_1'+\cdots +x_{m_0'}'}\in {\cal U}'$ est
$L_{x_1'}\otimes\cdots\otimes L_{x_{m_0'}'}$. Par suite on a
$$
T_0(W_L')|{\cal V}\cong (V_L\boxtimes W_L)|{\cal V}.
$$

Pour conclure, il nous reste \`a v\'erifier que $T_0(W_L')$
est pervers et est le prolongement interm\'ediaire de
sa restriction \`a l'ouvert ${\cal V}$
(en effet $V_L\boxtimes W_L$ est un faisceau pervers
irr\'eductible de support ${\cal C}{\it oh}_X^{0,n_0}\times {\cal C}_0$
tout entier). Or, d'apr\`es le th\'eor\`eme 2.1, $W_L'$ est facteur
direct du faisceau pervers $R\pi_*'\rho'^*L^{\boxtimes m_0'}$ o\`u
les morphismes de champs
$$
\pi':\widetilde{\cal C}_0'\rightarrow {\cal C}_0'
$$
et
$$
\rho':\widetilde{\cal C}_0'\rightarrow X^{m_0'}
$$
sont d\'efinis comme suit: $\widetilde{\cal C}_0'$
est le champ des suites d'\'epimorphismes de
${\cal O}_X$-Modules coh\'erents de rang g\'en\'erique $0$,
$$
{\cal M}_0'^\bullet =\bigl({\cal M}_0'^{m_0'}\twoheadrightarrow
{\cal M}_0'^{m_0'-1}\twoheadrightarrow
\cdots\twoheadrightarrow {\cal M}_0'^1
\twoheadrightarrow {\cal M}_0'^0=(0)\bigr),
$$
avec ${\rm deg}({\cal M}_0'^i)=i$, $\pi'$ envoie
${\cal M}_0'^\bullet$ sur ${\cal M}_0'^{m_0'}$ et $\rho'$
envoie ${\cal M}_0'^\bullet$ sur $(x_1',\ldots ,x_{m_0'}')$ o\`u
on a not\'e par $x_i'$ le support de
${\rm Ker}({\cal M}_0'^i\twoheadrightarrow {\cal M}_0'^{i-1})$.
Il suffit donc de v\'erifier que
$$
T_0(R\pi_*'\rho'^*L^{\boxtimes m_0'})
$$
est pervers et est le prolongement interm\'ediaire de
sa restriction \`a l'ouvert ${\cal V}$.

Cette derni\`ere question est locale pour la topologie \'etale sur
${\cal C}{\it oh}_X^{0,n_0}\times {\cal C}_0$ et donc pour la
topologie \'etale sur $X$. Par cons\'equent, on peut supposer
que $X$ est la droite affine ${\bb A}^1$ sur ${\bb C}$ et que $L$
est le faisceau constant de fibre ${\bb C}^n$ et on est ramen\'e
\`a v\'erifier que $T_0(R\pi_*'{\bb C})$ est pervers
et est le prolongement interm\'ediaire de sa restriction
\`a l'ouvert ${\cal V}$.

La correspondance de Hecke admet alors la pr\'esentation suivante.
Notons ${\goth g}$ (resp. ${\goth g}'$, resp. ${\goth h}$)
l'alg\`ebre le Lie de $G=GL_{m_0}$ (resp.
$G'=GL_{m_0'}$, resp. $H=GL_{n_0}$)
sur ${\bb C}$~; notons ${\goth p}'$
la sous-alg\`ebre
parabolique de ${\goth g}'$ des matrices
$$
\xi'=\pmatrix{A&u\cr 0&a\cr}\in {\goth g}'
$$
avec $A\in {\goth g}$, $a\in {\goth h}$ et
$u\in {\rm Mat}_{m_0\times n_0}$
et notons $P'$ le sous-groupe parabolique de $G'$ correspondant.
On a ${\cal C}{\it oh}_X^{0,n_0}=[{\goth h}/H]$,
${\cal C}_0=[{\goth g}/G]$, ${\cal C}_0'=[{\goth g}'/G']$,
${\cal H}_0=[{\goth p}'/P']$, o\`u toutes les actions sont les actions
adjointes, et le morphisme $p_0$ (resp. $q_0$) est induit
par la projection
$$
{\goth p}'\rightarrow {\goth h}\times {\goth g},~
\pmatrix{A&u\cr 0&a\cr}\mapsto (a,A),
$$
(resp. par l'inclusion ${\goth p}'\subset {\goth g}'$). On a donc un
diagramme commutatif
$$
\matrix{&&\kern -2mm{\goth p}'\kern -2mm&&\cr
\noalign{\vskip 0.1mm}\cr
&\swarrow&\kern -2mm\downarrow\kern -2mm &\searrow&\cr
\noalign{\vskip 0.1mm}\cr
{\goth h}\times {\goth g}\kern -2mm&\leftarrow&\kern -2mm
[{\goth p}'/N']\kern -2mm&&\kern -2mm {\goth g}'\cr
\noalign{\vskip 0.1mm}\cr
\downarrow\kern -2mm&\square&\kern -2mm\downarrow
\kern -2mm&&\kern -2mm\downarrow\cr
\noalign{\vskip 0.1mm}\cr
[{\goth h}/H]\times [{\goth g}/G]\kern -2mm&\leftarrow&
\kern -2mm [{\goth p}'/P']\kern -2mm&\rightarrow &
\kern -2mm [{\goth g}'/G']\cr}
$$
o\`u $N'$ est le radical unipotent de $P'$. Dans ce diagramme,
la fl\`eche verticale ${\goth p}'\rightarrow [{\goth p}'/N']$
est repr\'esentable
et est une fibration \`a fibres toutes isomorphes
\`a l'espace affine ${\bb A}^{m_0n_0}$,
alors que la fl\`eche horizontale
$[{\goth p}'/N']\rightarrow {\goth h}\times {\goth g}$ a pour
fibre en $(a,A)$ le champ de Picard associ\'e au complexe
$$
 {\rm Mat}_{m_0\times n_0}\rightarrow
{\rm Mat}_{m_0\times n_0},~v\mapsto Av-va.
$$

Soient $B'$ le sous-groupe de Borel standard des matrices
triangulaires sup\'erieures dans $G'$ et ${\goth b}'$ son
alg\`ebre de Lie sur ${\bb C}$. Si l'on pose
$$
\widetilde{\goth g}'=\{(\xi',g'B')\in {\goth g}'\times G'/B'
\mid \xi'\in {}^{g'}{\goth b}'\},
$$
on a ${\cal C}_0'=[\widetilde{\goth g}'/G']$ et
le morphisme $\pi'$ est induit par la projection
$$
\widetilde{\goth g}'\rightarrow {\goth g}',~
(\xi',g'B')\mapsto\xi'.
$$
Par suite, si l'on pose
$$
\widetilde{\goth p}'=\{(\xi',g'B')\in {\goth g}'\times G'/B'
\mid \xi'\in {\goth p}'\cap {}^{g'}{\goth b}'\}
$$
et si l'on note $\sigma$ le morphisme compos\'e
$$
\widetilde{\goth p}'\rightarrow {\goth p}'\rightarrow
[{\goth p}'/N']\rightarrow {\goth h}\times {\goth g},
$$
i.e. la projection
$$
\widetilde{\goth p}'\rightarrow {\goth h}\times {\goth g},~
(\xi'=\pmatrix{A&u\cr 0&a\cr},g'B')\mapsto (a,A),
$$
la restriction de $T_0(R\pi_*'{\bb C})$ \`a
${\goth h}\times {\goth g}$, d\'ecal\'ee de
$n_0^2+m_0^2$, n'est autre que $R\sigma_!{\bb C}[m_0'^2]$
(l'image directe \`a supports propres du faisceau constant
${\bb C}$ par la projection ${\goth p}'\rightarrow [{\goth p}'/N']$
n'est autre que ${\bb C}[-2n_0m_0]$).

Soit $D_{P',B'}\subset {\goth S}_{m_0'}$
le syst\`eme de repr\'esentants des classes
dans $({\goth S}_{m_0}\times {\goth S}_{n_0})\backslash
{\goth S}_{m_0'}$ form\'e des \'el\'ements de longueur
minimale dans leurs classes et, pour chaque $w'\in D_{P',B'}$,
notons $\dot w'\in G'$ la matrice de permutation correspondante.
On a
$$
G'=\coprod_{w'\in D_{P',B'}}P'\dot w'B'
$$
et donc
$$
\widetilde{\goth p}'=\coprod_{w'\in D_{P',B'}}
\widetilde{\goth p}_{w'}',
$$
o\`u chaque
$$
\widetilde{\goth p}_{w'}'=\{(\xi',g'B')\in {\goth g}'\times
P'\dot w'B'/B'\mid \xi'\in {\goth p}'\cap {}^{g'}{\goth b}'\}
$$
est une sous-vari\'et\'e localement
ferm\'ee de $\widetilde{\goth p}'$. On a trivialement
$$
\widetilde{\goth p}_{w'}'\cong\{(\xi',p'(P'\cap {}^{\dot w'}B'))
\in {\goth g}'\times P'/(P'\cap {}^{\dot w'}B')
\mid \xi'^{p'}\in {\goth p}'\cap {}^{w'}{\goth b}'\}
$$
et il est bien connu que
$$
P'\cap {}^{\dot w'}B'=(M'\cap B')(N'\cap {}^{\dot w'}B'),
$$
o\`u $P'=M'N'$ est la d\'ecomposition de Levi standard de $P'$.
Par suite, la restriction $\sigma_{w'}$ de la projection
$\widetilde{\goth p}'\buildrel\sigma\over\longrightarrow
{\goth h}\times {\goth g}$ \`a $\widetilde{\goth p}_{w'}'$
se factorise en
$$
\widetilde{\goth p}_{w'}'\buildrel\tilde\sigma_{w'}\over
\longrightarrow\widetilde{\goth m}'\rightarrow {\goth m}'\cong
{\goth h}\times {\goth g},
$$
o\`u ${\goth m}'$ est l'alg\`ebre de Lie de $M'$ sur ${\bb C}$ et
$$
\widetilde{\goth m}'=\{(\mu',m'(M'\cap B))\in
{\goth m}'\times M'/(M'\cap B')
\mid \mu'\in {}^{m'}({\goth m}'\cap {\goth b})\}.
$$

Si l'on note ${\goth n}'$ l'alg\`ebre de Lie de $N'$ sur ${\bb C}$,
de sorte que ${\goth p}'={\goth m}'\oplus {\goth n}'$,
la fibre de $\widetilde\sigma_{w'}$ en $(\mu',m'(M'\cap B))$
s'identifie canoniquement \`a
$$
\{(\nu',n'(N'\cap {}^{\dot w'}B'))\in
{\goth n'}\times N'/(N'\cap {}^{\dot w'}B')\mid
\mu'^{m'n'}-\mu'^{m'}+\nu'^{m'n'}\in
{\goth n}'\cap {}^{\dot w'}{\goth b}'\}
$$
et est donc non canoniquement isomorphe \`a l'espace affine
${\goth n}'$. En particulier, les $\widetilde{\goth p}_{w'}'$
sont des vari\'et\'es lisses, connexes et de dimension \'egale \`a
celle de ${\goth p}'$, puisque $\widetilde{\goth m}'$
est lisse, connexe et de dimension \'egale \`a celle de ${\goth m}'$.
De plus, pour chaque $w'$, $R\widetilde\sigma_{w'!}{\bb C}$
est isomorphe \`a ${\bb C}[-2m_0n_0]$
et $R\sigma_{w'!}{\bb C}[m_0'^2]$ est donc
isomorphe \`a l'image directe \`a supports
propres de ${\bb C}[m_0^2+n_0^2]$ par la projection canonique
$\widetilde{\goth m}'\rightarrow {\goth m}'$.

Or cette projection est petite au sens
de Goresky-MacPherson. Par suite chaque
$R\sigma_{w'!}{\bb C}[m_0'^2]$ est un faisceau pervers
semi-simple,
prolongement interm\'ediaire de sa restriction \`a un ouvert
dense de ${\goth m}'$.
Comme $R\sigma_!{\bb C}[m_0'^2]$ se d\'evisse
en les $R\sigma_{w'!}{\bb C}[m_0'^2]$ dans
$D_{\rm c}^{\rm b}({\goth m}',{\bb C})$,
c'est n\'ecessairement un faisceau
pervers, extension successive des faisceaux pervers semi-simples
$R\sigma_{w'!}{\bb C}[m_0'^2]$.
Mais, d'apr\`es Deligne, les extensions sont toutes scind\'ees
car chaque $R\sigma_{w'!}{\bb C}[m_0'^2]$ est ``pur'', de poids
ind\'ependant de $w'$, et le faisceau pervers
$R\sigma_!{\bb C}[m_0'^2]$ est semi-simple et
prolongement interm\'ediaire de sa restriction
\`a un ouvert dense de ${\goth m}'$.
C'est ce que l'on voulait d\'emontrer.

\hfill\hfill\cqfd

Soient $\pi_i':{\cal E}_i'\rightarrow {\cal C}_i'$ le fibr\'e vectoriel
de fibre ${\rm Hom}_{{\cal O}_X}((\Omega_X^1)^{\otimes (i-1)},
{\cal M}_i')$ en ${\cal M}_i'$ et
$\pi_i'^\vee :{\cal E}_i'^\vee\rightarrow {\cal C}_i'$
son fibr\'e dual. On notera ${\cal E}_i'^\circ$ l'ouvert de
${\cal E}_i'$ correspondant aux morphismes
injectifs de $(\Omega_X^1)^{\otimes (i-1)}$ dans ${\cal M}_i'$ et
${\cal E}_{i-1}'^{\vee\circ}$
l'ouvert de ${\cal E}_{i-1}'^\vee$
correspondant aux extensions $(\Omega_X^1)^{\otimes (i-1)}
\hookrightarrow{\cal M}_i'\twoheadrightarrow {\cal M}_{i-1}'$
avec ${\rm Hom}_{{\cal O}_X}({\cal M}_i',
(\Omega_X^1)^{\otimes (2n-1)})=(0)$.
On notera enfin $j_i':{\cal E}_i'^\circ\hookrightarrow {\cal E}_i'$
et $j_{i-1}'^\vee:{\cal E}_i'^\circ\cong
{\cal E}_{i-1}'^{\vee\circ}\hookrightarrow
{\cal E}_{i-1}'^\vee$
les inclusions.

On peut alors consid\'erer le diagramme
$$
\def\normalbaselines{\baselineskip20pt\lineskip3pt
\lineskiplimit3pt}
\matrix{{\cal C}{\it oh}_X^{0,n_0}\times {\cal E}_i&
\mapleftover{p_i}&{\cal J}_i&
\maprightover{q_i}&{\cal E}_i'\cr
\mapdownleft{{\rm id}\times\pi_i}&
\square&\Big\downarrow&+&
\mapdownright{\pi_i'}\cr
{\cal C}{\it oh}_X^{0,n_0}\times {\cal C}_i
&\mapleftunder{p_i}&{\cal H}_i&
\maprightunder{q_i}&{\cal C}_i'\cr}
$$
o\`u ${\cal J}_i$ est le champ des triplets
$(e_i,{\cal M}_i',\alpha_i)$ avec
$(\Omega_X^1)^{\otimes (i-1)}\buildrel{e_i}\over\longrightarrow
{\cal M}_i$ dans ${\cal E}_i$ et $({\cal M}_i,{\cal M}_i',\alpha_i)$
dans ${\cal H}_i$, o\`u le premier carr\'e
est le carr\'e cart\'esien \'evident
et o\`u $q_i$ envoie $(e_i,{\cal M}_i',\alpha_i)
\in {\rm ob}\,{\cal J}_i$ sur le morphisme compos\'e
$(\Omega_X^1)^{\otimes (i-1)}\maprightover{e_i}{\cal M}_i
\maprightover{\alpha_i}{\cal M}_i'$. Il est clair
que
$$
p_i^{-1}({\cal C}{\it oh}_X^{0,n_0}\times{\cal E}_i^\circ)=
q_i^{-1}({\cal E}_i'^\circ).
$$
On notera ${\cal J}_i^\circ$ cet ouvert de ${\cal J}_i$ et
$$
\matrix{&&{\cal J}_i^\circ&&\cr
\noalign{\smallskip}
&{}^{p_i^\circ}\kern -1.5mm\swarrow&&\searrow^{q_i^\circ}\cr
\noalign{\smallskip}
{\cal C}{\it oh}_X^{0,n_0}\times {\cal E}_i^\circ\kern -3mm&&&&
\kern -3mm{\cal E}_i'^\circ\cr}
$$
les restrictions \`a ${\cal J}_i^\circ$ des projections $p_i$ et
$q_i$ d\'efinies ci-dessus.

On a des actions naturelles de ${\bb G}_{\rm m}$ sur
${\cal C}{\it oh}_X^{0,n_0}\times {\cal E}_i$, ${\cal J}_i$ et
${\cal E}_i'$ qui respectent les ouverts ${\cal C}{\it oh}_X^{0,n_0}
\times {\cal E}_i^\circ$, ${\cal J}_i^\circ$ et ${\cal E}_i'^\circ$
et pour lesquelles les projections $p_i$ et $q_i$ sont
\'equivariantes.
\vskip 2mm

On d\'efinit l'op\'erateur de Hecke
$$
T_i:D_{\rm c}^{\rm b}({\cal E}_i',{\bb C})\rightarrow
D_{\rm c}^{\rm b}({\cal C}{\it oh}_X^{0,n_0}\times {\cal E}_i,{\bb C})
$$
par
$$
T_i(K')=Rp_{i!}q_i^*K'[in_0]
$$
et l'op\'erateur de Hecke
$$
T_i^\circ :D_{\rm c}^{\rm b}({\cal E}_i'^\circ,{\bb C})\rightarrow
D_{\rm c}^{\rm b}({\cal C}{\it oh}_X^{0,n_0}\times {\cal E}_i^\circ,{\bb C})
$$
par
$$
T_i^\circ (K')=Rp_{i!}^\circ q_i^{\circ *}K'[in_0].
$$

\th LEMME 4.2
\enonce
{\rm (i)} Pour tout $K'\in {\rm ob}\,D_{\rm c}^{\rm b}({\cal C}_0',
{\bb C})$, on a un isomorphisme canonique
$$
T_1^\circ (j_0'^{\vee *}\pi_0'^{\vee *}K'[m_0'])\cong
({\rm id}\times j_0^\vee )^*
({\rm id}\times\pi_0^\vee )^*T_0(K')[m_0]
$$
dans $D_{\rm c}^{\rm b}({\cal C}{\it oh}_X^{0,n_0}\times
{\cal E}_1^\circ,{\bb C})$.

\decale{\rm (ii)} Pour tout $i=2,\ldots ,n$ et tout
$K'\in {\rm ob}\,D_{\rm c}^{\rm b}({\cal E}_{i-1}'^\circ,{\bb C})$
muni d'un isomorphisme $\mu_{i-1}'^*K'\cong {\bb C}\boxtimes K'$,
o\`u $\mu_{i-1}'$ est l'action de ${\bb G}_{\rm m}$ sur
${\cal E}_{i-1}'$, on a un isomorphisme canonique
$$
T_i^\circ \bigl(j_{i-1}'^{\vee *}
{\cal F}_{i-1}'(j_{i-1,!}'K')\bigr)\cong
({\rm id}\times j_{i-1}^\vee )^*
{\cal F}_{i-1,{\cal C}oh_X^{0,n_0}}
\bigl(({\rm id}\times j_{i-1})_!T_{i-1}^\circ (K')\bigr),
$$
dans $D_{\rm c}^{\rm b}({\cal C}{\it oh}_X^{0,n_0}\times
{\cal E}_i^\circ,{\bb C})$, o\`u ${\cal F}_{i-1}'$ et
${\cal F}_{i-1,{\cal C}oh_X^{0,n_0}}$ sont les
transformations de Fourier pour les fibr\'es
$\pi_{i-1}':{\cal E}_{i-1}'\rightarrow {\cal C}_{i-1}'$ et
${\rm id}\times\pi_{i-1}:{\cal C}{\it oh}_X^{0,n_0}\times
{\cal E}_{i-1}\rightarrow {\cal C}{\it oh}_X^{0,n_0}\times
{\cal C}_{i-1}$ respectivement.

\decale{\rm (iii)} Pour tout
$K'\in {\rm ob}\,D_{\rm c}^{\rm b}({\cal C}_n',{\bb C})$,
on a un isomorphisme canonique
$$
T_n^\circ (j_n'^*\pi_n'^*K')\cong
({\rm id}\times j_n)^*({\rm id}\times\pi_n)^*T_n(K')
$$
dans $D_{\rm c}^{\rm b}({\cal C}{\it oh}_X^{0,n_0}\times
{\cal E}_n^\circ,{\bb C})$.
\endth

\rem Preuve
\endrem
Pour chaque $i=1,\ldots ,n$,
introduisons le champ ${\cal K}_i$ des triplets
$(e_i',{\cal M}_i,\alpha_i)$ avec
$(\Omega_X^1)^{\otimes (i-1)}\buildrel
{e_i'}\over\longrightarrow {\cal M}_i'$ dans ${\cal E}_i'$ et
$({\cal M}_i,{\cal M}_i',\alpha_i)$ dans ${\cal H}_i$.
On peut identifier ${\cal J}_i$ au sous-champ de
${\cal K}_i$ form\'e des triplets
$(e_i',{\cal M}_i,\alpha_i)$ tels que $e_i'$ se factorise en
$$
(\Omega_X^1)^{\otimes (i-1)}
\buildrel {e_i}\over\longrightarrow {\cal M}_i
\buildrel {\alpha_i}\over\longrightarrow
{\cal M}_i'.
$$
Pour chaque $\alpha_i$, on a la suite exacte courte
$$\displaylines{
\qquad 0\rightarrow
{\rm Hom}_{{\cal O}_X}((\Omega_X^1)^{\otimes (i-1)},{\cal M}_i)
\rightarrow
{\rm Hom}_{{\cal O}_X}((\Omega_X^1)^{\otimes (i-1)},{\cal M}_i')
\hfill\cr\hfill
\rightarrow
{\rm Hom}_{{\cal O}_X}((\Omega_X^1)^{\otimes (i-1)},
{\rm Coker}(\alpha_i))\rightarrow 0\qquad}
$$
(on a ${\rm Hom}_{{\cal O}_X}({\cal M}_i,
(\Omega_X^1)^{\otimes i})=(0)$) et ${\cal J}_i$ est donc
un sous-fibr\'e vectoriel de ${\cal K}_i$. On notera
$k_i:{\cal J}_i\hookrightarrow {\cal K}_i$ l'inclusion.
Bien entendu, si on note $\widetilde q_i$ la projection canonique
de ${\cal K}_i$ sur ${\cal E}_i'$, on a $\widetilde q_i\circ k_i=q_i$.

Dualement, pour chaque $i=0,\ldots ,n-1$,
${\cal J}_i^\vee$ est un fibr\'e vectoriel quotient de
${\cal K}_i^\vee$. Plus pr\'ecis\'ement

\decale{$\bullet$} ${\cal J}_i^\vee$ est le champ des triplets
$(e_i^\vee ,{\cal M}_i',\alpha_i)$ avec
$e_i^\vee =((\Omega_X^1)^{\otimes i}\hookrightarrow {\cal M}_{i+1}
\twoheadrightarrow {\cal M}_i)$ dans ${\cal E}_i^\vee$ et
$({\cal M}_i,{\cal M}_i',\alpha_i)$ dans ${\cal H}_i$,

\decale{$\bullet$} ${\cal K}_i^\vee$ est le champ des triplets
$(e_i^\vee ,e_i'^\vee ,\alpha_i)$
avec $e_i^\vee =((\Omega_X^1)^{\otimes i}\hookrightarrow
{\cal M}_{i+1}\twoheadrightarrow{\cal M}_i)$
dans ${\cal E}_i^\vee$,  $e_i'^\vee =
((\Omega_X^1)^{\otimes i}\hookrightarrow {\cal M}_{i+1}'
\twoheadrightarrow {\cal M}'_i)$ dans ${\cal E}_i'^\vee$,
$({\cal M}_i,{\cal M}_i',\alpha_i)$ dans ${\cal H}_i$ et
$e_i^\vee =\alpha_i^*(e_i'^\vee )$, i.e. avec un diagramme
commutatif
$$
\def\normalbaselines{\baselineskip20pt\lineskip3pt
\lineskiplimit3pt}
\matrix{(\Omega_X^1)^{\otimes i}&\hookrightarrow&
{\cal M}_{i+1}&\twoheadrightarrow&{\cal M}_i\cr
\Vert&&\Big\downarrow&\square&\mapdownright{\alpha_i}\cr
(\Omega_X^1)^{\otimes i}&\hookrightarrow&{\cal M}_{i+1}'&
\twoheadrightarrow&{\cal M}_i'\,,\cr}
$$

\decale{$\bullet$} le morphisme
$k_i^\vee :{\cal K}_i^\vee\twoheadrightarrow
{\cal J}_i^\vee$ envoie
$(e_i^\vee ,e_i'^\vee,\alpha_i)$ sur
$(e_i^\vee ,{\cal M}_i',\alpha_i)$.
\vskip 1mm
\noindent De plus, $\widetilde p_i=p_i^\vee\circ k_i^\vee $
est la projection canonique de
${\cal K}_i^\vee$ sur $X\times {\cal E}_i$.

On a alors un diagramme commutatif de champs
\vskip 1mm
$$
\def\normalbaselines{\baselineskip20pt\lineskip3pt
\lineskiplimit3pt}
\matrix{{\cal C}{\it oh}_X^{0,n_0}\times {\cal E}_i&
\mapleftover{p_i}&{\cal J}_i&
\buildrel{k_i}\over\hookrightarrow &{\cal K}_i
&\maprightover{\tilde q_i}&{\cal E}_i'\cr
\mapdownleft{{\rm id}\times\pi_i}&\square &\Big\downarrow
&&\Big\downarrow
&\square&
\mapdownright{\pi_i'}\cr
{\cal C}{\it oh}_X^{0,n_0}\times{\cal C}_i&\mapleftover{p_i}
&{\cal H}_i&{=\!=}&{\cal H}_i&
\maprightover{q_i}&{\cal C}_i'\cr}
$$
et dualement un diagramme commutatif de champs
\vskip 1mm
$$
\def\normalbaselines{\baselineskip20pt\lineskip3pt
\lineskiplimit3pt}
\matrix{{\cal C}{\it oh}_X^{0,n_0}\times {\cal E}_i^\vee&
\mapleftover{p_i^\vee}&{\cal J}_i^\vee&
\buildrel{k_i^\vee}\over\twoheadleftarrow &{\cal K}_i^\vee
&\maprightover{\tilde q_i^\vee}&{\cal E}_i'^\vee\cr
\mapdownleft{{\rm id}\times\pi_i^\vee}&\square &
\Big\downarrow&&\Big\downarrow&\square&
\mapdownright{\pi_i'^\vee}\cr
{\cal C}{\it oh}_X^{0,n_0}\times{\cal C}_i&\mapleftover{p_i}&
{\cal H}_i&{=\!=}&{\cal H}_i&\maprightover{q_i}&{\cal C}_i'\,.\cr}
$$

On notera  $T_{i-1}^\vee$ l'op\'erateur de Hecke
$$
D_{\rm c}^{\rm b}({\cal E}_{i-1}'^\vee,{\bb C})\rightarrow
D_{\rm c}^{\rm b}({\cal C}{\it oh}_X^{0,n_0}\times{\cal E}_{i-1}^\vee,{\bb C}),
{}~K'\mapsto R\widetilde p_{i-1,!}^\vee
\widetilde q_{i-1}^{\vee *}K'[in_0].
$$
On peut identifier ${\cal J}_i^\circ$ \`a un ouvert de
${\cal K}_{i-1}^\vee$ et on a le diagramme commutatif de champs
\vskip 1mm
$$
\def\normalbaselines{\baselineskip20pt\lineskip3pt
\lineskiplimit3pt}
\matrix{{\cal C}{\it oh}_X^{0,n_0}\times {\cal E}_i^\circ&
\mapleftover{p_i^\circ}&{\cal J}_i^\circ&
\maprightover{q_i^\circ}&{\cal E}_i'^\circ\cr
\mapdownleft{{\rm id}\times j_{i-1}^\vee}&\square
&\Big\downarrow &&\mapdownright{j_{i-1}'^\vee}\cr
{\cal C}{\it oh}_X^{0,n_0}\times {\cal E}_{i-1}^\vee
&\mapleftover{{\tilde p}_{i-1}^\vee}&{\cal K}_{i-1}^\vee
&\maprightover{{\tilde q}_{i-1}^\vee} &
{\cal E}_{i-1}'^\vee\,.\cr}
$$
Il est donc clair que
$$
T_i^\circ\circ j_{i-1}'^{\vee *}\cong
({\rm id}\times j_{i-1}^\vee )^*\circ T_{i-1}^\vee .
$$

Pour la partie (i) du lemme on conclut en remarquant que
$$
T_0^\vee (\pi_0'^{\vee *}K'[m_0'])\cong
({\rm id}\times\pi_0^\vee )^*T_0(K')[m_0]
$$
puisque $Rk_{0!}^\vee {\bb C}={\bb C}[-2n_0]$.

Pour la partie (ii), comme on a trivialement
$$
T_{i-1}\circ j_{i-1,!}'\cong
({\rm id}\times j_{i-1})_!\circ T_{i-1}^\circ ,
$$
il ne nous reste plus qu'\`a d\'emontrer que
$$
T_{i-1}^\vee\circ {\cal F}_{i-1}'\cong
{\cal F}_{i-1,{\cal C}oh_X^{0,n_0}}\circ T_{i-1}.
$$
Or, on a
$$
T_{i-1}^\vee\circ {\cal F}_{i-1}(K')\cong
R\widetilde p_{i-1,!}^\vee
{\cal F}_{{\cal K}_{i-1}}(\widetilde q_{i-1}^*K')[in_0],
$$
o\`u ${\cal F}_{{\cal K}_{i-1}}$ est la transformation de Fourier
pour le fibr\'e ${\cal K}_{i-1}\rightarrow {\cal H}_{i-1}$,
on a
$$
{\cal F}_{i-1,{\cal C}oh_X^{0,n_0}}\circ T_{i-1}(K')\cong
Rp_{i-1,!}^\vee {\cal F}_{{\cal J}_{i-1}}(q_{i-1}^*K')[(i-1)n_0],
$$
o\`u ${\cal F}_{{\cal J}_{i-1}}$ est la transformation de Fourier
pour le fibr\'e ${\cal J}_{i-1}\rightarrow {\cal H}_{i-1}$,
et on a
$$\eqalign{
R\widetilde p_{i-1,!}^\vee
{\cal F}_{{\cal K}_{i-1}}(\widetilde q_{i-1}^*K')
&\cong Rp_{i-1,!}^\vee Rk_{i-1,!}^\vee
{\cal F}_{{\cal K}_{i-1}}(\widetilde q_{i-1}^*K')\cr
&\cong R p_{i-1,!}^\vee
{\cal F}_{{\cal J}_{i-1}}(k_{i-1}^*\widetilde q_{i-1}^*K')[-n_0]
\cong Rp_{i-1,!}^\vee {\cal F}_{{\cal J}_{i-1}}
(q_{i-1}^*K')[-n_0],\cr}
$$
d'o\`u la conclusion.

La partie (iii) est imm\'ediate.

\hfill\hfill\cqfd

\th COROLLAIRE 4.3
\enonce
Pour chaque $i=1,\ldots ,n$, on a
$$
T_i^\circ (W_{L,i}'^\circ )\cong V_L
\boxtimes W_{L,i}^\circ
$$
dans $D_{\rm c}^{\rm b}({\cal C}{\it oh}_X^{0,n_0}
\times{\cal E}_i^\circ,{\bb C})$, o\`u $V_L$ est le faisceau
pervers irr\'eductible de Whittaker associ\'e \`a $L$ sur
${\cal C}{\it oh}_X^{0,n_0}$ et
o\`u $W_{L,i}^\circ$ et $W_{L,i}'^\circ$ sont
les objets de $D_{\rm c}^{\rm b}({\cal E}_i^\circ,{\bb C})$ et
$D_{\rm c}^{\rm b}({\cal E}_i'^\circ,{\bb C})$ respectivement
qui ont \'et\'e associ\'es \`a $L$ dans la section $3$.
\endth

\th CONJECTURE 4.4
\enonce
{\rm (i)} Pour chaque $i=1,\ldots ,n$, on a
$$
T_i(W_{L,i}')\cong V_L \boxtimes W_{L,i}
$$
dans $D_{\rm c}^{\rm b}({\cal C}{\it oh}_X^{0,n_0}
\times{\cal E}_i,{\bb C})$, o\`u $W_{L,i}$ et $W_{L,i}'$ sont les
objets de $D_{\rm c}^{\rm b}({\cal E}_i,{\bb C})$ et
$D_{\rm c}^{\rm b}({\cal E}_i',{\bb C})$ respectivement
qui ont \'et\'e associ\'es \`a $L$ dans la section $3$.

\decale{\rm (ii)} On a
$$
T_n (K_L')\cong V_L[-n_0]\boxtimes K_L
$$
dans $D_{\rm c}^{\rm b}({\cal C}{\it oh}_X^{0,n_0}\times
{\cal C}_n,{\bb C})$, o\`u $K_L$ et $K_L'$ sont les objets de
$D_{\rm c}^{\rm b}({\cal C}_n,{\bb C})$ et
$D_{\rm c}^{\rm b}({\cal C}_n',{\bb C})$ respectivement
qui sont associ\'es \`a $L$ par la conjecture $3.1$.
\endth

\vfill\eject
\centerline{\douzebf Expos\'e II}
\vskip 20mm
{\bf 1. Passage \`a la caract\'eristique $p$.}
\vskip 5mm

Dans ce deuxi\`eme expos\'e on fixe un corps fini ${\bb F}_q$ de
caract\'eristique $p$, une cl\^oture alg\'ebrique $k$ de ce corps
fini et une courbe $X$ connexe, projective et lisse sur $k$, de genre
$g\geq 2$ et d\'efinie sur ${\bb F}_q$. Soit $|X|$ l'ensemble des
orbites de ${\rm Gal}(k/{\bb F}_q)$ dans  $X(k)$. Pour chaque
$x\in |X|$, on note ${\rm deg}(x)$ le cardinal de $x$,
$\kappa (x)$ l'unique sous-corps  \`a $q_x=q^{{\rm deg}(x)}$
\'el\'ements de $k$ et ${\rm Frob}_x$ l'\'el\'evation \`a la
puissance $q_x^{-1}$ dans $k$ (c'est un g\'en\'erateur topologique
de ${\rm Gal}(k/\kappa (x))$). Pour chaque $x\in |X|$, on
choisit une fois pour toute un repr\'esentant $\overline x\in X(k)$
de cette orbite; ${\rm Gal}(k/\kappa (x))\subset
{\rm Gal}(k/{\bb F}_q)$ est le fixateur de $\overline x$.

On fixe les entiers $n\geq 1$ et $m_0>n(3n-1)(g-1)$ comme dans
l'expos\'e I et on construit comme pr\'ec\'edemment le diagramme
fondamental
$$
\matrix{&&\kern -3mm{\cal E}_n\supset {\cal E}_n^\circ\cong
{\cal E}_{n-1}^{\vee\circ}\subset
{\cal E}_{n-1}^\vee\kern -3mm &\cr
&\swarrow&&\searrow \cr
{\cal C}_n\kern -3mm&&&\cr}\quad\cdots\quad
\matrix{&&&\kern -3mm{\cal E}_1\supset {\cal E}_1^\circ\cong
{\cal E}_0^{\vee\circ}
\subset {\cal E}_0^\vee\kern -3mm&&\cr
\searrow &&\swarrow &&\searrow &\cr
&\kern -3mm{\cal C}_1\kern -3mm &&&&\kern -3mm{\cal C}_0\cr}
$$
qui est maintenant un diagramme de champs alg\'ebriques sur $k$,
d\'efini sur ${\bb F}_q$.

On fixe de plus un nombre premier $\ell\not= p$ et
une cl\^oture alg\'ebrique $\overline{\bb Q}_\ell$ du corps des
nombres $\ell$-adiques. Alors, pour tout
$\overline{\bb Q}_\ell$-faisceau $L$,
lisse, irr\'eductible de rang $n$ sur $X$ et
d\'efini sur ${\bb F}_q$, on peut former, comme
dans l'expos\'e I, le faisceau de Whittaker $W_L$, qui est
maintenant un $\overline{\bb Q}_\ell$-faisceau pervers
irr\'eductible sur ${\cal C}_0$, d\'efini sur ${\bb F}_q$.

On fixe en outre un caract\`ere additif non trivial
$\psi :{\bb F}_q\rightarrow\overline{\bb Q}_\ell^\times$
et, pour tout $i=1,\ldots ,n-1$, on note
$$
{\cal F}_{\psi ,i}: D_{\rm c}^{\rm b}({\cal E}_i,\overline{\bb Q}_\ell )
\rightarrow D_{\rm c}^{\rm b}({\cal E}_i^\vee,\overline{\bb Q}_\ell )
$$
la transformation de Fourier pour le fibr\'e
vectoriel ${\cal E}_i\rightarrow {\cal C}_i$ d\'efinie par
$$
{\cal F}_{\psi ,i}(K_i)=R{\rm pr}_{i,!}^\vee({\rm pr}_i^*K_i\otimes
{\cal L}_\psi (\langle\cdot ,\cdot\rangle_i ))[m_0-i^2(g-1)],
$$
o\`u ${\rm pr}_i^\vee$ et ${\rm pr}_i$ sont les deux projections
canoniques de ${\cal E}_i^\vee\times_{{\cal C}_i}{\cal E}_i$, o\`u
$\langle\cdot ,\cdot\rangle_i:
{\cal E}_i^\vee\times_{{\cal C}_i}{\cal E}_i
\rightarrow {\bb A}_k^1$
est l'accouplement canonique et o\`u ${\cal L}_\psi$ est le
$\overline{\bb Q}_\ell$-faisceau lisse de rang un sur
${\bb A}_k^1$ d'Artin-Schreier
associ\'e \`a $\psi$.
Cette transformation de Fourier envoie
$D_{\rm mon}^{\rm b}({\cal E}_i,\overline{\bb Q}_\ell )$
dans $D_{\rm mon}^{\rm b}({\cal E}_i^\vee,\overline{\bb Q}_\ell )$
et, si $K_i\in {\rm ob}\,D_{\rm mon}^{\rm b}({\cal E}_i,
\overline{\bb Q}_\ell )$ est muni d'un isomorphisme
$\mu_i^*K_i\cong \overline{\bb Q}_\ell\boxtimes K_i$,
${\cal F}_{\psi ,i}(K_i)$ est aussi muni d'un isomorphisme
$$
\mu_i^{\vee *}{\cal F}_{\psi ,i}(K_i)\cong
\overline{\bb Q}_\ell\boxtimes{\cal F}_{\psi ,i}(K_i)
$$
et ${\cal F}_{\psi ,i}(K_i)$ est ind\'ependant du choix de $\psi$.
On peut former \`a l'aide de ces transformations de Fourier
les objets $W_{L,i}$ et $W_{L,i}^\circ$.
Les conjectures 3.1 et 3.2 gardent un sens.

On peut enfin d\'efinir les op\'erateurs de Hecke $T_i$ et
$T_i^\circ$ et les r\'esultats de la section 4 de l'expos\'e I valent
sans modification, \`a des torsions \`a la Tate
pr\`es que l'on va pr\'eciser ci-dessous.

\vskip 5mm
{\bf 2. Fonction trace de Frobenius de $W_L$.}
\vskip 5mm

Rappelons qu'une partition $\mu =(\mu_1,\ldots ,\mu_n)$ de
longueur $\leq n$ est une suite d'entiers telle que
$\mu_1\geq\mu_2\geq\cdots\geq\mu_n\geq 0$.
On munit l'ensemble de ces partitions de la relation d'ordre usuelle
$$
\lambda\geq\mu\Longleftrightarrow
\left\{\matrix{\lambda_1\geq\mu_1,\hfill\cr
\noalign{\smallskip}
\lambda_1+\lambda_2\geq\mu_1+\mu_2,\hfill\cr
\noalign{\smallskip}
\cdots\hfill\cr
\noalign{\smallskip}
\lambda_1+\cdots +\lambda_n=\mu_1+\cdots +\mu_n.\hfill\cr}
\right.
$$
Pour toute paire $(\lambda, \mu )$ de partitions de longueur
$\leq n$, on note
$$
K_{\lambda\mu}(t)
$$
le polyn\^ome de Kostka-Foulkes correspondant.
Rappelons que $K_{\lambda\mu}(t)$
est un polyn\^ome unitaire \`a coefficients entiers positifs ou nuls
et de degr\'e $n(\mu )-n(\lambda )$ si
$\lambda\geq\mu$, o\`u on a pos\'e
$n(\nu )=\sum_{i=1}^n(i-1)\nu_i$ pour
toute partition $\nu$ de longueur $\leq n$,
et que $K_{\lambda\mu}(t)=0$ sinon (cf. [Ma] Ch. III, \S 6).
Pour toute paire $(\lambda, \mu )$ de partitions
de longueur $\leq n$, on pose
$$
\widetilde K_{\lambda\mu}(t)=t^{n(\mu )}
K_{\lambda\mu}(t^{-1}),
$$
de sorte que $\widetilde K_{\lambda\mu}(t)$ est un polyn\^ome de
la forme
$$
t^{n(\lambda )}+\hbox{termes d'ordre sup\'erieur},
$$
\`a coefficients entiers positifs ou nuls et de degr\'e
$\leq n(\mu )$.

Pour tout espace vectoriel $V$ de dimension $n$ et toute partition
$\lambda$ de longueur $\leq n$, on note
$$
R^{(\lambda_1-\lambda_2,\ldots ,
\lambda_{n-1}-\lambda_n,\lambda_n)}V
$$
la repr\'esentation de
plus haut poids $(\lambda_1-\lambda_2)\omega_1+\cdots +
(\lambda_{n-1}-\lambda_n)\omega_{n-1}+\lambda_n\omega_n$
de $GL(V)$, o\`u  $\omega_i$ est bien entendu le plus haut poids
de $\bigwedge^iV$.

Pour tout $x\in |X|$ et toute partition $\mu =(\mu_1,\ldots ,\mu_n)$
de longueur $\leq n$, on pose
$$
w_{L,x}^{\rm top}(\mu_1-\mu_2,\ldots ,\mu_{n-1}-\mu_n,\mu_n)
=q_x^{n(\mu )}{\rm tr}({\rm Frob}_x,
R^{(\mu_1-\mu_2,\ldots ,\mu_{n-1}-\mu_n,\mu_n)}L_{\overline x})
$$
et
$$
w_{L,x}(\mu_1-\mu_2,\ldots ,\mu_{n-1}-\mu_n,\mu_n)
=\sum_{\lambda\geq\mu}\widetilde K_{\lambda\mu}(q_x)
{\rm tr}({\rm Frob}_x,R^{(\lambda_1-
\lambda_2,\ldots ,\lambda_{n-1}-\lambda_n,
\lambda_n)}L_{\overline x}).
$$

Alors, pour tout uplet $(D_1,\ldots ,D_n)$ de
diviseurs effectifs sur $X$, d\'efinis sur ${\bb F}_q$, avec
$\sum_{i=1}^ni{\rm deg}(D_i)=m_0$ et
$D_i=\sum_{x\in |X|}d_{i,x}x$, on introduit les expressions
$$
w_L^{\rm top}(D_1,\ldots ,D_n)=\prod_{x\in |X|}w_{L,x}^{\rm top}
(d_{1,x},\ldots ,d_{n,x})
$$
et
$$
w_L(D_1,\ldots ,D_n)=\prod_{x\in |X|}w_{L,x}(d_{1,x},\ldots ,d_{n,x}).
$$

\th TH\'EOR\`EME 2.1 ([La] Lemme (3.3.6) et Thm. (3.3.8))
\enonce
Soit $(D_1,\ldots ,D_n)$ un $n$-uplet de diviseurs effectifs
sur $X$, d\'efinis sur ${\bb F}_q$,
avec $\sum_{i=1}^ni{\rm deg}(D_i)=m_0$. Alors :
\decale{\rm (i)} la fibre du $\overline{\bb Q}_\ell$-faisceau
pervers irr\'eductible $W_L$ en ${\cal O}_{D_1+\cdots +D_n}\oplus
\cdots\oplus {\cal O}_{D_n}
\in {\rm ob}\,{\cal C}_0({\bb F}_q)$  a sa cohomologie
concentr\'ee en degr\'es compris entre $0$ et
$\sum_{i=1}^ni(i-1){\rm deg}(D_i)$,
\decale{\rm (ii)} $w_L^{\rm top}(D_1,\ldots ,D_n)$ est la trace de
l'endomorphisme
de Frobenius g\'eom\'etrique  sur le groupe de cohomologie de
degr\'e $\sum_{i=1}^ni(i-1){\rm deg}(D_i)$ de cette fibre,
\decale{\rm (iii)} $w_L(D_1,\ldots ,D_n)$ est la somme altern\'ee
des traces de l'endomorphisme de Frobenius
g\'eom\'etrique sur les groupes de cohomologie ce cette fibre.
\hfill\hfill\cqfd
\endth
\vskip 2mm

{\pc REMARQUE} 2.2\pointir
L'\'enonc\'e du lemme (3.3.6) de [La] contient une erreur~:
il faut y remplacer
$K_{\underline m',\underline m}(q)$ par
$\widetilde K_{\underline m',\underline m}(q)$
et donc, pour chaque $a$,
$K_{\underline m',\underline m,a}$ par
$\widetilde K_{\underline m',\underline m,a}$
(cf. [Ma] Ch. III, \S 7, Example 9).

\hfill\hfill\cqfd

\vskip 5mm
{\bf 3. La fonction trace de Frobenius de $W_{L,n}^\circ$.}
\vskip 5mm

Pour chaque $i=1,\ldots ,n$, introduisons d'une part l'ensemble
$$
E_i=\{({\cal A}_i\hookrightarrow {\cal L}_i\twoheadrightarrow
{\cal L}_{i-1};D_i,\ldots ,D_n)\}
$$
des uplets tels que~:

\decale{$\bullet$} chaque $D_j=\sum_{x\in |X|}d_{j,x}x$
($j=i,\ldots ,n$)
est un diviseur effectif sur $X$, d\'efini sur ${\bb F}_q$,

\decale{$\bullet$} ${\cal A}_i
=(\Omega_X^1)^{\otimes (i-1)}(D_i+\cdots +D_n)$,

\decale{$\bullet$} ${\cal L}_{i-1}$ est un
${\cal O}_X$-Module localement libre de rang $i-1$ sur $X$, d\'efini
sur ${\bb F}_q$, tel que ${\rm Hom}_{{\cal O}_X}({\cal L}_{i-1},
(\Omega_X^1)^{\otimes (2n-1)})=(0)$,

\decale{$\bullet$} ${\cal A}_i\hookrightarrow {\cal L}_i
\twoheadrightarrow{\cal L}_{i-1}$ est une extension
(\`a isomorphisme pr\`es) de ${\cal O}_X$-Modules, d\'efinie sur
${\bb F}_q$, de sorte que
${\cal L}_i$ est un ${\cal O}_X$-Module
localement libre de rang $i$ sur $X$,
d\'efini sur ${\bb F}_q$,

\decale{$\bullet$} ${\rm deg}({\cal L}_i)+{\rm deg}(D_{i+1})+
2{\rm deg}(D_{i+2})+
\cdots +(n-i){\rm deg}(D_n)=m_i$.
\vskip 1mm

Introduisons d'autre part la fonction
$$
w_{L,i}^{\rm top}:E_i\rightarrow \overline{\bb Q}_\ell
$$
d\'efinie par r\'ecurrence sur $i$ de la fa\c con suivante.
Pour $i=1$, on pose
$$
w_{L,1}^{\rm top}({\cal A}_1\buildrel\sim\over\rightarrow
{\cal L}_1\rightarrow 0;D_1,\ldots ,D_n)
=w_L^{\rm top}(D_1,\ldots ,D_n).
$$
Pour $i=2,\ldots ,n$, on pose
$$\displaylines{
\qquad w_{L,i}^{\rm top}({\cal A}_{i}\hookrightarrow {\cal L}_{i}
\twoheadrightarrow
{\cal L}_{i-1};D_i,\ldots ,D_n)
\hfill\cr\hfill
=\sum_{e_{i-1}}\psi(\langle e_{i-1}^\vee ,e_{i-1}\rangle )
w_{L,i-1}^{\rm top}({\cal A}_{i-1}\hookrightarrow {\cal L}_{i-1}
\twoheadrightarrow{\cal L}_{i-2};D_{i-1},\ldots ,D_n),\qquad}
$$
o\`u:

\decale{$\bullet$} $e_{i-1}$ parcourt les homomorphismes de
$(\Omega_X^1)^{\otimes (i-2)}(D_i+\cdots +D_n)$ dans
${\cal L}_{i-1}$ qui sont non nuls (et donc
injectifs) et qui sont d\'efinis sur ${\bb F}_q$,

\decale{$\bullet$} $D_{i-1}$ est l'unique diviseur effectif sur
$X$ (d\'efini sur ${\bb F}_q$) tel que $e_{i-1}$ se factorise par
l'inclusion
$$
(\Omega_X^1)^{\otimes (i-2)}(D_i+\cdots +D_n)\hookrightarrow
(\Omega_X^1)^{\otimes (i-2)}(D_{i-1}+D_i+\cdots +D_n)
={\cal A}_{i-1}
$$
et que le quotient ${\cal L}_{i-2}$ de
${\cal L}_{i-1}$ par ${\cal A}_{i-1}$ soit un
${\cal O}_X$-Module localement libre (de rang $i-2$),

\decale{$\bullet$}  $e_{i-1}^\vee$ est la classe de l'extension
${\cal A}_i\hookrightarrow {\cal L}_i
\twoheadrightarrow {\cal L}_{i-1}$ dans
${\rm Ext}_{{\cal O}_X}^1({\cal L}_{i-1},{\cal A}_i)$,

\decale{$\bullet$} $\langle\cdot ,\cdot\rangle$
est l'accouplement pour la dualit\'e de Serre entre
$$
{\rm Hom}_{{\cal O}_X}((\Omega_X^1)^{\otimes (i-2)}
(D_i+\cdots +D_n),{\cal L}_{i-1})
$$
et
$$
{\rm Ext}_{{\cal O}_X}^1({\cal L}_{i-1},{\cal A}_i).
$$
\vskip 2mm

{\pc REMARQUE} 3.1\pointir
Les ensembles $E_i$ et les fonctions $w_L^{\rm top}$,
$w_{L,i}^{\rm top}$ sont des variantes des ensembles
${\cal P}_i$ et des fonctions $W_E$, $V_{E,i}$ introduits
dans [La] \S 4: les notations et les normalisations
sont l\'eg\`erement diff\'erentes et les transformations
de Radon ont \'et\'e remplac\'ees par des
transformations de Fourier, ce qui revient au m\^eme car
les fonctions $w_{L,i}^{\rm top}$ sont
constantes sur les droites \'epoint\'ees.

\hfill\hfill\cqfd

Pour $i=1,\ldots ,n$, consid\'erons le sous-ensemble
$$
E_i^\circ\subset E_i
$$
d\'efini par la condition
$$
{\rm Hom}_{{\cal O}_X}({\cal L}_i,(\Omega_X^1)^{\otimes (2n-1)})
=(0).
$$
A chaque
$({\cal A}_i\hookrightarrow {\cal L}_i\twoheadrightarrow
{\cal L}_{i-1};D_i,\ldots ,D_n)
\in E_i^\circ$, on associe l'objet
$$
{\cal M}_i\buildrel {\rm dfn}\over {=\!=}
{\cal L}_i\oplus {\cal O}_{D_{i+1}+\cdots +D_n}\oplus\cdots\oplus
{\cal O}_{D_n}
$$
de ${\cal C}_i({\bb F}_q)$ et on d\'efinit une application
$$
E_i^\circ\rightarrow {\rm ob}\,{\cal E}_i^\circ ({\bb F}_q)
$$
en envoyant
$({\cal A}_i\hookrightarrow {\cal L}_i\twoheadrightarrow
{\cal L}_{i-1};D_i,\ldots ,D_n)$
sur le momomorphisme $(\Omega_X^1)^{\otimes (i-1)}
\hookrightarrow {\cal M}_i$
somme directe de l'application compos\'ee
$$
(\Omega_X^1)^{\otimes (i-1)}\hookrightarrow
(\Omega_X^1)^{\otimes (i-1)}(D_i+\cdots +D_n)={\cal A}_i
\hookrightarrow {\cal L}_i
$$
et de l'application nulle de $(\Omega_X^1)^{\otimes (i-1)}$ dans
${\cal O}_{D_{i+1}+\cdots +D_n}\oplus\cdots\oplus {\cal O}_{D_n}$.

La conjecture suivante est certainement plus accessible que les
pr\'ec\'edentes. Une m\'ethode permettant probablement de
la prouver est indiqu\'ee dans la section suivante (voir la remarque
4.5).

\th CONJECTURE 3.2
\enonce
Il existe une constante $c_n\in\overline{\bb Q}_\ell^\times$ telle
que, pour tout \'el\'ement $({\cal A}_n\hookrightarrow {\cal L}_n
\twoheadrightarrow {\cal L}_{n-1};D_n)$ dans $E_n^\circ$ d'image
$(\Omega_X^1)^{\otimes (n-1)}\hookrightarrow {\cal M}_n$
dans ${\rm ob}\,{\cal E}_n^\circ ({\bb F}_q)$, la trace de
l'endomorphisme de Frobenius g\'eom\'etrique sur la fibre de
$W_{L,n}^\circ\in {\rm ob}\,D_{\rm c}^{\rm b}({\cal E}_n^\circ ,
\overline{\bb Q}_\ell )$ en $(\Omega_X^1)^{\otimes (n-1)}
\hookrightarrow {\cal M}_n$ soit \'egale \`a
$$
c_nw_{L,n}^{\rm top}({\cal A}_n\hookrightarrow {\cal L}_n
\twoheadrightarrow {\cal L}_{n-1};D_n).
$$
\endth

Pour $n=2$, il est possible de prouver cette conjecture par un
calcul direct. C'est ce que nous allons faire pour terminer cette
section. Nous aurons besoin des cas  particuliers $n=2$ des deux
r\'esultats suivants.

\th LEMME 3.3
\enonce
Soit ${\cal M}_i={\cal L}_i\oplus {\cal T}_i$ un
${\cal O}_X$-Module  coh\'erent de rang
g\'en\'erique $i\in\{1,\ldots ,n\}$, avec ${\cal L}_i$ localement
libre et ${\cal T}_i={\cal O}_{D_{i+1}+\cdots +D_n}\oplus\cdots
\oplus {\cal O}_{D_n}$ pour des diviseurs effectifs
$D_{i+1},\ldots ,D_n$ sur $X$. Soit $e:(\Omega_X^1)^{\otimes (i-1)}
\hookrightarrow {\cal M}_i$ un homomorphisme injectif de
${\cal O}_X$-Modules, somme directe d'homomorphismes de
${\cal O}_X$-Modules $e_i:(\Omega_X^1)^{\otimes (i-1)}
\hookrightarrow {\cal L}_i$ et $e_j:(\Omega_X^1)^{\otimes (i-1)}
\rightarrow {\cal O}_{D_j+\cdots +D_n}$, $j=i+1,\ldots ,n$.
On note $A_i$ l'unique diviseur effectif sur $X$ tel que
${\rm Coker}(e_i)$ soit isomorphe \`a
${\cal L}_{i-1}\oplus {\cal O}_{A_i}$, avec ${\cal L}_{i-1}$
localement libre de rang
$i-1$, et, pour chaque $j=i+1,\ldots ,n$, on note $A_j$ l'unique
diviseur effectif sur $X$
tel que ${\rm Coker}(e_j)$ soit isomorphe \`a ${\cal O}_{A_j}$
{\rm (}bien entendu, on a
$A_j\leq D_j+\cdots +D_n${\rm )}.

Alors, ${\rm Coker}(e)$ est isomorphe \`a
$$
{\cal L}_{i-1}\oplus {\cal O}_{D'_i+\cdots +D'_n}\oplus\cdots
\oplus {\cal O}_{D'_n}
$$
o\`u les diviseurs effectifs $D'_i,\ldots ,D'_n$ sont d\'efinis
par les relations
$$
D'_i+2D'_{i+1}+\cdots +(n-i+1)D'_n=
A_i+(D_{i+1}+2D_{i+2}+\cdots +(n-i)D_n)
$$
et, pour $j=i+1,\ldots ,n$,
$$\displaylines{
\qquad
D'_j+2D'_{j+1}+\cdots +(n-j+1)D'_n
\hfill\cr\hfill
\eqalign{={\rm Inf}\{&A_k+(D_{j+1}+2D_{j+2}+\cdots +(n-j)D_n),\cr
&A_\ell +(D_j+2D_{j+1}+\cdots +(n-j+1)D_n)-
(D_\ell+\cdots +D_n)\},\cr}\qquad}
$$
l'indice $k$ parcourant $\{i,\ldots ,j\}$ et
l'indice $\ell$ parcourant $\{j+1,\ldots ,n\}$.
\endth

\rem Preuve
\endrem
Si ${\cal O}$ est un anneau de valuation discr\`ete, d'uniformisante
$\varpi$, et si
$a_i,\ldots ,a_n$ et $d_{i+1},\ldots ,d_n$ sont des entiers positifs
ou nuls avec
$a_j\leq d_j+\cdots +d_n$ pour $j=i+1,\ldots ,n$, les mineurs
non nuls de la matrice
$$
\pmatrix{
\varpi^{a_i}&0&\ldots&\ldots&0\cr
\varpi^{a_{i+1}}&\varpi^{d_{i+1}+\cdots +d_n}&\ddots&&\vdots\cr
\varpi^{a_{i+2}}&0&\ddots&\ddots&\vdots\cr
\vdots&\vdots &\ddots&\ddots&0\cr
\varpi^{a_n}&0&\ldots&0&\varpi^{d_n}\cr}
$$
sont (au signe pr\`es) les
$$
M_{I,I}=\varpi^{a_i+\sum_{j\in I-\{i\}}(d_j+\cdots +d_n)}
$$
avec $i\in I\subset \{i,\ldots ,n\}$, les
$$
M_{J,J}=\varpi^{\sum_{j\in J}(d_j+\cdots +d_n)}
$$
avec $J\subset\{i+1,\ldots ,n\}$ et les
$$
M_{J,K}=\varpi^{a_k+\sum_{j\in J-\{k\}}(d_j+\cdots +d_n)}
$$
avec $J\subset \{i+1,\ldots ,n\}$ et $K=J\cup\{i\}-\{k\}$
pour un $k\in J$. Par suite, les
diviseurs \'el\'ementaires de cette matrice sont
$$
\varpi^{d'_n}|\varpi^{d'_{n-1}+d'_n}|\cdots |\varpi^{d'_i+\cdots +d'_n}
$$
avec
$$
d'_i+2d'_{i+1}+\cdots +(n-i+1)d'_n=
a_i+(d_{i+1}+2d_{i+2}+\cdots +(n-i)d_n)
$$
et, pour $j=i+1,\ldots ,n$,
$$\displaylines{
\qquad
d'_j+2d'_{j+1}+\cdots +(n-j+1)d'_n
\hfill\cr\hfill
\eqalign{={\rm Inf}\{&a_k+(d_{j+1}+2d_{j+2}+\cdots +(n-j)d_n),\cr
&a_\ell +(d_j+2d_{j+1}+\cdots +(n-j+1)d_n)-
(d_\ell+\cdots +d_n)\},\cr}\qquad}
$$
l'indice $k$ parcourant
$\{i,\ldots ,j\}$ et l'indice
$\ell$ parcourant $\{j+1,\ldots ,n\}$ (on a
$$
\sum_{j\in J}(d_j+\cdots +d_n)\geq
 a_k+\sum_{j\in J-\{k\}}(d_j+\cdots +d_n)
$$
pour tout $J\subset\{i+1,\ldots ,n\}$ et tout $k\in J$ et on a
$$
\sum_{j\in J}(d_j+\cdots +d_n)\geq
d_{n-|J|+1}+2d_{n-|J|+2}+\cdots +|J|d_n
$$
pour tout $J\subset\{i+1,\ldots ,n\}$).

\hfill\hfill\cqfd

\rem Preuve de la conjecture {\rm 3.2} pour $n=2$, avec
$c_2=(-1)^{g-1}$
\endrem
Compte tenu du lemme 3.3, pour chaque
$$
e^\vee =(\Omega_X^1\hookrightarrow {\cal L}_2
\twoheadrightarrow {\cal M}_1)
\in {\rm ob}\,{\cal E}_2^\circ ({\bb F}_q)
$$
image de $({\cal A}_2\hookrightarrow {\cal L}_2
\twoheadrightarrow {\cal L}_1;D_2)
\in E_2^\circ$, la trace de l'endomorphisme de Frobenius
g\'eom\'etrique sur la fibre de $W_{L,2}^\circ$ en $e^\vee$ est
$$
w_{L,2}^\circ (e^\vee)=(-1)^{g-1}\kern -2mm
\sum_{(A_1,\alpha_1),(A_2,\alpha_2)}
\psi (\langle e^\vee ,e\rangle )
w_L(A_1+D_2-2\,{\rm Inf}(A_1,A_2),{\rm Inf}(A_1,A_2)),
$$
o\`u $(A_1,\alpha_1)$ parcourt les paires form\'ees d'un
diviseur effectif $A_1$ et d'un isomorphisme
$\alpha_1:{\cal O}_X(A_1)
\buildrel\sim\over\rightarrow {\cal L}_1$, o\`u
$(A_2,\alpha_2)$ parcourt les paires form\'ees
d'un diviseur effectif $A_2\leq D_2$ et d'un \'epimorphisme
$\alpha_2:{\cal O}_X\twoheadrightarrow {\cal O}_X(-A_2)/
{\cal O}_X(-D_2)$ et o\`u $e$ est la somme de
$$
e_1:{\cal O}_X\hookrightarrow
{\cal O}_X(A_1)\buildrel{\scriptstyle\alpha_1\atop
\scriptstyle\sim}\over\longrightarrow {\cal L}_1
$$
et de
$$
e_2:{\cal O}_X\buildrel\alpha_2\over
\twoheadrightarrow {\cal O}_X(-A_2)/{\cal O}_X(-D_2)
\hookrightarrow {\cal O}_X/{\cal O}_X(-D_2)={\cal O}_{D_2}.
$$
On a encore
$$
w_{L,2}^\circ (e^\vee)=(-1)^{g-1}\kern -4mm
\sum_{\scriptstyle (A_1,\alpha_1)\atop\scriptstyle 0\leq
D'_2\leq {\rm Inf}(A_1,D_2)}
\kern -4mm\psi (\langle e_1^\vee ,e_1\rangle )
w_L(A_1+D_2-2D'_2,D'_2)
\kern -3mm\sum_{{\scriptstyle (A_2,\alpha_2)\atop
\scriptstyle D'_2\leq A_2\leq D_2}
\atop\scriptstyle D'_2={\rm Inf}(A_1,A_2)}\kern -3mm
\psi (\langle e_2^\vee ,e_2\rangle )
$$
o\`u on a d\'ecompos\'e $e^\vee$ en $e_1^\vee+e_2^\vee$ avec
$e_1^\vee$ une extension de
${\cal L}_1$ par $\Omega_X^1$ et  $e_2^\vee$ une extension de
${\cal O}_{D_2}$ par $\Omega_X^1$.

Maintenant, pour chaque diviseur effectif $D'_2\leq
{\rm Inf}(A_1,D_2)$,
la somme
$$
\sum_{{\scriptstyle (A_2,\alpha_2)\atop\scriptstyle D'_2\leq
A_2\leq D_2}\atop\scriptstyle D'_2={\rm Inf}(A_1,A_2)}
\kern -3mm\psi (\langle e_2^\vee ,e_2\rangle )
$$
est le produit sur les $x\in |X|$ de sommes locales
$$
S_x=\sum_{\alpha_{2,x}}\psi (\langle e_{2,x}^\vee ,
e_{2,x}\rangle ),
$$
o\`u $\alpha_{2,x}$ parcourt les homomorphismes de
${\cal O}_{X,x}$-modules
${\cal O}_{X,x}\rightarrow {\cal O}_{X,x}\varpi_x^{d'_{2,x}}/
{\cal O}_{X,x}\varpi_x^{d_{2,x}}$
qui sont surjectifs si
$a_{1,x}>d'_{2,x}$ et arbitraires si $a_{1,x}=d'_{2,x}$. Bien entendu,
on a pos\'e $A_1=\sum_{x\in |X|}a_{1,x}x$, ..., $\varpi_x$ est une
uniformisante de l'anneau de
valuation discr\`ete ${\cal O}_{X,x}$,
$e_{2,x}:{\cal O}_X\rightarrow {\cal O}_{X,x}/
{\cal O}_{X,x}\varpi_x^{d_{2,x}}$
est le compos\'e de $\alpha_{2,x}$ et de
l'inclusion ${\cal O}_{X,x}\varpi_x^{d'_{2,x}}/{\cal O}_{X,x}
\varpi_x^{d_{2,x}}\subset
{\cal O}_{X,x}/{\cal O}_{X,x}\varpi_x^{d_{2,x}}$ et $e_{2,x}^\vee$
est l'extension de
${\cal O}_{X,x}/{\cal O}_{X,x}\varpi_x^{d_{2,x}}$ par $\Omega_X^1$
induite par $e_2^\vee$.
Or, comme ${\cal L}_2$ est sans torsion, il en est de m\^eme de
l'extension de ${\cal O}_{D_2}$
par $\Omega_X^1$ de classe $e_2^\vee$ et la restriction de
$e_{2,x}^\vee$ \`a ${\cal O}_{X,x}\varpi_x^\delta/
{\cal O}_{X,x}\varpi_x^{d_{2,x}}\subset
{\cal O}_{X,x}/{\cal O}_{X,x}\varpi_x^{d_{2,x}}$ n'est pas scind\'ee
quels que soient
$x\in |X|$ et $\delta\in\{0,1,\ldots ,d_{2,x}-1\}$. Par suite, on a
$$
S_x=1
$$
si $d'_{2,x}=d_{2,x}$,
$$
S_x=-1
$$
si $a_{1,x}>d'_{2,x}$ et $d'_{2,x}=d_{2,x}-1$ et
$$
S_x=0
$$
dans tous les autres cas. Il s'en suit que
$$
\sum_{{\scriptstyle (A_2,\alpha_2)\atop
\scriptstyle D'_2\leq A_2\leq D_2}
\atop\scriptstyle D'_2={\rm Inf}(A_1,A_2)}\kern -3mm
\psi (\langle e_2^\vee ,e_2\rangle )
=(-1)^{{\rm deg}(D_2-D'_2)}
$$
si $A_1\geq D_2$ et $d_{2,x}-1\leq d'_{2,x}\leq d_{2,x}$,
$\forall x\in |X|$, et que
$$
\sum_{{\scriptstyle (A_2,\alpha_2)\atop
\scriptstyle D'_2\leq A_2\leq D_2}
\atop\scriptstyle D'_2={\rm Inf}(A_1,A_2)}\kern -3mm
\psi (\langle e_2^\vee ,e_2\rangle )
=0
$$
sinon.

Si l'on revient \`a l'expression de $w_{L,2}^\circ (e^\vee)$,
on voit que
$$\displaylines{
\qquad w_{L,2}^\circ (e^\vee)=(-1)^{g-1}\kern -1mm
\sum_{\scriptstyle (A_1,\alpha_1)\atop\scriptstyle A_1\geq D_2}
\psi (\langle e_1^\vee ,e_1\rangle )\cdot
\hfill\cr\hfill
\cdot\kern -10mm\sum_{\scriptstyle 0\leq D'_2\leq D_2\atop
\scriptstyle d_{2,x}-1\leq d'_{2,x}\leq d_{2,x},~\forall x}
\kern -10mm
(-1)^{{\rm deg}(D_2-D'_2)}w_L(A_1+D_2-2D'_2,D'_2).\qquad}
$$
Mais, pour chaque $A_1\geq D_2\geq 0$, on a
$$\displaylines{
\quad\sum_{\scriptstyle 0\leq D'_2\leq D_2\atop
\scriptstyle d_{2,x}-1\leq d'_{2,x}\leq d_{2,x},~\forall x}
\kern -5mm (-1)^{{\rm deg}(D_2-D'_2)}w_L(A_1+D_2-2D'_2,D'_2)
\hfill\cr\hfill
=\prod_{\scriptstyle x\in |X|\atop\scriptstyle d_{2,x}\geq 1}
\bigl(w_{L,x}(a_{1,x}-d_{2,x},d_{2,x})
-w_{L,x}(a_{1,x}-d_{2,x}+2,d_{2,x}-1)\bigr)\cdot
\hfill\cr\hfill
\cdot\kern -2mm
\prod_{\scriptstyle x\in |X|\atop\scriptstyle d_{2,x}=0}
w_{L,x}(a_{1,x},0),\qquad}
$$
o\`u, pour chaque $x\in |X|$ et chaque couple $(m_1,m_2)$
d'entiers positifs ou nuls, on a
$$
w_{L,x}(m_1,m_2)=\kern -4mm
\sum_{{\scriptstyle \ell_1,\ell_2\atop
\scriptstyle \ell_1+\ell_2\geq m_1+m_2}\atop\scriptstyle
\ell_1+2\ell_2= m_1+2m_2}\kern -4mm
\widetilde K_{(\ell_1+\ell_2,\ell_2),(m_1+m_2,m_2)}(q_x)
{\rm tr}({\rm Frob}_x,R^{(\ell_1,\ell_2)}L_{\overline x}).
$$
Or, pour toutes partitions $\lambda$ et $\mu$ de longueur $\leq 2$
d'un m\^eme entier, on peut calculer le polyn\^ome de
Kostka-Foulkes $K_{\lambda\mu}(t)$ par la r\`egle de Lascoux et
Sch\"utzenberger (cf. [Ma] Ch. III, \S 6). Il y a
exactement un tableau $T$ de forme $\lambda$ et de poids $\mu$
si $\lambda\geq\mu$ et aucun sinon. Dans le
premier cas, le mot $w(T)$ associ\'e \`a $T$ est
$$
2^{\lambda_1-\mu_1}1^{\mu_1}2^{\lambda_2}
$$
et sa charge $c(T)=c(w(T))$ est
$$
\lambda_2c(12)+(\lambda_1-\mu_1)c(21)+(\mu_1-\mu_2)c(1)
=\lambda_1-\mu_1
$$
(on a $\mu_1\geq\mu_2=(\lambda_1-\mu_1)+\lambda_2$,
$c(12)=c(1)=0$ et $c(21)=1$), de sorte que
$$
K_{\lambda\mu}(t)=t^{\lambda_1-\mu_1}.
$$
On en d\'eduit d'une part que, pour tout $m_2\geq 1$, on a
$$
\widetilde K_{(\ell_1+\ell_2,\ell_2),(m_1+m_2,m_2)}(t)=
\widetilde K_{(\ell_1+\ell_2,\ell_2),(m_1+m_2+1,m_2-1)}(t)
$$
et
$\widetilde K_{(m_1+m_2,m_2),(m_1+m_2,m_2)}(t)=t^{m_2}$,
de sorte que
$$
w_{L,x}(m_1,m_2)-w_{L,x}(m_1+2,m_2-1)
=q_x^{m_2}{\rm tr}({\rm Frob}_x,R^{(m_1,m_2)}L_{\overline x}).
$$
On en d\'eduit d'autre part que l'on a
$$
\widetilde K_{(\ell_1+\ell_2,\ell_2),(m_1,0)}(t)=
\left\{\matrix{1&\hbox{si }\ell_1=m_1\hbox{ et }
\ell_2=0,\cr
\noalign{\smallskip}
0&\hbox{sinon,}\hfill\cr}\right.
$$
de sorte que
$$
w_{L,x}(m_1,0)={\rm tr}({\rm Frob}_x,R^{(m_1,0)}L_{\overline x}).
$$
Comme, par d\'efinition,
$$
w_{L,x}^{\rm top}(m_1,m_2)=q_x^{m_2}{\rm tr}({\rm Frob}_x,
R^{(m_1,m_2)}L_{\overline x}),
$$
on obtient finalement que
$$
w_{L,2}^\circ (e^\vee)=(-1)^{g-1}\kern -1mm
\sum_{\scriptstyle (A_1,\alpha_1)\atop\scriptstyle A_1\geq D_2}
\psi (\langle e_1^\vee ,e_1\rangle )w_L^{\rm top}(A_1-D_2,D_2),
$$
i.e.
$$
w_{L,2}^\circ (e^\vee)=(-1)^{g-1}w_{L,2}^{\rm top}
({\cal A}_2\hookrightarrow
{\cal L}_2\twoheadrightarrow {\cal L}_1;D_2).
$$
\hfill\hfill\cqfd

\vskip 5mm
{\bf 4. Op\'erateurs de Hecke et fonctions trace de Frobenius.}
\vskip 5mm

On a vu dans le premier expos\'e que ``$W_L$ est vecteur propre
de l'op\'erateur de Hecke $T_0$ avec valeur propre $V_L$''.
Traduisons ce r\'esultat en termes de fonctions. Avec des
notations \'evidentes, pour toutes suites finies de diviseurs
effectifs $(D_1,\ldots ,D_s)$ et $(E_1,\ldots ,E_t)$ sur $X$,
d\'efinies sur ${\bb F}_q$, avec $\sum_{j=1}^sj{\rm deg}(D_j)=
m_0$ et $\sum_{k=1}^tk{\rm deg}(E_k)=n_0$, la trace de Frobenius
en $(E_1,\ldots ,E_t,D_1,\ldots ,D_s)$ sur $T_0W_L'$ est \'egale
\`a
$$\displaylines{
\qquad (T_0w_L')(E_1,\ldots ,E_t,D_1,\ldots ,D_s)
\hfill\cr\hfill
\buildrel{\rm dfn}\over{=\!=}
{1\over |{\rm Hom}_{{\cal O}_X}({\cal N}_0,{\cal M}_0)|}
\sum_{(D_1',\ldots ,D_{s'}')}\kern -2mm N(D_1',\ldots ,D_{s'}')
w_L'(D_1',\ldots ,D_{s'}'),\qquad}
$$
o\`u  ${\cal M}_0={\cal O}_{D_1+\cdots +D_s}\oplus \cdots
\oplus {\cal O}_{D_s}$ et ${\cal N}_0={\cal O}_{E_1+\cdots +E_t}
\oplus\cdots\oplus {\cal O}_{E_t}$,  o\`u $(D_1',\ldots ,D_{s'}')$
parcourt les suites finies de diviseurs effectifs sur $X$, d\'efinis
sur ${\bb F}_q$, avec $\sum_{j'=1}^{s'}j'{\rm deg}(D_{j'}')=m_0'$,
et o\`u $N(D_1',\ldots ,D_{s'}')$ est le nombre des classes
d'extensions $({\cal M}_0\hookrightarrow {\cal M}_0'
\twoheadrightarrow {\cal N}_0)$  dans ${\rm Ext}_{{\cal O}_X}^1
({\cal N}_0,{\cal M}_0)$ telles que ${\cal M}_0'$ soit isomorphe \`a
${\cal O}_{D_1'+\cdots +D_{s'}'}\oplus \cdots\oplus
{\cal O}_{D_{s}'}$. On a donc
$$
(T_0w_L')(E_1,\ldots ,E_t,D_1,\ldots ,D_s)=
v_L(E_1,\ldots ,E_t)
w_L(D_1,\ldots ,D_s).\leqno (4.1)
$$

Donnons une d\'emonstration directe de cette formule.
Nous aurons besoin du lemme suivant, o\`u on utilise la
terminologie et les notations de [Ma].

\th LEMME 4.2
\enonce
Pour toute partition $\lambda$, soit
$$
Q^\lambda (t)=\sum_{\lambda'\geq\lambda}
\widetilde K_{\lambda',\lambda}(t)\chi^{\lambda'}
$$
le polyn\^ome de Green correspondant,
$\chi^{\lambda'}$ \'etant le caract\`ere de la repr\'esentation
irr\'educ-\break tible de ${\goth S}_{|\lambda'|}$ correspondant
\`a la partition $\lambda'$ de $|\lambda'|=\sum_i\lambda_i'$.

Fixons des entiers $m,n\geq 0$ et identifions
${\goth S}_m\times{\goth S}_n$ au sous-groupe de
${\goth S}_{m+n}$ form\'e
des permutations de $\{1,\ldots ,m+n\}$ qui respectent les
sous-ensembles $\{1,\ldots ,m\}$ et $\{m+1,\ldots ,m+n\}$.
Alors, pour toute partition $\lambda$ de $m+n$, on a
$$
{\rm Res}_{{\goth S}_{m+n}}^{{\goth S}_m\times{\goth S}_n}
Q^\lambda (t)=
\sum_{\scriptstyle\mu ,\nu\atop\scriptstyle |\mu |=m,|\nu |=n}
\kern -1mm g_{\mu\nu}^\lambda (t)(Q^\mu (t)\times Q^\nu (t)),
$$
o\`u les $g_{\mu\nu}^\lambda (t)$ sont les polyn\^omes de Hall,
et, pour toutes partitions $\mu$ et $\nu$ de $m$ et
$n$ respectivement, on a
$$
{\rm Ind}_{{\goth S}_m\times{\goth S}_n}^{{\goth S}_{m+n}}
(Q^\mu (t)\times Q^\nu (t))=
\sum_{\scriptstyle\lambda\atop\scriptstyle |\lambda |=m+n}
\kern -1mm g_\lambda^{\mu\nu}(t)Q^{\lambda}(t),
$$
o\`u on a pos\'e
$$\displaylines{
\qquad g_\lambda^{\mu\nu}(t)=t^{2n(\mu )+2n(\nu )-2n(\lambda )}
g_{\mu\nu}^\lambda (t)\cdot
\hfill\cr\hfill
\cdot\prod_{i\geq 1}{(1-t^{-1})\cdots (1-t^{-m_i(\mu )})
(1-t^{-1})\cdots (1-t^{-m_i(\nu )})\over
(1-t^{-1})\cdots (1-t^{-m_i(\lambda )})}.\qquad}
$$
\endth

\rem Preuve
\endrem
Pour toutes partitions $\mu'$ et $\nu'$ de $m$ et $n$
respectivement, les polyn\^omes de Schur
$s_\bullet (x_1,x_2,\ldots )$ satisfont la relation
(cf. [Ma] Ch. I, \S 9)
$$
s_{\mu'}(x_1,x_2,\ldots )s_{\nu'}(x_1,x_2,\ldots )
=\sum_{\scriptstyle\lambda'\atop
\scriptstyle |\lambda'|=|\mu'|+|\nu'|}c_{\mu'\nu'}^{\lambda'}
s_{\lambda'}(x_1,x_2,\ldots )
$$
o\`u les constantes $c_{\bullet\bullet}^\bullet$ sont d\'efinies
par
$$
{\rm Res}_{{\goth S}_{m+n}}^{{\goth S}_m\times{\goth S}_n}
\chi^{\lambda'}=\sum_{\scriptstyle\mu',\nu'\atop
\scriptstyle |\mu'|=m,|\nu'|=n}c_{\mu'\nu'}^{\lambda'}
(\chi^{\mu'}\times\chi^{\nu'})\qquad
(\forall\lambda'\hbox{ avec }|\lambda'|=m+n)
$$
et, pour toutes partitions $\mu$ et $\nu$ de $m$ et
$n$ respectivement, les polyn\^omes de Hall-Littlewood
$P_\bullet (x_1,x_2,\ldots ;t)$ satisfont la relation
(cf. [Ma] Ch. III, \S 3)
$$
P_\mu (x_1,x_2,\ldots ;t)P_\nu (x_1,x_2,\ldots ;t)=
\sum_{\scriptstyle\lambda\atop
\scriptstyle|\lambda|=|\mu|+|\nu|}f_{\mu\nu}^\lambda (t)
P_\lambda (x_1,x_2,\ldots ;t),
$$
o\`u on a pos\'e
$$
f_{\mu\nu}^\lambda (t)=t^{n(\lambda )-n(\mu )-n(\nu )}
g_{\mu\nu}^\lambda (t^{-1}).
$$
Comme les polyn\^omes de Kostka-Foulkes $K_{\bullet\bullet}(t)$
sont d\'efinis par les relations
$$
s_{\lambda'}(x_1,x_2,\ldots )=\sum_{\lambda\leq\lambda'}
K_{\lambda'\lambda}(t)P_\lambda (x_1,x_2,\ldots ;t),
$$
la premi\`ere formule du lemme s'en suit.

Consid\'erons maintenant,
pour chaque entier $\ell\geq 0$,
le ${\bb Z}$-module libre $R^\ell$
engendr\'e par les caract\`eres irr\'eductibles de ${\goth S}_\ell$
et introduisons le produit scalaire sur $R^\ell$ \`a
valeurs dans ${\bb Q}[t]$, d\'efini par
$$
(f,g)_t=\sum_{\scriptstyle\lambda\atop\scriptstyle
|\lambda |=\ell}z_\lambda (t)^{-1}f_\lambda g_\lambda,
$$
o\`u, pour chaque partition $\lambda$ de $\ell$,
$f_\lambda$ (resp. $g_\lambda$) est la valeur de $f\in R^\ell$
(resp. $g\in R^\ell$) sur la classe de conjugaison de l'\'el\'ement
$$
(1,\ldots ,\lambda_1)(\lambda_1+1,\ldots ,\lambda_1+\lambda_2)
\cdots\in {\goth S}_\ell
$$
et o\`u
$$
z_\lambda (t)=
\prod_{i\geq 1}{i^{m_i(\lambda )}\cdot m_i(\lambda )!\over
(1-t^i)^{m_i(\lambda )}}
$$
(comme dans [Ma], on a not\'e $m_i(\lambda )$ le nombre des
entiers $j\geq 1$ tels que $\lambda_j=i$). Si $m$ et $n$
sont des entiers $\geq 0$, $R^m\otimes R^n$ est le
${\bb Z}$-module libre engendr\'e par les caract\`eres
irr\'eductibles de ${\goth S}_m\times {\goth S}_n$ et
on note encore $(\cdot ,\cdot )_t$
le produit scalaire sur $R^m\otimes R^n$ d\'efini par
$$
(f,g)_t=\sum_{\scriptstyle\mu ,\nu\atop\scriptstyle
|\mu |=m,|\nu |=n}z_\mu (t)^{-1}z_\nu (t)^{-1}
f_{\mu ,\nu} g_{\mu ,\nu} .
$$
On a alors
$$
(f',{\rm Ind}_{{\goth S}_m\times{\goth S}_n}^{{\goth S}_{m+n}}g)_t
=
({\rm Res}_{{\goth S}_{m+n}}^{{\goth S}_m\times{\goth S}_n}f',g)_t
$$
pour tous $f'\in R^{m+n}$ et $g\in R^m\otimes R^n$.
En effet, si on pose
$$
{\rm Res}_{{\goth S}_{m+n}}^{{\goth S}_m\times{\goth S}_n}f'=f
$$
et
$$
{\rm Ind}_{{\goth S}_m\times{\goth S}_n}^{{\goth S}_{m+n}}g=g'
$$
et si, pour toutes partitions $\mu$ et $\nu$ de $m$ et $n$
respectivement, on note $\lambda (\mu ,\nu )$ la partition de
$m+n$ d\'efinie par
$$
m_i(\lambda (\mu ,\nu )) =m_i(\mu )+m_i(\nu )\qquad
(\forall i\geq 1),
$$
on a
$$
f_{\mu ,\nu} =f'_{\lambda (\mu ,\nu )},
$$
$$
g_\lambda' =\left\{\matrix{\displaystyle\prod_{i\geq 1}
\pmatrix{m_i(\mu )+m_i(\nu)\cr m_i(\mu )\cr}g_{\mu ,\nu}
\hfill&\hbox{si }
\lambda =\lambda (\mu ,\nu )\hbox{ pour une paire }(\mu ,\nu ),\cr
\noalign{\medskip}
0\hfill&\hbox{sinon}\hfill\cr}\right.
$$
et
$$
z_{\lambda (\mu ,\nu )}(t)=z_\mu (t)z_\nu (t)
\prod_{i\geq 1}\pmatrix{m_i(\mu )+m_i(\nu)\cr m_i(\mu )\cr}.
$$
Bien s\^ur, cette propri\'et\'e reste conserv\'ee si on \'etend le
produit scalaire ci dessus et les op\'erations
${\rm Res}_{{\goth S}_{m+n}}^{{\goth S}_m\times{\goth S}_n}$ et
${\rm Ind}_{{\goth S}_m\times{\goth S}_n}^{{\goth S}_{m+n}}$
par ${\bb Z}[t]$-lin\'earit\'e aux ${\bb Z}[t]$-modules libres
$R^\ell [t]$.
Or, si on pose
$$
X^\lambda (t)=t^{n(\lambda )}Q^\lambda (t^{-1}),
$$
on a la relation d'orthogonalit\'e
$$
(X^\lambda (t),X^{\lambda'} (t))_t=
\delta_\lambda^{\lambda'}\prod_{i\geq 1}(1-t)
\cdots (1-t^{m_i(\lambda )})
$$
pour toutes partitions $\lambda$, $\lambda'$
(cf. [Ma] Ch. III, $(7.3)'$). De cette
relation d'orthogonalit\'e et de la relation
$$
(X^\lambda ,
{\rm Ind}_{{\goth S}_m\times{\goth S}_n}^{{\goth S}_{m+n}}
(X^\mu\times X^\nu ))_t=
({\rm Res}_{{\goth S}_{m+n}}^{{\goth S}_m\times{\goth S}_n}
X^\lambda ,X^\mu\times X^\nu )_t,
$$
on d\'eduit facilement que la seconde formule du lemme est
\'equivalente \`a la premi\`ere.

\hfill\hfill\cqfd

\rem Preuve directe de la formule $(4.1)$
\endrem
D\'ecomposant chaque $D_j$ en $\sum_{x\in |X|}d_{j,x}x$
et chaque $E_k$ en $\sum_{x\in |X|}e_{k,x}x$,
on v\'erifie imm\'ediatement que
$$
(T_0w_L')(E_1,\ldots ,E_t,D_1,\ldots ,D_s)=
\prod_x\sum_{(d_1',\ldots ,d_{s'}')}\kern -2mm
N_x(d_1',\ldots ,d_{s'}')w_{L,x}(d_1',\ldots ,d_{s'}'),
$$
o\`u $(d_1',\ldots ,d_{s'}')$ parcourt les suites finies d'entiers
$\geq 0$ et o\`u $N_x(d_1',\ldots ,d_{s'}')$ est le nombre
des classes d'extensions
$$
{\cal O}_{(d_{1,x}+\cdots +d_{s ,x})x}\oplus\cdots\oplus
{\cal O}_{d_{s,x}x}\hookrightarrow M_x'
\twoheadrightarrow {\cal O}_{(e_{1,x}+\cdots +e_{t,x})x}
\oplus\cdots\oplus {\cal O}_{e_{t,x}x}
$$
dans
${\rm Ext}_{{\cal O}_{X,x}}^1\bigl(
{\cal O}_{(e_{1,x}+\cdots +e_{t,x})x}\oplus\cdots\oplus
{\cal O}_{e_{t,x}x},{\cal O}_{(d_{1,x}+\cdots +d_{s,x})x}\oplus
\cdots\oplus {\cal O}_{d_{s,x}x}\bigr)$
telles que le ${\cal O}_{X,x}$-module
$M_x'$ soit isomorphe \`a
${\cal O}_{(d_1'+\cdots +d_{s'}')x}\oplus\cdots
\oplus {\cal O}_{d_{s'}'x}$.

Maintenant fixons $x\in |X|$ et posons
$$
M_0={\cal O}_{(d_{1,x}+\cdots +d_{s,x})x}\oplus\cdots\oplus
{\cal O}_{d_{s,x}x},
$$
$$
N_0={\cal O}_{(e_{1,x}+\cdots +e_{t,x})x}\oplus\cdots\oplus
{\cal O}_{e_{t,x}x},
$$
$$
\lambda =(d_1'+\cdots +d_{s'}',\ldots ,d_{s'}'),
$$
$$
\mu =(d_{1,x}+\cdots +d_{s,x},\ldots ,d_{s,x})
$$
et
$$
\nu =(e_{1,x}+\cdots +e_{t ,x},\ldots ,e_{t,x}).
$$
Alors, d'apr\`es [Ma] Ch. II, (4.3),
le nombre des sous-${\cal O}_{X,x}$-modules $M$ de
$$
M_0'={\cal O}_{(d_1'+\cdots +d_{s'}')x}
\oplus\cdots\oplus {\cal O}_{d_{s'}'x}
$$
qui sont isomorphes \`a $M_0$ et tels que $M_0'/M$
soit isomorphe \`a $N_0$ est \'egal \`a
$$
g_{\nu ,\mu}^\lambda (q_x)=g_{\mu ,\nu}^\lambda (q_x).
$$
Par suite, le nombre $N_x(d_1',\ldots ,d_{s'}')$
des classes d'extensions $(M_0\hookrightarrow M'
\twoheadrightarrow N_0)$ telles que $M'$ soit isomorphe
\`a $M_0'$ est \'egal \`a
$$
{|{\rm Hom}_{{\cal O}_{X,x}}(N_0,M_0)|\over
|{\rm Aut}_{{\cal O}_{X,x}}(M_0')|}
g_{\mu ,\nu}^\lambda (q_x)|{\rm Aut}_{{\cal O}_{X,x}}(M_0)|
\,|{\rm Aut}_{{\cal O}_{X,x}}(N_0)|.
$$
En effet, on a une bijection naturelle entre cet ensemble de
classes d'extensions et l'ensemble quotient
$$
{\rm Aut}_{{\cal O}_{X,x}}(M'_0)\big\backslash
\coprod_{M\subset M'_0}
\bigl({\rm Isom}_{{\cal O}_{X,x}}(M_0,M)\times
{\rm Isom}_{{\cal O}_{X,x}}(M'_0/M,N_0)\bigr),
$$
le stabilisateur dans ${\rm Aut}_{{\cal O}_{X,x}}(M'_0)$
d'un point arbitraire $(M,\alpha ,\beta)$ dans l'ensemble
$\coprod_{M\subset M'_0}\bigl({\rm Isom}_{{\cal O}_{X,x}}(M_0,M)
\times {\rm Isom}_{{\cal O}_{X,x}}(M'_0/M,N_0)\bigr)$
n'est autre que
$$
{\rm Id}_{M_0'}+u\circ {\rm Hom}_{{\cal O}_{X,x}}(M_0'/M,M)
\circ v,
$$
o\`u $u:M\hookrightarrow M_0'$ est l'inclusion et $v:M_0'
\twoheadrightarrow M_0'/M$ est l'application quotient, et ce
stabilisateur a donc $|{\rm Hom}_{{\cal O}_{X,x}}(N_0,M_0)|$
\'el\'ements.  Comme, pour toute partition $\pi =(\pi_1,\pi_2,
\ldots )$, $\pi_1\geq\pi_2\geq \cdots\geq 0$,  on a
$$
|{\rm Aut}_{{\cal O}_{X,x}}({\cal O}_{\pi_1x}
\oplus{\cal O}_{\pi_2x}\oplus\cdots )|
=a_\pi (q_x),
$$
o\`u
$$
a_\pi (t)=t^{|\pi |+2n(\pi )}\prod_{i\geq 1}(1-t^{-1})\cdots
(1-t^{-m_i(\pi )})
$$
(cf. [Ma] Ch. II, (1.6)), on obtient que
$N_x(d_1',\ldots ,d_{s'}')$ est \'egal \`a
$$
{|{\rm Hom}_{{\cal O}_{X,x}}(N_0,M_0)|\over a_\lambda (q_x)}
g_{\mu ,\nu}^\lambda (q_x)a_\mu (q_x)a_\nu (q_x).
$$

Par suite, l'expression
$$
{1\over |{\rm Hom}_{{\cal O}_{X,x}}(N_0,M_0)|}
\sum_{(d_1',\ldots ,d_{s'}')}
N_x(d_1',\ldots ,d_{s'}')
w_{L,x}(d_1',\ldots ,d_{s'}')
$$
est \'egale \`a
$$
\sum_\lambda
{g_{\mu ,\nu}^\lambda (q_x)
a_\mu (q_x)a_\nu (q_x)\over a_{\lambda}(q_x)}
w_{L,x}(\lambda_1-\lambda_2,\lambda_2-\lambda_3,\ldots ,
\lambda_{s'}),
$$
i.e. \`a
$$
\sum_\lambda g_\lambda^{\mu ,\nu}(q_x)
w_{L,x}(\lambda_1-\lambda_2,\lambda_2-\lambda_3,\ldots ,
\lambda_{s'}),
\leqno (*)
$$
puisque $g_{\mu\nu}^\lambda (t)\equiv 0$ si
$|\lambda |\not=|\mu |+|\nu |$.

Mais
$$
w_{L,x}(\lambda_1-\lambda_2,\lambda_2-\lambda_3,\ldots ,
\lambda_{s'})=\sum_{\lambda'\geq\lambda}
\widetilde K_{\lambda'\lambda}(q_x){\rm tr}({\rm Frob}_x,
R^{(\lambda_1'-\lambda_2',\lambda_2'-\lambda'_3,\ldots )}
L_{\overline x})
$$
est la trace de ${\rm Frob}_x$ sur
$$
\bigl(Q^\lambda (q_x)\otimes
L_{\overline x}^{\otimes |\lambda |}\bigr)^{{\goth S}_{|\lambda |}}
$$
et donc, compte tenu de la deuxi\`eme formule du lemme 4.2,
l'expression $(*)$ est la trace de ${\rm Frob}_x$ sur
$$
\bigl({\rm Ind}_{{\goth S}_{|\mu|}\times
{\goth S}_{|\nu|}}^{{\goth S}_{|\mu |+|\nu |}}
(Q^\mu (q_x)\times Q^\nu (q_x))\otimes
L_{\overline x}^{\otimes (|\mu |+|\nu |)}
\bigr)^{{\goth S}_{|\mu |+|\nu |}}
$$
ou, ce qui revient au m\^eme, sur
$$
\bigl(Q^\mu (q_x)\otimes
L_{\overline x}^{\otimes |\mu |}\bigr)^{{\goth S}_{|\mu |}}\otimes
\bigl(Q^\nu (q_x)\otimes
L_{\overline x}^{\otimes |\nu |}\bigr)^{{\goth S}_{|\nu |}}.
$$
Ceci ach\`eve la preuve de la formule (4.1).

\hfill\hfill\cqfd

Pr\'ecisons le lemme 4.2 du premier expos\'e en tenant compte
des torsions \`a la Tate.

\th LEMME 4.3
\enonce
{\rm (i)} Pour tout $K'\in {\rm ob}\,D_{\rm c}^{\rm b}({\cal C}_0',
\overline{\bb Q}_\ell )$, on a un isomorphisme canonique
$$
T_1^\circ (j_0'^{\vee *}\pi_0'^{\vee *}K'[m_0])(n_0)\cong
({\rm id}\times j_0^\vee )^*({\rm id}\times\pi_0^\vee )^*
T_0 (K')[m_0]
$$
dans $D_{\rm c}^{\rm b}({\cal C}{\it oh}_X^{0,n_0}\times
{\cal E}_1^\circ, \overline{\bb Q}_\ell )$.

\decale{\rm (ii)} Pour tout
$K'\in {\rm ob}\,D_{\rm c}^{\rm b}({\cal E}_{i-1}'^\circ,
\overline{\bb Q}_\ell )$, on a un isomorphisme canonique
$$
T_i^\circ \bigl(j_{i-1}'^{\vee *}
{\cal F}_{\psi ,i-1}'(j_{i-1,!}'K')\bigr)(n_0)\cong
({\rm id}\times j_{i-1}^\vee )^*
{\cal F}_{\psi ,i-1,{\cal C}oh_X^{0,n_0}}
\bigl(({\rm id}\times j_{i-1})_!T_{i-1}^\circ (K')\bigr),
$$
dans  $D_{\rm c}^{\rm b}({\cal C}{\it oh}_X^{0,n_0}\times
{\cal E}_i^\circ,\overline{\bb Q}_\ell )$, o\`u ${\cal F}_{\psi ,i-1}'$ et
${\cal F}_{\psi ,i-1,{\cal C}oh_X^{0,n_0}}$ sont les transformations
de Fourier pour les fibr\'es $\pi_{i-1}':{\cal E}_{i-1}'\rightarrow
{\cal C}_{i-1}'$ et ${\rm id}\times\pi_{i-1}:
{\cal C}{\it oh}_X^{0,n_0}\times {\cal E}_{i-1}\rightarrow
{\cal C}{\it oh}_X^{0,n_0}\times {\cal C}_{i-1}$ respectivement.

\decale{\rm (iii)} Pour tout
$K'\in {\rm ob}\,D_{\rm c}^{\rm b}({\cal C}_n',\overline{\bb Q}_\ell )$,
on a un isomorphisme canonique
$$
T_n^\circ (j_n'^*\pi_n'^*K')\cong
({\rm id}\times j_n)^*({\rm id}\times\pi_n)^*T_n(K')
$$
dans $D_{\rm c}^{\rm b}({\cal C}{\it oh}_X^{0,n_0}\times
{\cal E}_n^\circ,\overline{\bb Q}_\ell )$.
\endth

\rem Preuve
\endrem
On raisonnera en termes de fonctions et on utilisera des
notations \'evidentes pour les op\'erateurs de Hecke, les
transformations de Fourier, ..., au niveau des fonctions.

Pour toute fonction $\varphi'$ sur ${\cal C}_0'({\bb F}_q)$,
tout ${\cal N}_0\in {\rm ob}\,{\cal C}{\it oh}_X^{0,n_0}({\bb F}_q)$
et tout monomorphisme $e_1=({\cal O}_X\hookrightarrow
{\cal M}_1)$ dans ${\rm ob}\,{\cal E}_1^\circ ({\bb F}_q)$, on a
$$
T_1^\circ (j_0'^{\vee *}\pi_0'^{\vee *}
\varphi'[m_0])(n_0)({\cal N}_0,e)
={(-1)^{m_0}\over q^{n_0}|{\rm Hom}_{{\cal O}_X}({\cal N}_0,
{\cal M}_1)|}\sum_{u_1}\varphi' ({\cal M}_0'),
$$
o\`u $u_1=({\cal M}_1\buildrel\alpha_1\over\hookrightarrow
{\cal M}_1'\twoheadrightarrow {\cal N}_0)$
parcourt ${\rm Ext}_{{\cal O}_X}^1({\cal N}_0,{\cal M}_1)$ et o\`u
${\cal M}_0'$ est le conoyau de $\alpha_1\circ e_1:{\cal O}_X
\hookrightarrow {\cal M}_1'$,
et on a
$$
({\rm id}\times j_0^\vee )^*({\rm id}\times\pi_0^\vee )^*
T_0 (\varphi')[m_0]({\cal N}_0,e)=
{(-1)^{m_0}\over |{\rm Hom}_{{\cal O}_X}({\cal N}_0,{\cal M}_0)|}
\sum_{u_0}\varphi' ({\cal M}_0'),
$$
o\`u $u_0=({\cal M}_0\buildrel\alpha_0\over\hookrightarrow
{\cal M}_0'\twoheadrightarrow {\cal N}_0)$ parcourt
${\rm Ext}_{{\cal O}_X}^1({\cal N}_0,{\cal M}_0)$, ${\cal M}_0$
\'etant le conoyau de $e_1:{\cal O}_X\hookrightarrow {\cal M}_1$.

La partie (i) du lemme en d\'ecoule puisque l'on a la suite exacte
longue
$$\matrix{
&0&\rightarrow &{\rm Hom}_{{\cal O}_X}({\cal N}_0,{\cal M}_1)&
\rightarrow& {\rm Hom}_{{\cal O}_X}({\cal N}_0,{\cal M}_0)&
\rightarrow &\cr
\noalign{\medskip}
\rightarrow &{\rm Ext}_{{\cal O}_X}^1({\cal N}_0,{\cal O}_X)&
\rightarrow &{\rm Ext}_{{\cal O}_X}^1({\cal N}_0,{\cal M}_1)&
\rightarrow& {\rm Ext}_{{\cal O}_X}^1({\cal N}_0,{\cal M}_0)&
\rightarrow& 0,\cr}
$$
avec ${\rm dim}_{{\bb F}_q}
{\rm Ext}_{{\cal O}_X}^1({\cal N}_0,{\cal O}_X)=n_0$.

Pour toute fonction $\varphi'$ sur ${\cal E}_{i-1}'^\circ({\bb F}_q)$,
tout ${\cal N}_0\in {\rm ob}\,{\cal C}{\it oh}_X^{0,n_0}({\bb F}_q)$
et tout monomorphisme $e_i=((\Omega_X^1)^{\otimes (i-1)}
\hookrightarrow {\cal M}_i)$ dans ${\rm ob}\,{\cal E}_i^\circ
({\bb F}_q)$, on a
$$\displaylines{
\qquad T_i^\circ \bigl(j_{i-1}'^{\vee *}
{\cal F}_{\psi ,i-1}'(j_{i-1,!}'\varphi')\bigr)(n_0)({\cal N}_0,e_i)
\hfill\cr\hfill
={(-1)^{in_0+m_0'-(i-1)^2(g-1)}\over
q^{n_0}|{\rm Hom}_{{\cal O}_X}({\cal N}_0,{\cal M}_i)|}
\sum_{u_i}\sum_{e_{i-1}'}\psi (\langle j_{i-1}'^\vee
(\alpha_i\circ e_i), e_{i-1}'\rangle)\varphi' (e_{i-1}'),\qquad}
$$
o\`u $u_i=({\cal M}_i\buildrel\alpha_i\over\hookrightarrow
{\cal M}_i'\twoheadrightarrow {\cal N}_0)$ parcourt
${\rm Ext}_{{\cal O}_X}^1({\cal N}_0,{\cal M}_i)$ et o\`u
$e_{i-1}'$ parcourt l'ensemble des monomorphismes de
$(\Omega_X^1)^{\otimes (i-2)}$ dans le conoyau
${\cal M}_{i-1}'$ de $\alpha_i\circ e_i$, et on a
$$\displaylines{
\qquad ({\rm id}\times j_{i-1}^\vee )^*
{\cal F}_{\psi ,i-1,{\cal C}oh_X^{0,n_0}}\bigr(({\rm id}
\times j_{i-1})_!T_{i-1}^\circ (\varphi')\bigl)({\cal N}_0,e_i)
\hfill\cr\hfill
=(-1)^{m_0-(i-1)^2(g-1)}
\sum_{e_{i-1}}\psi (\langle j_{i-1}^\vee (e_i),e_{i-1}\rangle)\cdot
\hfill\cr\hfill
\cdot {(-1)^{(i-1)n_0}\over |{\rm Hom}_{{\cal O}_X}({\cal N}_0,
{\cal M}_{i-1})|}
\sum_{u_{i-1}}\varphi' (\alpha_{i-1}\circ e_{i-1}),\qquad}
$$
o\`u $e_{i-1}$ parcourt l'ensemble des monomorphismes de
$(\Omega_X^1)^{\otimes (i-2)}$ dans le conoyau ${\cal M}_{i-1}$
de $e_i$ et o\`u $u_{i-1}=({\cal M}_{i-1}\buildrel\alpha_{i-1}
\over\hookrightarrow {\cal M}_{i-1}'\twoheadrightarrow
{\cal N}_0)$ parcourt ${\rm Ext}_{{\cal O}_X}^1({\cal N}_0,
{\cal M}_{i-1})$.

Or, d'une part, on a la suite exacte
longue
$$\matrix{
0&\rightarrow &{\rm Hom}_{{\cal O}_X}({\cal N}_0,{\cal M}_i)&
\rightarrow& {\rm Hom}_{{\cal O}_X}({\cal N}_0,{\cal M}_{i-1})&
\rightarrow &\cr
\noalign{\medskip}
{\rm Ext}_{{\cal O}_X}^1({\cal N}_0,(\Omega_X^1)^{\otimes (i-1)})
&\rightarrow &{\rm Ext}_{{\cal O}_X}^1({\cal N}_0,{\cal M}_i)&
\rightarrow& {\rm Ext}_{{\cal O}_X}^1({\cal N}_0,{\cal M}_{i-1})&
\rightarrow& 0\cr}
$$
avec ${\rm dim}\,{\rm Ext}_{{\cal O}_X}^1({\cal N}_0,
(\Omega_X^1)^{\otimes (i-1)})=n_0$. D'autre part,
pour chaque $u_{i-1}$, on a la suite exacte courte
$$\displaylines{
\qquad 0\rightarrow
{\rm Hom}_{{\cal O}_X}((\Omega_X^1)^{\otimes (i-2)},
{\cal M}_{i-1})\rightarrow
{\rm Hom}_{{\cal O}_X}((\Omega_X^1)^{\otimes (i-2)},
{\cal M}_{i-1}')
\hfill\cr\hfill
\rightarrow {\rm Hom}_{{\cal O}_X}((\Omega_X^1)^{\otimes (i-2)},
{\cal N}_0)\rightarrow 0,\qquad}
$$
puisque
${\rm Hom}_{{\cal O}_X}({\cal M}_{i-1},
(\Omega_X^1)^{\otimes (i-1)})=(0)$. Enfin, on a
$$
\langle j_{i-1}'^\vee (\alpha_i\circ e_i),
\alpha_{i-1}\circ e_{i-1}\rangle
=\langle j_{i-1}^\vee (e_i),e_{i-1}\rangle
$$
pour chaque diagramme commutatif de ${\cal O}_X$-Modules
coh\'erents
$$
\def\normalbaselines{\baselineskip20pt\lineskip3pt
\lineskiplimit3pt}
\matrix{(\Omega_X^1)^{\otimes (i-1)}&\buildrel{e_i}\over
\hookrightarrow&{\cal M}_i&\twoheadrightarrow &{\cal M}_{i-1}\cr
\Big|\!\Big|&&\mapdownleft{\alpha_i}&\square&
\mapdownright{\alpha_{i-1}}\cr
(\Omega_X^1)^{\otimes (i-1)}&\hookrightarrow
&{\cal M}_i'&\twoheadrightarrow &{\cal M}_{i-1}'\cr}
$$
et chaque homomorphisme $e_{i-1}:(\Omega_X^1)^{\otimes (i-2)}
\rightarrow {\cal M}_{i-1}$.

La partie (ii) du lemme s'ensuit.

La partie (iii) est triviale.

\hfill\hfill\cqfd

\th COROLLAIRE 4.4
\enonce
Pour chaque $i=1,\ldots ,n$, on a
$$
T_i^\circ (w_{L,i}'^\circ )=v_L
\boxtimes w_{L,i}^\circ
$$
sur ${\rm ob}\,{\cal C}{\it oh}_X^{0,n_0}({\bb F}_q)\times
{\rm ob}\,{\cal E}_i({\bb F}_q)$,
o\`u $w_{L,i}^\circ$ et $w_{L,i}'^\circ$ sont les fonctions traces de
Frobenius sur les fibres des complexes
$W_{L,i}^\circ$ et $W_{L,i}'^\circ$ d\'efinis comme dans la section
$3$ du premier expos\'e.
\endth
\vskip 2mm

{\pc REMARQUE} 4.5\pointir
D'apr\`es le corollaire 4.4, la fonction trace de  Frobenius sur
les fibres de $W_{L,n}^\circ$ est fonction propre des
op\'erateurs de Hecke avec comme valeur propre $v_L$. Il
en est bien s\^ur de m\^eme de sa restriction \`a $E_n^\circ$.
Par suite, on s'attend \`a ce que cette restriction soit
proportionnelle \`a $w_{L,n}^{\rm top}$ qui est vecteur propre
des op\'erateurs de Hecke avec les m\^emes valeurs propres par
construction (cf. [Sh]). Le seul probl\`eme est que la notion
d'op\'erateur de Hecke introduite dans ces expos\'es
n'est pas exactement celle qui est utilis\'ee classiquement.
Dans notre d\'efinition de ${\cal H}{\it ecke}_X^{n,m_n,m_n'}$ nous
permettons des modifications du type $({\cal M}_n,{\cal M}_n',
\alpha_n)$ o\`u ${\cal M}_n$ n'a pas de torsion mais o\`u
${\cal M}_n'$ en a, alors que, dans la traduction g\'eom\'etrique
des op\'erateurs de Hecke utilis\'es par Shalika, ce n'est pas
permis. Pour prouver la conjecture 3.2 en utilisant le corollaire
4.4, il faudrait donc d'abord comparer ces deux notions
d'op\'erateurs de Hecke.

\hfill\hfill\cqfd

\vfill\eject

\centerline{\douzebf Compl\'ements}
\vskip 20mm
{\bf 1. Variante sch\'ematique du diagramme fondamental.}
\vskip 5mm

On munit $X$ d'un fibr\'e en droites tr\`es ample ${\cal O}_X(1)$
et on note $\delta >0$ le degr\'e de ${\cal O}_X(1)$.

Soit $\nu$ un entier tel que $\nu\delta\geq 2g$ et soit
$N=m_0+n((n-2)(g-1)+\nu\delta )$.
Pour chaque $i=0,\ldots ,n$, on consid\`ere le sch\'ema
$$
{\rm Quot}_{{\cal O}_X(-\nu )^N/X}^{i,m_i}
$$
des ${\cal O}_X$-Modules coh\'erents quotients
$$
{\cal O}_X(-\nu )^N\twoheadrightarrow {\cal M}_i
$$
qui sont de rang g\'en\'erique $i$ et de degr\'e $m_i$.
D'apr\`es Grothendieck, c'est un sch\'ema projectif
sur ${\bb C}$ et on a un morphisme (repr\'esentable) de champs
alg\'ebriques
$$
{\rm Quot}_{{\cal O}_X(-\nu )^N/X}^{i,m_i}\rightarrow
 {\cal C}{\it oh}_{X}^{i,m_i},~
({\cal O}_X(-\nu )^N\twoheadrightarrow {\cal M}_i)
\mapsto {\cal M}_i.
$$

Pour chaque $({\cal O}_X(-\nu )^N\twoheadrightarrow {\cal M}_i)
\in {\rm Quot}_{{\cal O}_X(-\nu )^N/X}^{i,m_i}$, le
${\cal O}_X$-Module coh\'erent ${\cal M}_i(\nu )$ est
engendr\'e par ses sections globales. Ceci nous am\`ene
\`a introduire l'ouvert
$$
{\cal C}_{i,\nu}\subset {\cal C}_i\subset{\cal C}{\it oh}_{X}^{i,m_i}
$$
form\'e des ${\cal M}_i\in {\cal C}_i$ tels que ${\cal M}_i(\nu )$
est engendr\'e par ses sections globales.

\th LEMME 1.1
\enonce
On a ${\cal C}_{0,\nu}={\cal C}{\it oh}_{X}^{0,m_0}$ et, pour
tout $i=1,\ldots ,n$, l'ouvert ${\cal C}_{i,\nu}$ de
${\cal C}{\it oh}_{X}^{i,m_i}$ contient l'ouvert
${\cal F}{\it ib}_{X}^{i,m_i,{\rm ss}}$ des
${\cal O}_X$-Modules localement libres semi-stables
de rang $i$ et de degr\'e $m_i$.
\endth

\rem Preuve
\endrem
Compte tenu du lemme 1.1 de l'expos\'e I, il s'agit de voir que,
pour chaque
${\cal L}_i\in {\rm ob}\,{\cal F}{\it ib}_{X}^{i,m_i,{\rm ss}}$, le
${\cal O}_X$-Module ${\cal L}_i(\nu )$
est engendr\'e par ses sections globales. Or, on a
$$
m_i+i\nu\delta >i\nu\delta > i(2g-1)
$$
par hypoth\`ese et donc on peut appliquer le lemme 5.2 de [Ne] au
${\cal O}_X$-Module semi-stable
${\cal L}_i(\nu )$.

\hfill\hfill\cqfd
\vskip 2mm

{\pc REMARQUE} 1.2\pointir
Pour $\nu\geq 2g/\delta$ variable, on a
$$
\cdots\subset {\cal C}_{i,\nu}\subset {\cal C}_{i,\nu +1}
\subset\cdots
$$
et
$$
\bigcup_{\nu\geq 2g/\delta}{\cal C}_{i,\nu}={\cal C}_i.
$$
\hfill\hfill\cqfd

On note
$$
Q_i\subset {\rm Quot}_{{\cal O}_X(-\nu )^N/X}^{i,m_i}
$$
l'ouvert des $({\cal O}_X(-\nu )^N\twoheadrightarrow {\cal M}_i)$
tels que ${\cal M}_i\in {\cal C}_i$.
Pour $({\cal O}_X(-\nu )^N\twoheadrightarrow {\cal M}_i)\in Q_i$,
on a
$$
H^1(X,{\cal M}_i(\nu ))=\bigl({\rm Hom}_{{\cal O}_X}({\cal M}_i,
\Omega_X^1(-\nu ))\bigr)^\vee=(0),
$$
puisque $H^0(X,{\cal O}_X(\nu ))\not= (0)$ ($\nu\delta\geq 2g$)
et que ${\rm Hom}_{{\cal O}_X}({\cal M}_i,\Omega_X^1)=(0)$.
Par suite, les $H^0(X,{\cal M}_i(\nu ))$ s'organisent en un fibr\'e
vectoriel de rang $m_0+i((i-2)(g-1)+n\delta )$ sur $Q_i$. On note
$$
Q_{i,\nu}\subset Q_i
$$
l'ouvert form\'e des $({\cal O}_X(-\nu )^N\twoheadrightarrow
{\cal M}_i)\in Q_i$ tels que l'application lin\'eaire induite
$$
H^0(X,{\cal O}_X)^N\rightarrow H^0(X,{\cal M}_i(\nu ))
$$
soit surjective. Bien entendu, le morphisme de champs
${\rm Quot}_{{\cal O}_X(-\nu)/X}^{i,m_i}\rightarrow
{\cal C}{\it oh}_{X}^{i,m_i}$ induit un morphisme
(repr\'esentable) de champs
$$
Q_{i,\nu}\rightarrow {\cal C}_{i,\nu}.
$$

\th LEMME 1.3
\enonce
Le morphisme
$Q_{i,\nu}\rightarrow {\cal C}_{i,\nu}$ est
quasi-projectif, lisse, surjectif et
\`a fibres connexes de dimension
$\bigl(m_0+n((n-2)(g-1)+\nu\delta )\bigr)
\bigl(m_0+i((i-2)(g-1)+\nu\delta )\bigr)$.
En particulier, le sch\'ema quasi-projectif $Q_{i,\nu}$ est
connexe et lisse de dimension
$$
\bigl(m_0+n((n-2)(g-1)+\nu\delta )\bigr)
\bigl(m_0+i((i-2)(g-1)+\nu\delta )\bigr)+i^2(g-1)
$$
sur ${\bb C}$.
\endth

\rem Preuve
\endrem
Se donner un point de $Q_{i,\nu}$ revient \`a se donner un
${\cal O}_X$-Module coh\'erent
${\cal M}_i$ dans ${\cal C}_{i,\nu}$ et une application
lin\'eaire surjective
$$
{\bb C}^N\twoheadrightarrow H^0(X,{\cal M}_i(\nu )).
$$

\hfill\hfill\cqfd

Pour chaque $i=0,\ldots ,n$, notons
$$
\pi_{i,\nu}:E_{i,\nu}\rightarrow Q_{i,\nu}
$$
le fibr\'e vectoriel
de fibre ${\rm Hom}_{{\cal O}_X}((\Omega_X^1)^{\otimes (i-1)},
{\cal M}_i)$ en
$({\cal O}_X(-\nu )^N\twoheadrightarrow {\cal M}_i)\in Q_{i,\nu}$
et notons
$$
\pi_{i,\nu}^\vee :E_{i,\nu}^\vee\rightarrow Q_{i,\nu}
$$
le fibr\'e vectoriel dual, de fibre ${\rm Ext}_{{\cal O}_X}^1
({\cal M}_i,(\Omega_X^1)^{\otimes i})$
en $({\cal O}_X(-\nu )^N\twoheadrightarrow {\cal M}_i)$.
Les fibr\'es vectoriels $\pi_{i,\nu}$ et
$\pi_{i,\nu}^\vee$ sont tous les deux de rang $m_0-i^2(g-1)$.

Pour chaque $i=1,\ldots ,n$, soit
$$
E_{i,\nu}^\circ\subset E_{i,\nu}
$$
l'ouvert des diagrammes
$$
\matrix{&&{\cal O}_X(-\nu )^N\cr
&&\downarrow\cr
(\Omega_X^1)^{\otimes (i-1)}&\rightarrow\kern -2mm &
{\cal M}_i\cr}
$$
tels que la fl\`eche horizontale soit injective
(d'apr\`es le lemme 1.1, cet ouvert est non vide) et soit
$$
E_{i-1,\nu}^{\vee\circ}\subset E_{i-1,\nu}^\vee
$$
l'ouvert des diagrammes
$$
\matrix{&&&&{\cal O}_X(-\nu )^N\cr
&&&&\downarrow\cr
(\Omega_X^1)^{\otimes (i-1)}&\hookrightarrow
&{\cal M}_i&\twoheadrightarrow\kern -2mm  &{\cal M}_{i-1}\cr}
$$
tels que ${\cal M}_i\in {\cal C}_i$.

On d\'efinit alors un morphisme de sch\'emas
$$
\rho_{i,\nu}:E_{i,\nu}^\circ\rightarrow E_{i-1,\nu}^{\vee\circ}
$$
en envoyant le diagramme
$$
\matrix{&&{\cal O}_X(-\nu )^N\cr
&&\downarrow\cr
(\Omega_X^1)^{\otimes (i-1)}&\hookrightarrow\kern -2mm &
{\cal M}_i\cr}
$$
sur le diagramme
$$
\matrix{&&&&{\cal O}_X(-\nu )^N\cr
&&&&\downarrow\cr
(\Omega_X^1)^{\otimes (i-1)}&\hookrightarrow &{\cal M}_i&
\twoheadrightarrow\kern -2mm  &{\cal M}_{i-1},\cr}
$$
o\`u ${\cal M}_{i-1}$ est le quotient de ${\cal M}_i$ par
$(\Omega_X^1)^{\otimes (i-1)}$ et o\`u la fl\`eche
verticale est le quotient compos\'e des quotients
${\cal O}_X(-\nu )^N\twoheadrightarrow {\cal M}_i$ et ${\cal M}_i
\twoheadrightarrow {\cal M}_{i-1}$. Ce dernier diagramme
est bien dans $E_{i-1,\nu}^{\vee\circ}$ car on a le plongement
$$
{\rm Hom}_{{\cal O}_X}({\cal M}_{i-1},
(\Omega_X^1)^{\otimes (2n-1)})\hookrightarrow
{\rm Hom}_{{\cal O}_X}({\cal M}_i,(\Omega_X^1)^{\otimes (2n-1))})
=(0)
$$
et la suite exacte
$$
H^0(X,{\cal M}_i(\nu ))\rightarrow H^0(X,{\cal M}_{i-1}(\nu ))
\rightarrow
H^1(X,(\Omega_X^1)^{\otimes (i-1)}(\nu ))=(0)
$$
puisque $\nu\delta\geq 2g$.

\th LEMME 1.4
\enonce
Le morphisme $\rho_{i,\nu}:E_{i,\nu}^\circ\rightarrow
 E_{i-1,\nu}^{\vee\circ}$ est lisse,
surjectif et \`a fibres connexes de dimension
$N((2i-3)(g-1)+\nu\delta )$.
\endth

\rem Preuve
\endrem
Soit
$$
\matrix{&&&&{\cal O}_X(-\nu )^N\cr
&&&&\downarrow\cr
(\Omega_X^1)^{\otimes (i-1)}&\hookrightarrow &{\cal M}_i&
\twoheadrightarrow\kern -2mm  &{\cal M}_{i-1}\cr}
$$
un point de $E_{i-1,\nu}^{\vee\circ}$. Comme ${\cal M}_{i-1}(\nu )$
et $(\Omega_X^1)^{\otimes (i-1)}(\nu )$ sont engendr\'es par leurs
sections globales et que la fl\`eche $H^0(X,{\cal M}_i(\nu ))
\rightarrow H^0(X,{\cal M}_{i-1}(\nu ))$ est surjective,
${\cal M}_i(\nu )$ est aussi engendr\'e par ses sections globales.
La fibre de $\rho_{i,\nu}$ en ce point est donc isomorphe
\`a la vari\'et\'e des applications
lin\'eaires surjectives
$$
{\bb C}^N\twoheadrightarrow H^0(X,{\cal M}_i(\nu ))
$$
qui rel\`event l'application lin\'eaire surjective
$$
{\bb C}^N\twoheadrightarrow H^0(X,{\cal M}_{i-1}(\nu ))
$$
induite par
l'\'epimorphisme de ${\cal O}_X$-Modules ${\cal O}_X(-\nu )^N
\twoheadrightarrow
{\cal M}_{i-1}(\nu )$.

\hfill\hfill\cqfd

La variante sch\'ematique du diagramme fondamental est alors
$$\displaylines{
\quad\matrix{&&\kern -3mm E_{n,\nu}\supset E_{n,\nu}^\circ
\twoheadrightarrow E_{n-1,\nu}^{\vee\circ}\subset
E_{n-1,\nu}^\vee\kern -3mm&\cr
&\swarrow&&\searrow \cr
Q_{n,\nu}\kern -3mm&&&\cr}\quad\cdots
\hfill\cr\hfill
\cdots\quad\matrix{&&&\kern -3mm E_{1,\nu}\supset
E_{1,\nu}^\circ\twoheadrightarrow E_{0,\nu}^{\vee\circ}
\subset E_{0,\nu}^\vee\kern -3mm&&\cr
\searrow &&\swarrow&&\searrow &\cr
&\kern -3mm Q_{1,\nu}\kern -3mm &&&&\kern -3mm Q_{0,\nu}.\cr}
\quad}
$$

Ce diagramme de sch\'emas s'envoie dans le diagramme
fondamental de champs de l'expos\'e I
de mani\`ere \'evidente (oubli de l'application quotient
${\cal O}_X(-\nu )^N\twoheadrightarrow$).
En outre, le groupe alg\'ebrique $GL_{N,{\bb C}}$
agit sur ce diagramme de sch\'emas (composition
d'un automorphisme de ${\cal O}_X(-\nu )^N$ avec l'application
quotient ${\cal O}_X(-\nu )^N\twoheadrightarrow$) et le
morphisme du diagramme de sch\'emas vers le diagramme de
champs est invariant pour cette action. Enfin, la preuve du
lemme 1.3 montre que le champ ${\cal C}_{n,\nu}$ est le quotient
du sch\'ema $Q_{n,\nu}$ pour cette action de $GL_{N,{\bb C}}$.

Les constructions de $W_L$ et des $W_{L,i}$, $W_{L,i}^\circ$
admettent des variantes $GL_{N,{\bb C}}$-\'equivarian-\break
tes sur les sch\'emas $Q_{0,\nu}$ et $E_{i,\nu}$, $E_{i,\nu}^\circ$.

\vskip 5mm
{\bf 2. Vari\'et\'es caract\'eristiques.}
\vskip 5mm

Pour $i=0,\ldots ,n$, le fibr\'e cotangent $T^*{\cal C}_i$ au champ
alg\'ebrique lisse ${\cal C}_i$
est le champ des paires $({\cal M}_i,\theta_i)$, o\`u
${\cal M}_i$ est un objet de ${\cal C}_i$ et o\`u
$\theta_i:{\cal M}_i\rightarrow {\cal M}_i\otimes\Omega_X^1$
est un homomorphisme ${\cal O}_X$-lin\'eaire. Notons
$$
\Lambda_{i,n}\subset T^*{\cal C}_i
$$
le sous-champ ferm\'e dont les points sont
les paires $({\cal M}_i,\theta_i)$ telles que l'homomorphisme
$\theta_i$ soit nilpotent de niveau $n$, i.e. telles que
l'homomorphisme compos\'e
$$\displaylines{
\qquad {\cal M}_i\buildrel\theta_i\over\longrightarrow
{\cal M}_i\otimes_{{\cal O}_X}\Omega_X^1
\buildrel\theta_i\otimes{\rm id}\over\longrightarrow
{\cal M}_i\otimes_{{\cal O}_X}(\Omega_X^1)^{\otimes 2}
\rightarrow\cdots
\hfill\cr\hfill
\cdots\rightarrow
{\cal M}_i\otimes_{{\cal O}_X}(\Omega_X^1)^{\otimes (n-1)}
\buildrel\theta_i\otimes{\rm id}\over\longrightarrow
{\cal M}_i\otimes_{{\cal O}_X}(\Omega_X^1)^{\otimes n}\qquad}
$$
soit identiquement nul.

\th CONJECTURE 2.1
\enonce
Pour tout syst\`eme local irr\'eductible $L$ de rang $n$ sur $X$,
la vari\'et\'e caract\'eristique
$$
{\rm Car}(K_L)\subset T^*{\cal C}_n
$$
du faisceau pervers conjectural $K_L$ est \'egale \`a
$\Lambda_{n,n}$.
\endth

Pour motiver cette conjecture,
commen\c cons par \'etudier la vari\'et\'e caract\'eristique
$$
{\rm Car}(W_L)\subset T^*{\cal C}_0
$$
du faisceau pervers $W_L$.

\th LEMME 2.2 ([La] Th\'eor\`eme (3.1.13)(i))
\enonce
Le ferm\'e $\Lambda_{0,n}$ de $T^*{\cal C}_0$ est lagrangien.

Plus pr\'ecis\'ement, pour chaque suite d'entiers $\geq 0$
de longueur $n$, ${\bf d}=(d_1,\ldots ,d_n)$, telle que
$$
\sum_{i=1}^nid_i=m_0,
$$
notons
$$
X_{\bf d}^{(m_0)}\subset X^{(m_0)}
$$
la sous-vari\'et\'e localement ferm\'ee dont les points sont les
$D\in X^{(m_0)}$ de la forme
$$
D=\sum_{i=1}^niD_i,
$$
o\`u chaque $D_i$ est un diviseur effectif sur $X$, de degr\'e $d_i$
et sans multiplicit\'e \($D_i=x_{i,1}+\cdots +x_{i,d_i}$ avec
$x_{i,j}\not= x_{i,k}$, $\forall j\not= k$\), et o\`u les
supports des $D_i$ sont deux \`a deux disjoints \($x_{i',j'}\not=
x_{i'',j''}$ d\`es que $i'\not= i''$\),
et notons
$$
{\cal C}_{0,{\bf d}}\subset {\cal C}_0
$$
le sous-champ localement ferm\'e image de $X_{\bf d}^{(m_0)}$
par le morphisme lisse $\iota_{(m_0)}:X^{(m_0)}\rightarrow
{\cal C}_0$.  Alors, $\Lambda_{0,n}$ est la r\'eunion des
adh\'erences dans $T^*{\cal C}_0$ des sous-fibr\'es conormaux
$$
T_{{\cal C}_{0,{\bf d}}}^*{\cal C}_0\subset T^*{\cal C}_0.
$$
\endth

\rem Preuve
\endrem
L'\'enonc\'e du lemme est local pour la topologie \'etale sur
${\cal C}_0$ et donc pour la topologie \'etale sur $X$.
On peut donc supposer que $X={\bb A}^1$. Dans ce cas
$T^*{\cal C}_0$ et $\Lambda_{0,n}$ admettent la description
suivante. On a la pr\'esentation
$$
{\cal C}_0\cong [{\goth g}/G],
$$
o\`u $G=GL_{m_0}$ agit par conjugaison sur ${\goth g}=gl_{m_0}$,
et on peut identifier $\Lambda_{0,n}\subset T^*{\cal C}_0$
au quotient de
$$
\{(\xi,\xi^*)\in {\goth g}\times {\goth g}
\mid [\xi ,\xi^*]=0\hbox{ et }(\xi^*)^n=0\}
\subset \{(\xi,\xi^*)\in {\goth g}\times {\goth g}
\mid [\xi ,\xi^*]=0\}
$$
par l'action
$$
((\xi,\xi^*),g)\mapsto (g^{-1}\xi g,g^{-1}\xi^*g)
$$
de $G$ (comme d'habitude,
on a identifi\'e ${\goth g}^*$ \`a ${\goth g}$ par
l'isomorphisme
$$
{\goth g}\rightarrow {\goth g}^*,~\xi^*\mapsto
(\xi\mapsto {\rm tr}(\xi\xi^*))\,).
$$

Pour chaque suite d'entiers ${\bf d}$ comme dans l'\'enonc\'e
du lemme, notons ${\cal O}_{\bf d}\subset {\goth g}$ l'orbite de
la matrice de Jordan nilpotente ayant, pour chaque $i=1,\ldots ,n$,
exactement $d_i$ blocs de Jordan de taille $i$ et notons
${\goth g}_{\bf d}$ le ferm\'e de l'ouvert ${\goth g}^{\rm reg}$
des \'el\'ements r\'eguliers de
${\goth g}$ form\'e des $\xi$ dont le
polyn\^ome minimal est de la forme
$$
\prod_{i=1}^n\prod_{j=1}^{d_i}(T-\alpha_{i,j})^i
$$
avec des $\alpha_{i,j}$ deux \`a deux distincts.
On conclut en remarquant que
$$
\{(\xi,\xi^*)\in {\goth g}\times {\goth g}
\mid [\xi ,\xi^*]=0\hbox{ et }(\xi^*)^n=0\}=
\bigcup_{\bf d}\{(\xi,\xi^*)\in {\goth g}\times {\cal O}_{\bf d}
\mid [\xi ,\xi^*]=0\}
$$
(r\'eunion disjointe) et que, pour chaque ${\bf d}$,
$\{(\xi,\xi^*)\in {\goth g}\times {\cal O}_{\bf d}
\mid [\xi ,\xi^*]=0\}$ est l'adh\'erence
dans $ T^*{\goth g}$ du fibr\'e conormal
$$
T_{{\goth g}_{\bf d}}^*{\goth g}^{\rm reg}\subset
T^*{\goth g}^{\rm reg}\subset T^*{\goth g}.
$$

\hfill\hfill\cqfd

\th TH\'EOR\`EME 2.3 ([La] Th\'eor\`eme (3.1.13)(ii))
\enonce
Pour tout syst\`eme local irr\'eductible $L$ de rang $n$ sur $X$,
la vari\'et\'e caract\'eristique du faisceau pervers irr\'eductible de
Whittaker $W_L$ est exactement le ferm\'e lagrangien
$\Lambda_{0,n}$ de $T^*{\cal C}_0$.

\hfill\hfill\cqfd
\endth

Consid\'erons maintenant, pour chaque $i$, les champs cotangents
$T^*{\cal E}_i$ et $T^*{\cal E}_i^\vee$ aux champs alg\'ebriques
et lisses sur ${\bb C}$, ${\cal E}_i$ et ${\cal E}_i^\vee$
respectivement. Le premier est le champ des paires $(e_i,\eta_i)$
o\`u
$$
e_i: (\Omega_X^1)^{\otimes (i-1)}\rightarrow {\cal M}_i
$$
est un objet du champ ${\cal E}_i$ et o\`u
$$
\eta_i:{\cal M}_i\rightarrow [(\Omega_X^1)^{\otimes (i-1)}
\buildrel{e_i}\over\longrightarrow {\cal M}_i]
\otimes\Omega_X^1
$$
est un morphisme dans la cat\'egorie d\'eriv\'ee des
${\cal O}_X$-Modules ($[{\cal A}\rightarrow {\cal B}]$ d\'esigne
l'objet
$$
\cdots\rightarrow 0\rightarrow {\cal A}\rightarrow {\cal B}
\rightarrow 0\rightarrow\cdots
$$
de $D^{\rm b}({\cal O}_X)$, o\`u le ${\cal O}_X$-Module ${\cal A}$
est plac\'e en degr\'e $-1$ et le ${\cal O}_X$-Module ${\cal B}$
est plac\'e en degr\'e $0$). Le second est le champ des paires
$(e_i^\vee,\eta_i^\vee)$ o\`u
$$
e_i^\vee=\bigl((\Omega_X^1)^{\otimes i}\hookrightarrow
{\cal M}_{i+1}\twoheadrightarrow{\cal M}_i\bigr)
$$
est un objet du champ ${\cal E}_i^\vee$ et o\`u
$$
\eta_i^\vee :{\cal M}_{i+1}\rightarrow {\cal M}_i\otimes\Omega_X^1
$$
est un morphisme de ${\cal O}_X$-Modules.

L'isomorphisme de champs
${\cal E}_i^\circ\cong {\cal E}_{i-1}^{\vee\circ}$ induit un
isomorphisme de champs
$$
T^*{\cal E}_i^\circ\cong T^*{\cal E}_{i-1}^{\vee\circ}
$$
qui envoie
$(e_i,\eta_i)$ sur $(e_{i-1}^\vee ,\eta_{i-1}^\vee )$, o\`u
$$
e_{i-1}^\vee =\bigl((\Omega_X^1)^{\otimes (i-1)}
\buildrel{e_i}\over\hookrightarrow
{\cal M}_i\twoheadrightarrow{\cal M}_{i-1}\bigr)
$$
et
$$
\eta_{i-1}^\vee :
{\cal M}_i\buildrel{\eta_i}\over\longrightarrow
[(\Omega_X^1)^{\otimes (i-1)}\buildrel{e_i}\over\hookrightarrow
{\cal M}_i]\otimes\Omega_X^1\buildrel\sim\over\longrightarrow
{\cal M}_{i-1}\otimes\Omega_X^1.
$$

Les fibr\'es cotangents \`a un fibr\'e vectoriel et \`a son
fibr\'e dual sont toujours canoniquement isomorphes. Ici on a:

\th LEMME 2.4
\enonce
L'isomorphisme canonique
$$
T^*{\cal E}_i\cong T^*{\cal E}_i^\vee
$$
est donn\'e  de la fa\c con  suivante:

\decale{\rm (i)} \`a $(e_i,\eta_i)$ dans $T^*{\cal E}_i$
on associe $(e_i^\vee,\eta_i^\vee)$ dans $T^*{\cal E}_i^\vee$,
o\`u $e_i^\vee$ est l'extension
$$
(\Omega_X^1)^{\otimes i}\hookrightarrow
{\cal M}_{i+1}\twoheadrightarrow{\cal M}_i
$$
correspondant au morphisme compos\'e
$$
{\cal M}_i\buildrel{\eta_i}\over\longrightarrow
[(\Omega_X^1)^{\otimes (i-1)}
\buildrel{e_i}\over\longrightarrow
{\cal M}_i]\otimes\Omega_X^1\rightarrow
[(\Omega_X^1)^{\otimes (i-1)}\rightarrow 0]\otimes\Omega_X^1=
(\Omega_X^1)^{\otimes i}[1]
$$
dans $D^{\rm b}({\cal O}_X)$ et o\`u
$$
\eta_i^\vee :{\cal M}_{i+1}\rightarrow
[0\rightarrow {\cal M}_i]\otimes\Omega_X^1
={\cal M}_i\otimes\Omega_X^1
$$
est l'unique homomorphisme rendant commutatif le diagramme
$$
\matrix{
(\Omega_X^1)^{\otimes i}&{=\!=}&
[0\rightarrow (\Omega_X^1)^{\otimes (i-1)}]\otimes\Omega_X^1\cr
\noalign{\smallskip}
\downarrow&&\downarrow\cr
\noalign{\smallskip}
{\cal M}_{i+1}&\buildrel{\eta_i^\vee}\over\longrightarrow&
[0\rightarrow {\cal M}_i]\otimes\Omega_X^1\cr
\noalign{\smallskip}
\downarrow&&\downarrow\cr
\noalign{\smallskip}
{\cal M}_i&\buildrel{\eta_i}\over\longrightarrow
& [(\Omega_X^1)^{\otimes (i-1)}
\buildrel{e_i}\over\rightarrow {\cal M}_i]\otimes\Omega_X^1\cr
\noalign{\smallskip}
\downarrow&&\downarrow\cr
\noalign{\smallskip}
(\Omega_X^1)^{\otimes i}[1]&{=\!=}&
[(\Omega_X^1)^{\otimes (i-1)}\rightarrow 0]\otimes
\Omega_X^1,\cr}
$$

\decale{\rm (ii)} inversement, \`a
$(e_i^\vee ,\eta_i^\vee)$ dans $T^*{\cal E}_i^\vee$ on associe
$(e_i ,\eta_i)$ dans $T^*{\cal E}_i$, o\`u
$e_i\otimes\Omega_X^1$ est l'homomorphisme compos\'e
$$
(\Omega_X^1)^{\otimes i}\hookrightarrow {\cal M}_{i+1}
\buildrel{\eta_i^\vee}\over\longrightarrow
{\cal M}_i\otimes\Omega_X^1
$$
et o\`u $\eta_i$ est le morphisme compos\'e
$$
{\cal M}_i=[(\Omega_X^1)^{\otimes i}
\hookrightarrow {\cal M}_{i+1}]\buildrel\eta_i^\vee\over
\longrightarrow [(\Omega_X^1)^{\otimes i}\rightarrow
{\cal M}_i\otimes\Omega_X^1]=
[(\Omega_X^1)^{\otimes (i-1)}\buildrel{e_i}\over
\longrightarrow
{\cal M}_i]\otimes\Omega_X^1.
$$
\endth

\rem Preuve
\endrem
Seule la d\'efinition de $\eta_i^\vee$ dans la partie (i)
n\'ecessite quelques \'eclaircissements.
L'existence de $\eta_i^\vee$ r\'esulte de l'axiome (TR3) des
cat\'egories triangul\'ees, puisque les deux verticales du
diagramme sont des triangles distingu\'es. Quant \`a l'unicit\'e,
elle vient ce que ${\rm Hom}_{{\cal O}_X}({\cal M}_i,
(\Omega_X^1)^{\otimes i})=(0)$ par d\'efinition de ${\cal C}_i$.
En effet, comme ${\rm Hom}_{{\cal O}_X}((\Omega_X^1)^{\otimes i},
(\Omega_X^1)^{\otimes i})$ est de dimension $1$, soit on a aussi
${\rm Hom}_{{\cal O}_X}({\cal M}_{i+1},
(\Omega_X^1)^{\otimes i})=(0)$ et l'unicit\'e  est imm\'ediate
(la diff\'erence entre deux $\eta_i^\vee$ possibles
se factorise n\'ecessairement en
$$
{\cal M}_{i+1}\rightarrow
[0\rightarrow (\Omega_X^1)^{\otimes (i-1)}]\otimes\Omega_X^1
\rightarrow
[0\rightarrow {\cal M}_i]\otimes\Omega_X^1
$$
et est donc automatiquement nulle),
soit l'extension $e_i^\vee$ est scind\'ee et la diff\'erence entre
deux $\eta_i^\vee$ possibles se factorise en
$$
{\cal M}_{i+1}\rightarrow {\cal M}_i\rightarrow
[0\rightarrow (\Omega_X^1)^{\otimes (i-1)}]\otimes\Omega_X^1
\rightarrow
[0\rightarrow {\cal M}_i]\otimes\Omega_X^1
$$
et est aussi automatiquement nulle.

\hfill\hfill\cqfd

Dans la suite, on notera simplement $\Omega^i$ et
$(\cdot )\otimes\Omega^i$ les expressions
$(\Omega_X^1)^{\otimes i}$
et $(\cdot )\otimes_{{\cal O}_X}(\Omega_X^1)^{\otimes i}$.

Consid\'erons le champ ${\cal U}$ des diagrammes commutatifs de
${\cal O}_X$-Modules coh\'erents,
$$
({\cal M}_\bullet ,e_\bullet ,\eta_\bullet )=\left\{\matrix{
\Omega^{n-1}&\buildrel{e_n}\over
{\lhook\kern -1mm\joinrel{\hbox to 12mm{\rightarrowfill}}}&
{\cal M}_n
&\twoheadrightarrow&{\cal M}_{n-1}\cr
\noalign{\smallskip}
\big|\!\big|&&\mapdownright{\eta_n}&&
\mapdownright{\eta_{n-1}}\cr
\noalign{\smallskip}
\Omega^{n-1}&\buildrel{e_{n-1}\otimes\Omega}\over
{\lhook\kern -1mm\joinrel{\hbox to 12mm{\rightarrowfill}}}
&{\cal M}_{n-1}\otimes\Omega
&\twoheadrightarrow&{\cal M}_{n-2}\otimes\Omega\cr
\noalign{\smallskip}
\big|\!\big|&&\mapdownright{\eta_{n-1}\otimes\Omega}&&
\mapdownright{\eta_{n-2}\otimes\Omega}\cr
\noalign{\smallskip}
\Omega^{n-1}&\buildrel{e_{n-2}\otimes\Omega^2}\over
{\lhook\kern -1mm\joinrel{\hbox to 12mm{\rightarrowfill}}}
&{\cal M}_{n-2}\otimes\Omega^2
&\twoheadrightarrow&{\cal M}_{n-3}\otimes\Omega^2\cr
\cdot&&\cdot&&\cdot\cr
\noalign{\vskip -2mm}
\cdot&&\cdot&&\cdot\cr
\noalign{\vskip -2mm}
\cdot&&\cdot&&\cdot\cr
\noalign{\smallskip}
\big|\!\big|&&\mapdownright{\eta_2\otimes\Omega^{n-2}}&&
\mapdownright{\eta_1\otimes\Omega^{n-2}}\cr
\noalign{\smallskip}
\Omega^{n-1}&\buildrel{e_1\otimes\Omega^{n-1}}\over
{\lhook\kern -1mm\joinrel{\hbox to 12mm{\rightarrowfill}}}
&{\cal M}_1\otimes\Omega^{n-1}
&\twoheadrightarrow&{\cal M}_0\otimes\Omega^{n-1}\cr
\noalign{\smallskip}
&&\mapdownright{\eta_1\otimes\Omega^{n-1}}&&\cr
\noalign{\smallskip}
&&{\cal M}_0\otimes\Omega^n&&\cr}\kern 3mm\right\},
$$
o\`u chaque ${\cal M}_i$ est de
rang g\'en\'erique $i$ et de degr\'e $m_i$ et satisfait \`a la
condition
$$
{\rm Hom}_{{\cal O}_X}({\cal M}_i,\Omega^{2n-1})=(0)
$$
et o\`u chaque ligne est une suite exacte courte.

Il est facile de voir que ${\cal U}$ est alg\'ebrique et est en fait
un ouvert commun \`a tous les champs $T^*{\cal E}_i^\vee\cong
T^*{\cal E}_i$. Plus pr\'ecis\'ement, pour chaque $i=0,\ldots ,n-1$,
la donn\'e de la partie
$$\matrix{
\Omega^{n-1}&\buildrel{e_{i+1}\otimes\Omega^{n-i-1}}\over
{\lhook\kern -1mm\joinrel{\hbox to 18mm{\rightarrowfill}}}
&{\cal M}_{i+1}\otimes\Omega^{n-i-1}
&\twoheadrightarrow&{\cal M}_i\otimes\Omega^{n-i-1}\cr
\noalign{\smallskip}
&&\mapdownright{\eta_{i+1}\otimes\Omega^{n-i-1}}&&\cr
\noalign{\smallskip}
&&{\cal M}_i\otimes\Omega^{n-i}
&&\cr}
$$
du diagramme ci-dessus permet de reconstruire ce diagramme
tout entier et tout diagramme de ${\cal O}_X$-Modules coh\'erents,
$$\matrix{
\Omega^i&\buildrel{e_{i+1}}\over
{\lhook\kern -1mm\joinrel{\hbox to 9mm{\rightarrowfill}}}
&{\cal M}_{i+1}
&\twoheadrightarrow&{\cal M}_i\cr
\noalign{\smallskip}
&&\mapdownright{\eta_{i+1}}&&\cr
\noalign{\smallskip}
&&{\cal M}_{i}\otimes\Omega , &&\cr}
$$
o\`u ${\cal M}_i$ est de
rang g\'en\'erique $i$ et de degr\'e $m_i$, o\`u
${\rm Hom}_{{\cal O}_X}({\cal M}_{i+1},\Omega^{2n-1})=(0)$
et o\`u la ligne sup\'erieure est une suite exacte courte,
d\'efinit un point de
$$
T^*{\cal E}_{i+1}^\circ\cong T^*{\cal E}_i^{\vee\circ}\subset
T^*{\cal E}_i^\vee .
$$

On note $\Sigma_{0,n}$ le sous-champ ferm\'e de
$T^*{\cal E}_0^\vee$ image de
${\cal E}_0^\vee\times_{{\cal C}_0}\Lambda_{0,n}$
par l'immersion ferm\'ee
$$
T^*\pi_0^\vee :{\cal E}_0^\vee\times_{{\cal C}_0}T^*{\cal C}_0
\hookrightarrow T^*{\cal E}_0^\vee
$$
et on note $\Sigma_{n,n}$ le sous-champ ferm\'e de
$T^*{\cal E}_n$ image de
${\cal E}_n\times_{{\cal C}_n}\Lambda_{n,n}$
par l'immersion ferm\'ee
$$
T^*\pi_n :{\cal E}_n\times_{{\cal C}_n}T^*{\cal C}_n
\hookrightarrow T^*{\cal E}_n.
$$

\th LEMME 2.5
\enonce
Les intersections des ferm\'es $\Sigma_{0,n}$ et $\Sigma_{n,n}$
avec l'ouvert ${\cal U}$,  commun \`a $T^*{\cal E}_0^\vee$ et
$T^*{\cal E}_n$, coincident.
\endth

\rem Preuve
\endrem
Par d\'efinition, $\Sigma_{0,n}\cap {\cal U}$ (resp.
$\Sigma_{n,n}\cap {\cal U}$) est form\'e des diagrammes
$({\cal M}_\bullet ,e_\bullet ,\eta_\bullet )$ tels que
$\eta_1$ (resp. $\eta_n$) se factorise en
$$
{\cal M}_1\twoheadrightarrow {\cal M}_0\buildrel\theta_0\over
\longrightarrow {\cal M}_0\otimes\Omega
$$
(resp.
$$
{\cal M}_n\buildrel\theta_n\over\longrightarrow
{\cal M}_n\otimes\Omega\twoheadrightarrow
{\cal M}_{n-1}\otimes\Omega\,),
$$
avec $\theta_0$ (resp. $\theta_n$) nilpotent de niveau $n$.
Remarquons que, si une telle factorisation existe, elle est
n\'ecessairement unique (c'est \'evident pour $\eta_0$ et cela
resulte pour $\eta_n$ de
${\rm Hom}_{{\cal O}_X}({\cal M}_n,\Omega^{2n-1})=(0)$).

En consid\'erant le diagramme
$$\matrix{
{\cal M}_n\kern -4mm &\twoheadrightarrow &
{\cal M}_{n-1}\kern -5mm  &
\twoheadrightarrow &\cdots &
\twoheadrightarrow &{\cal M}_1\kern -5mm&
\twoheadrightarrow &{\cal M}_0\cr
\noalign{\smallskip}
\mapdownleft{\eta_n}\kern -4mm &&\mapdownleft{\eta_{n-1}}
\kern -5mm &&&&\mapdownleft{\eta_1}\kern -5mm&
\swarrow_{\theta_0}\kern -2mm&\cr
\noalign{\smallskip}
{\cal M}_{n-1}\otimes\Omega\kern -4mm
&\twoheadrightarrow &
{\cal M}_{n-2}\otimes\Omega\kern -5mm &
\twoheadrightarrow&\cdots &
\twoheadrightarrow &
{\cal M}_0\otimes\Omega\kern -5mm &&\cr
\noalign{\smallskip}
\downarrow\kern -4mm && \downarrow\kern -5mm &&&
\swarrow &&&\cr
\noalign{\medskip}
\cdot\kern -4mm &&\cdot\kern -5mm &&\cdot&&&&\cr
\noalign{\medskip}
\mapdownleft{\eta_2\otimes\Omega^{n-2}}\kern -4mm &&
\mapdownleft{\eta_1\otimes\Omega^{n-2}}\kern -5mm &
\swarrow_{\theta_0\otimes\Omega^{n-2}}\kern -9mm &&&&&
\cr
\noalign{\smallskip}
{\cal M}_1\otimes\Omega^{n-1}\kern -4mm &
\twoheadrightarrow &
{\cal M}_0\otimes\Omega^{n-1}\kern -5mm&&&&&&\cr
\noalign{\smallskip}
\mapdownleft{\eta_1\otimes\Omega^{n-1}}\kern -4mm&
\swarrow_{\theta_0\otimes\Omega^{n-1}}\kern -10mm&&&
&&&&\cr
\noalign{\smallskip}
{\cal M}_0\otimes\Omega^n\kern -4mm &&&&&&&&\cr}
$$
on voit que $({\cal M}_\bullet ,e_\bullet ,\eta_\bullet )\in
{\cal U}$ est dans $\Sigma_{0,n}$ si et seulement si
$$
(\eta_1\otimes\Omega^{n-1})\circ\cdots\circ
(\eta_{n-1}\otimes\Omega )\circ\eta_n =0
$$
(si cette condition est satisfaite, on a
$$
\eta_1\circ e_1=\eta_1\circ\cdots\circ
(\eta_{n-1}\otimes\Omega^{2-n})\circ
(\eta_n\otimes\Omega^{1-n})\circ (e_n\otimes\Omega^{1-n})=0
$$
et donc $\eta_1$ se factorise bien en
${\cal M}_1\twoheadrightarrow {\cal M}_0\buildrel\theta_0\over
\longrightarrow {\cal M}_0\otimes\Omega$).

De m\^eme,
en consid\'erant le diagramme
$$\matrix{
 &&&&&&&&\kern -4mm{\cal M}_n\cr
\noalign{\smallskip}
&&&&&&&\kern -11mm{}{\theta_n}\swarrow
&\kern -4mm\mapdownright{\eta_n}\cr
\noalign{\smallskip}
&&&&&&\kern -7mm{\cal M}_n\otimes\Omega &
\kern -4mm\twoheadrightarrow  &\kern -4mm
{\cal M}_{n-1}\otimes\Omega\cr
\noalign{\smallskip}
 &&&&&\kern -9mm {}^{\theta_n\otimes\Omega}\swarrow
&\kern -7mm\mapdownright{\eta_n\otimes\Omega}
&&\kern -4mm\mapdownright{\eta_{n-1}\otimes\Omega}\cr
\noalign{\medskip}
&&&&\cdot&&\kern -7mm\cdot&&\kern -4mm\cdot\cr
\noalign{\medskip}
&&&\swarrow&&&\kern -7mm \downarrow&&\kern -4mm
\downarrow \cr
\noalign{\smallskip}
 &&\kern -5mm {\cal M}_n\otimes\Omega^{n-1}&
\twoheadrightarrow&\cdots &\twoheadrightarrow&
{\cal M}_2\otimes\Omega^{n-1}&\kern -4mm\twoheadrightarrow &
\kern -4mm {\cal M}_1\otimes\Omega^{n-1}\cr
\noalign{\smallskip}
&\kern -13mm{}^{\theta_n\otimes\Omega^{n-1}}\swarrow&
\mapdownright{\eta_n\otimes\Omega^{n-1}} &&&&\kern -5mm
\mapdownright{\eta_2\otimes\Omega^{n-1}}&&\kern -4mm
\mapdownright{\eta_1\otimes\Omega^{n-1}}\cr
\noalign{\smallskip}
{\cal M}_n\otimes\Omega^n&\kern -5mm\twoheadrightarrow&
\kern -5mm{\cal M}_{n-1}\otimes\Omega^n&
\twoheadrightarrow &\cdots &\twoheadrightarrow  &\kern -5mm
{\cal M}_1\otimes\Omega^n&\kern -4mm\twoheadrightarrow &
\kern -4mm{\cal M}_0\otimes\Omega^n\cr}
$$
on voit que $({\cal M}_\bullet ,e_\bullet ,\eta_\bullet )\in
{\cal U}$ est dans $\Sigma_{n,n}$ si et seulement si
$$
(\eta_1\otimes\Omega^{n-1})\circ\cdots\circ
(\eta_{n-1}\otimes\Omega )\circ\eta_n =0
$$
(d'une part, si cette condition est satisfaite, on a
$$
(e_n^\vee\otimes\Omega )\circ\eta_n=(e_1^\vee\otimes\Omega^n)
\circ (\eta_1\otimes\Omega^{n-1})\circ\cdots\circ
(\eta_{n-1}\otimes\Omega )\circ\eta_n=0,
$$
o\`u $e_i^\vee :{\cal M}_{i-1}\rightarrow\Omega^{i-1}[1]$
est la classe de l'extension ${\cal M}_i$ de ${\cal M}_{i-1}$ par
$\Omega^{i-1}$ pour $i=1,\ldots ,n$,
et donc $\eta_n$ se factorise bien en
${\cal M}_n\buildrel\theta_n\over\longrightarrow
{\cal M}_n\otimes\Omega\twoheadrightarrow
{\cal M}_{n-1}\otimes\Omega$; d'autre part, on a
${\rm Hom}_{{\cal O}_X}({\cal M}_n,\Omega^i)=(0)$ pour tout
$i\leq 2n-1$ et donc
$$
{\rm Hom}_{{\cal O}_X}({\cal M}_n,{\rm Ker}({\cal M}_n
\twoheadrightarrow{\cal M}_0)\otimes\Omega^n)=(0)
$$
puisque ${\rm Ker}({\cal M}_n\twoheadrightarrow{\cal M}_0)$
est extension successive des $\Omega^i$ pour $i=0,\ldots ,n-1$).

\hfill\hfill\cqfd

On va maintenant d\'efinir, pour chaque $i=1,\ldots ,n-1$,
un sous-champ
$$
\Sigma_{i,n}\subset T^*{\cal E}_i\cong T^*{\cal E}_i^\vee
$$
dont on peut penser qu'il est exactement la vari\'et\'e
carat\'eristique du faisceau pervers $W_{L,i}$ pour n'importe
quel syst\`eme local irr\'eductible $L$ de rang $n$ sur $X$.

On utilisera la notion suivante. Soient $B$ et $C$ deux objets
d'une cat\'egorie triangul\'ee ${\cal K}$ et soient $f:B\rightarrow C$
et $\eta :B\rightarrow C$ deux morphismes entre ces objets.
On peut inclure le diagramme
$$
\matrix{B&\maprightover{f}&C\cr
\noalign{\medskip}
\mapdownright{\eta}&&\cr
\noalign{\medskip}
C&&\cr}
$$
dans un diagramme commutatif
$$
\matrix{
\cdot&&\cdot&&\cdot&&\cdot\cr
\noalign{\vskip -2mm}
\cdot&&\cdot&&\cdot&&\cdot\cr
\noalign{\vskip -2mm}
\cdot&&\cdot&&\cdot&&\cdot\cr
\big|\!\big|&&\mapdownright{\eta_3}&&
\mapdownright{\eta_2}&&\big|\!\big|&\cr
\noalign{\medskip}
A&\maprightover{e_2}&B_2&\maprightover{f_2}&B_1&
\maprightover{e_1^\vee}&A[1]\cr
\noalign{\medskip}
\big|\!\big|&&\mapdownright{\eta_2}&&
\mapdownright{\eta_1}&&\big|\!\big|&\cr
\noalign{\medskip}
A&\maprightover{e_1}&B_1&\maprightover{f_1}&B&
\maprightover{e^\vee}&A[1]\cr
\noalign{\medskip}
\big|\!\big|&&\mapdownright{\eta_1}&&
\mapdownright{\eta}&&\big|\!\big|&\cr
\noalign{\medskip}
A&\maprightover{e}&B&\maprightover{f}&C&
\maprightover{e_{-1}^\vee}&A[1]\cr
\noalign{\medskip}
\big|\!\big|&&\mapdownright{\eta}&&
\mapdownright{\eta_{-1}}&&\big|\!\big|&\cr
\noalign{\medskip}
A&\maprightover{e_{-1}}&C&\maprightover{f_{-1}}&C_1&
\maprightover{e_{-2}^\vee}&A[1]\cr
\noalign{\medskip}
\big|\!\big|&&\mapdownright{\eta_{-1}}&&
\mapdownright{\eta_{-2}}&&\big|\!\big|&\cr
\noalign{\medskip}
A&\maprightover{e_{-2}}&C_1&\maprightover{f_{-2}}&C_2&
\maprightover{e_{-3}^\vee}&A[1]\cr
\noalign{\medskip}
\big|\!\big|&&\mapdownright{\eta_{-2}}&&
\mapdownright{\eta_{-3}}&&\big|\!\big|&\cr
\noalign{\medskip}
\cdot&&\cdot&&\cdot&&\cdot\cr
\noalign{\vskip -2mm}
\cdot&&\cdot&&\cdot&&\cdot\cr
\noalign{\vskip -2mm}
\cdot&&\cdot&&\cdot&&\cdot\cr}
$$
o\`u toutes les lignes sont des triangles distingu\'es de ${\cal K}$,
ce diagramme \'etant unique \`a isomorphisme pr\`es, et
on dit que le morphisme $\eta$ est nilpotent de niveau $(n-i,i)$
(relativement \`a $f$) si
$$
\eta_{-i}\circ\cdots\circ\eta_{-1}\circ\eta\circ
\eta_1\circ\eta_2\circ\cdots\circ\eta_{n-i}=0.
$$
Plus g\'en\'eralement, si l'on dispose  d'une notion
de torsion dans ${\cal K}$, i.e. d'un
foncteur exact inversible $C\mapsto C(1)$ de ${\cal K}$ dans elle
m\^eme, d'inverse not\'e
$C\mapsto C(-1)$, et si $\eta$ est maintenant un morphisme de
$B$ dans $C(1)$, $f$ \'etant toujours un morphisme de $B$ dans
$C$, on dit encore que $\eta$ est nilpotent de niveau $(n-i,i)$
(relativement \`a $f$) si
$$
\eta_{-i}(i)\circ\cdots\circ\eta_{-1}(1)\circ\eta\circ
\eta_1(-1)\circ\eta_2(-2)\circ\cdots\circ\eta_{n-i}(i-n)=0,
$$
o\`u $C\mapsto C(j)$ est la puissance $j$-i\`eme de
$C\mapsto C(1)$ quel que soit $j\in {\bb Z}$.

On d\'efinit alors $\Sigma_{i,n}$ comme le ferm\'e de
$T^*{\cal E}_i\cong T^*{\cal E}_i^\vee$ form\'e des
$$
\matrix{\Omega^{i-1}&\maprightover{e_i}&\kern -10mm
{\cal M}_i\cr
\noalign{\medskip}
&&\kern -10mm\mapdownright{\eta_i}\cr
\noalign{\medskip}
&&\kern -10mm [\Omega^{i-1}\maprightover{e_i}{\cal M}_i]
\otimes\Omega\cr}
$$
tels que $\eta_i$ soit nilpotent de niveau
$(n-i,i)$ (relativement au morphisme canonique
$$
{\rm can}:{\cal M}_i\rightarrow
[\Omega^{i-1}\maprightover{e_i}{\cal M}_i]
$$
dans $D^{\rm b}({\cal O}_X)$). Par construction, il est clair que
$$
\Sigma_{i,n}\cap {\cal U}=\Sigma_{0,n}\cap {\cal U}
$$
et que
$$
\Sigma_{i,n}\cap T^*{\cal E}_i^\circ=
\Sigma_{i-1,n}\cap T^*{\cal E}_{i-1}^{\vee\circ}
$$
pour tout $i=1,\ldots,n$.

\vskip 5mm
{\bf 3. Une variante de la construction du complexe
$W_{L,n}^\circ$.}
\vskip 5mm

On se placera dans le contexte de l'expos\'e II. On a donc
fix\'e un corps alg\'ebriquement clos $k$ de caract\'eristique
$p>0$ et $X$ est une courbe connexe, projective, lisse et de genre
$g\geq 2$ sur $k$.

Pour chaque entier $i=1,\ldots ,n$, consid\'erons l'espace $D_i$
des classes d'isomorphie de drapeaux
$$
{\cal P}_{i,\bullet}=\bigl((0)={\cal P}_{i,0}\subset {\cal P}_{i,1}
\subset\cdots\subset {\cal P}_{i,i}\bigr)
$$
de fibr\'es vectoriels sur $X$, munis, pour chaque $j=1,\ldots ,i$,
d'un isomorphisme de la $j$-i\`eme composante du gradu\'e
$$
{\rm gr}_j{\cal P}_{i,\bullet}\buildrel{\rm dfn}\over{=\!=}
{\cal P}_{i,j}/{\cal P}_{i,j-1}
$$
sur le fibr\'e en droites $\Omega^{i-j}=
(\Omega_X^1)^{\otimes (i-j)}$.

\th LEMME 3.1
\enonce
Pour $i=1,\ldots ,n$, l'espace $D_i$ est repr\'esentable par un
sch\'ema affine sur $k$, non canoniquement isomorphe \`a
${\bb A}_k^{i-1}$.

Sur $X\times D_i$, on a un drapeau universel
$\widetilde{\cal P}_{i,\bullet}$.
\endth

\rem Preuve
\endrem
Soit $S$ un sch\'ema affine, connexe, r\'eduit et de type fini sur
$k$ et soit ${\cal L}$ un fibr\'e vectoriel sur $X\times S$, de dual
not\'e ${\cal L}^\vee$. On suppose que la fonction
$$
s\mapsto H^1(X\times s,{\cal L}^\vee |X\times s)
$$
est constante sur $S$, de sorte que $R^0{\rm pr}_{S*}{\cal L}^\vee$
et $R^1{\rm pr}_{S*}{\cal L}^\vee$ sont des fibr\'es vectoriels sur
$S$, dont la formation commute \`a tout changement de base
$T\rightarrow S$ et dont les fibr\'es duaux sont
$R^1{\rm pr}_{S*}(\Omega_X^1\otimes_{{\cal O}_X}{\cal L})$ et
$R^0{\rm pr}_{S*}(\Omega_X^1\otimes_{{\cal O}_X}{\cal L})$
respectivement. Soit $V$ le $S$-sch\'ema des sections
du fibr\'e vectoriel $R^1{\rm pr}_{S*}{\cal L}^\vee$, i.e.
$$
V={\bb V}\bigl(R^0{\rm pr}_{S*}(\Omega_X^1\otimes_{{\cal O}_X}
{\cal L})\bigr).
$$
Alors, sur $X\times V$, on a une extension universelle
de ${\cal O}_{X\times V}$ par ${\cal L}\otimes_{{\cal O}_S}
{\cal O}_V$. En effet, on a une section universelle
$$
{\cal O}_V\rightarrow R^1{\rm pr}_{V*}({\cal L}^\vee
\otimes_{{\cal O}_S}{\cal O}_V)
$$
qui se rel\`eve uniquement en une section
$$
{\cal O}_V\rightarrow R{\rm pr}_{V*}({\cal L}^\vee
\otimes_{{\cal O}_S}{\cal O}_V)[1]
$$
dans la cat\'egorie d\'eriv\'ee $D_{\rm coh}^{\rm b}({\cal O}_V)$
des complexes born\'es de ${\cal O}_V$-Modules \`a cohomologie
coh\'erente puisque l'on a le triangle distingu\'e
$$
R^0{\rm pr}_{V*}({\cal L}^\vee\otimes_{{\cal O}_S}{\cal O}_V)[1]
\rightarrow
R{\rm pr}_{V*}({\cal L}^\vee\otimes_{{\cal O}_S}{\cal O}_V)[1]
\rightarrow
R^1{\rm pr}_{V*}({\cal L}^\vee\otimes_{{\cal O}_S}{\cal O}_V)
\rightarrow
$$
et que, $V$ \'etant affine, on a
$$
H^n(V,R^0{\rm pr}_{V*}({\cal L}^\vee\otimes_{{\cal O}_S}
{\cal O}_V))=(0)\qquad (\forall n\geq 0).
$$
Par adjonction, ce rel\`evement \'equivaut \`a
une section
$$
{\cal O}_{X\times V}\rightarrow {\cal L}^\vee\otimes_{{\cal O}_S}
{\cal O}_V[1]
$$
dans $D_{\rm coh}^{\rm b}({\cal O}_{X\times V})$.
Cette derni\`ere section est l'extension duale de l'extension
cherch\'ee.

Pour d\'emontrer le lemme, on proc\`ede maintenant par
r\'ecurrence sur $i$.  L'espace $D_1$ est r\'eduit \`a un point.
Pour $i=2,\ldots ,n$, on suppose construit l'espace $D_{i-1}\cong
{\bb A}_k^{i-2}$ et le drapeau universel
$\widetilde{\cal P}_{i-1,\bullet}$.
On applique alors la construction pr\'ec\'edente \`a $S=D_{i-1}$
et \`a ${\cal L}=\Omega^{-(i-1)}\otimes_{{\cal O}_X}
\widetilde{\cal P}_{i-1,i-1}$. Les hypoth\`eses de cette
construction sont bien v\'erifi\'ees. En effet,
pour tout ${\cal P}_{i-1,\bullet}\in D_{i-1}$, la fl\`eche de
restriction
$$
{\rm Ext}_{{\cal O}_X}^1({\cal P}_{i-1}^{i-1},\Omega^{i-1})
\rightarrow
{\rm Ext}_{{\cal O}_X}^1({\cal P}_{i-1}^1,\Omega^{i-1})=
{\rm Ext}_{{\cal O}_X}^1(\Omega^{i-2},\Omega^{i-1})\cong
H^1(X,\Omega )\cong k
$$
est un isomorphisme puisque
${\rm Ext}_{{\cal O}_X}^1(\Omega^j,\Omega^{i-1})=(0)$
pour tout entier $j\leq i-3$ (on rappelle que l'on a suppos\'e
$g\geq 2$). Comme le fibr\'e en droites
$$
R^1{\rm pr}_{D_{i-1}*}{\cal H}{\it om}_{{\cal O}_X}
(\widetilde{\cal P}_{i-1,i-1},\Omega^{i-1})
$$
est n\'ecessairement trivial sur $D_{i-1}\cong {\bb A}_k^{i-2}$,
on a bien $D_i\cong {\bb A}_k^{i-1}$.

\hfill\hfill\cqfd
\vskip 2mm

On d\'efinit un morphisme de $k$-sch\'ema
$$
\theta_n :D_n\rightarrow {\bb A}_k^1
$$
en envoyant ${\cal P}_{n,\bullet}$ sur
$$
\sum_{i=1}^{n-1}u_i\in H^1(X,\Omega )\cong k,
$$
o\`u $u_i$ est la classe de l'extension
$$
\bigl(0\rightarrow {\rm gr}_i{\cal P}_{n,\bullet}\rightarrow
{\cal P}_{n,i+1}/{\cal P}_{n,i-1}\rightarrow
{\rm gr}_{i+1}{\cal P}_{n,\bullet}\rightarrow 0)\otimes
\Omega^{i+1-n}
$$
de ${\cal O}_X$ par $\Omega$.

Sur $D_n\times {\cal C}_n$, on a le fibr\'e vectoriel
$$
{\cal Z}_n\rightarrow D_n\times {\cal C}_n
$$
de fibre ${\rm Hom}_{{\cal O}_X}({\cal P}_{n,n},{\cal M}_n)$
en $({\cal P}_{n,\bullet} ,{\cal M}_n)$ (par d\'efinition de l'ouvert
${\cal C}_n$ du champ des ${\cal O}_X$-Modules coh\'erents, on a
$$
{\rm Hom}_{{\cal O}_X}({\cal M}_n,\Omega^j)=(0)
$$
pour tout $j=0,\ldots ,n$ et donc
$$
{\rm Ext}_{{\cal O}_X}^1({\cal P}_{n,n},{\cal M}_n)\cong
{\rm Hom}_{{\cal O}_X}({\cal M}_n,
{\cal P}_{n,n}\otimes\Omega)^\vee =(0)\,).
$$
La projection canonique de ce fibr\'e se factorise en
$$
{\cal Z}_n\buildrel{\alpha_n}\over\longrightarrow D_n\times
{\cal E}_n\buildrel{{\rm id}\times\pi_n}\over\longrightarrow
D_n\times {\cal C}_n
$$
puisque, pour chaque $({\cal P}_{n,n}\rightarrow {\cal M}_n)$
dans ${\cal Z}_n$, le compos\'e
$$
\Omega^{n-1}={\cal P}_{n,1}\hookrightarrow{\cal P}_{n,n}
\rightarrow {\cal M}_n
$$
est un \'el\'ement de ${\cal E}_n$. Le morphisme $\alpha_n$
ci-dessus est un fibr\'e affine sous le fibr\'e vectoriel sur
$D_n\times {\cal E}_n$ de fibre ${\rm Hom}_{{\cal O}_X}
({\cal P}_{n,n}/\Omega^{n-1},{\cal M}_n)$ en $({\cal P}_{n,n},
\Omega^{n-1}\hookrightarrow {\cal M}_n)$ puisque l'on a encore
${\rm Ext}_{{\cal O}_X}^1({\cal P}_{n,n}/\Omega^{n-1},{\cal M}_n)
=(0)$. En particulier, $\alpha_n$ est lisse purement de dimension
relative
$$
(n-1)m_0=(n-1)m_n-n\sum_{j=2}^n(n-j)(2g-2)+n(n-1)(1-g).
$$

Soit ${\cal Z}_n^\circ\subset {\cal Z}_n$ l'image inverse de
l'ouvert $D_n\times{\cal E}_n^\circ$ de $D_n\times{\cal E}_n$
et soit ${\cal Z}_n^{\circ\circ}\subset {\cal Z}_n^\circ$
l'ouvert de ${\cal Z}_n$ o\`u la fl\`eche ${\cal P}_{n,n}\rightarrow
{\cal M}_n$ est injective. On a un morphisme de champs
$$
\beta_n^{\circ\circ}:
{\cal Z}_n^{\circ\circ}\rightarrow D_n\times {\cal C}_0,
\quad ({\cal P}_{n,n}\hookrightarrow {\cal M}_n)\mapsto
({\cal P}_{n,n},{\cal M}_n/{\cal P}_{n,n}),
$$
qui fait de ${\cal Z}_n^{\circ\circ}$ un fibr\'e vectoriel sur $D_n
\times{\cal C}_0$ de fibre en $({\cal P}_{n,n},{\cal M}_0)$
l'espace vectoriel
$$
{\rm Ext}_{{\cal O}_X}^1({\cal M}_0,{\cal P}_{n,n})
$$
(on a trivialement ${\rm Hom}_{{\cal O}_X}({\cal M}_0,
{\cal P}_{n,n})=(0)$). En particulier, le morphisme
$\beta_n^{\circ\circ}$ est lisse, purement de dimension relative
$nm_0$. On remarquera que ce morphisme se factorise \`a travers
le morphisme
$$
D_n\times {\cal E}_0^\vee\buildrel{{\rm id}\times
\pi_0^\vee}\over\longrightarrow D_n\times {\cal C}_0
$$
car on peut pousser l'extension
$$
{\cal P}_{n,n}\hookrightarrow {\cal M}_n\twoheadrightarrow
{\cal M}_0
$$
par le morphisme quotient
$$
{\cal P}_{n,n}\twoheadrightarrow {\cal P}_{n,n}/{\cal P}_{n,n-1}
\cong {\cal O}_X
$$
et obtenir ainsi un point de ${\cal E}_0^\vee$.

On a construit un diagramme de champs alg\'ebriques sur $k$
$$
\matrix{&&&&&&{\cal Z}_n^{\circ\circ}&&\cr
\noalign{\smallskip}
&&&&&{}^{\alpha_n^{\circ\circ}}\kern -2mm\swarrow &&
\searrow^{\beta_n^{\circ\circ}}&\cr
\noalign{\smallskip}
{\bb A}_k^1&\mapleftover{\theta_n}&D_n&\mapleftover
{{\rm pr}_{D_n}}&D_n\times {\cal E}_n^\circ
\kern -4mm &&&&\kern -4mm
D_n\times {\cal C}_0\cr
\noalign{\smallskip}
&&&&\mapdownright{{\rm pr}_{{\cal E}_n^\circ}}&&&&\cr
\noalign{\smallskip}
&&&&{\cal E}_n^\circ &&&&\cr}
$$
o\`u $\alpha_n^{\circ\circ}$ est la restriction de $\alpha_n$
\`a l'ouvert ${\cal Z}_n^{\circ\circ}$.

\th PROPOSITION 3.2
\enonce
Avec les notations ci-dessus, pour tout
$\overline{\bb Q}_\ell$-faisceau $L$, lisse, irr\'eductible et de
rang $n$ sur $X$, le complexe $W_{L,n}^\circ$ est canoniquement
isomorphe dans $D_{\rm c}^{\rm b}({\cal E}_n^\circ ,\overline{\bb Q}_\ell )$
au complexe
$$
R{\rm pr}_{{\cal E}_n^\circ,!}\bigl(R{\alpha}_{n,!}^{\circ\circ}
\beta_n^{\circ\circ *}(\overline{\bb Q}_\ell[n-1]\boxtimes W_L)
\otimes{\rm pr}_{D_n}^*\theta_n^*{\cal L}_\psi\bigr)
[nm_0]
$$
d\'ecal\'e de $n-1+{n(n-1)(2n-1)\over 6}(g-1)$.
\endth

\rem Preuve
\endrem
Notons ${\cal D}_n$ le champ des drapeaux
$$
{\cal P}_{n,\bullet}=\bigl((0)={\cal P}_{n,0}\subset {\cal P}_{n,1}
\subset\cdots\subset {\cal P}_{n,n}\bigr)
$$
de fibr\'es vectoriels sur $X$, munis, pour chaque $j=1,\ldots ,n$,
d'un isomorphisme
$$
{\rm gr}_j{\cal P}_{n,\bullet}\cong\Omega^{n-j}.
$$
Si on note
$$
\rho_n :D_n\rightarrow {\cal D}_n
$$
le morphisme d\'efini par le drapeau universel
$\widetilde{\cal P}_{n,\bullet}$ sur $X\times D_n$
et si on note
$$
\sigma_n :{\cal D}_n\rightarrow D_n
$$
le morphisme qui envoie un objet sur sa classe d'isomorphie, alors
$\rho_n$ identifie ${\cal D}_n$  au quotient du sch\'ema affine
$D_n$ par l'action triviale du sch\'ema en groupes $U_n$ sur $D_n$
des automorphismes de $\widetilde{\cal P}_{n,\bullet}$ et
$\sigma_n$, qui est une section de $\rho_n$, identifie
$D_n$ \`a  l'espace grossier de ${\cal D}_n$. Il est facile de
v\'erifier par r\'ecurrence sur $n$ que
$U_n$ est une extension successive de copies du groupe
additif ${\bb G}_{{\rm a},D_n}$ et est de dimension
$$
\sum_{j=1}^{n-1}(n-j){\rm dim}_kH^0(X,\Omega^j)=n-1+
{n(n-1)(2n-1)\over 6}(g-1)
$$
(l'homomorphisme canonique
$$
U_n\rightarrow D_n\times_{D_{n-1}}U_{n-1}
$$
est un \'epimorphisme et son noyau n'est autre que
$$
{\bb G}_{{\rm a},D_n}\otimes_{{\cal O}_{D_n}}
{\rm pr}_{{D_n},*}{\cal H}{\it om}_{{\cal O}_{D_n}}
({\cal P}_{n,n-1},{\cal P}_{n,1})\,).
$$

Le champ
$$
{\cal E}_{n-1}^{\vee\circ}\times_{{\cal C}_{n-1}}
{\cal E}_{n-2}^{\vee\circ}\times_{{\cal C}_{n-2}}
\cdots {\cal E}_1^{\vee\circ}\times_{{\cal C}_1}
{\cal E}_0^{\vee\circ}
$$
s'identifie canoniquement au champ des triplets
$({\cal P}_{n,\bullet},{\cal M}_n,u_n)$, o\`u ${\cal P}_{n,\bullet}
\in {\cal D}_n$, o\`u ${\cal M}_n\in {\cal C}_n$ et o\`u $u_n:
{\cal P}_{n,n}\rightarrow {\cal M}_n$ est une application
${\cal O}_X$-lin\'eaire injective (pour chaque $i=0,\ldots ,n$,
le quotient
$$
{\cal M}_i\buildrel{\rm dfn}\over{=\!=}{\cal M}_n/{\cal P}_{n,n-i}
$$
est automatiquement dans ${\cal C}_i$).

L'application
$$
({\cal P}_{n,\bullet},{\cal M}_n,u_n)\mapsto {\cal P}_{n,\bullet}
$$
d\'efinit un morphisme de champs
$$
{\cal E}_{n-1}^{\vee\circ}\times_{{\cal C}_{n-1}}
{\cal E}_{n-2}^{\vee\circ}\times_{{\cal C}_{n-2}}
\cdots {\cal E}_1^{\vee\circ}\times_{{\cal C}_1}
{\cal E}_0^{\vee\circ}\rightarrow {\cal D}_n
$$
et l'application
$$
({\cal P}_{n,n}\hookrightarrow{\cal M}_n)
\mapsto \bigl(\Omega^i={\cal P}_{n,n-i}/{\cal P}_{n,n-i-1}
\hookrightarrow {\cal M}_n/{\cal P}_{n,n-i-1}\twoheadrightarrow
{\cal M}_n/{\cal P}_{n,n-i}\bigr)_{i=0,\ldots ,n-1}
$$
d\'efinit un morphisme de champs
$$
{\cal Z}_n^{\circ\circ}\rightarrow
{\cal E}_{n-1}^{\vee\circ}\times_{{\cal C}_{n-1}}
{\cal E}_{n-2}^{\vee\circ}\times_{{\cal C}_{n-2}}
\cdots {\cal E}_1^{\vee\circ}\times_{{\cal C}_1}
{\cal E}_0^{\vee\circ}.
$$
Le carr\'e
$$
\matrix{{\cal Z}_n^{\circ\circ}&\longrightarrow &D_n\cr
\noalign{\smallskip}
\Big\downarrow &&\mapdownright{\rho_n}\cr
\noalign{\smallskip}
{\cal E}_{n-1}^{\vee\circ}\times_{{\cal C}_{n-1}}
{\cal E}_{n-2}^{\vee\circ}\times_{{\cal C}_{n-2}}
\cdots {\cal E}_1^{\vee\circ}\times_{{\cal C}_1}
{\cal E}_0^{\vee\circ}&\longrightarrow &{\cal D}_n,\cr}
$$
o\`u la fl\`eche horizontale du haut est le compos\'e
${\rm pr}_{D_n}\circ\alpha_n^{\circ\circ}$, est cart\'esien.

Maintenant, par d\'efinition, $W_{L,n}^\circ$ se calcule comme
suit: on prend l'image inverse de
$$
\pi_0^{\vee *}W_L[m_0+\sum_{i=1}^{n-1}(m_0-i^2(g-1))]
$$
par la projection canonique de ${\cal E}_{n-1}^{\vee\circ}
\times_{{\cal C}_{n-1}}{\cal E}_{n-2}^{\vee\circ}
\times_{{\cal C}_{n-2}}\cdots {\cal E}_1^{\vee\circ}
\times_{{\cal C}_1}{\cal E}_0^{\vee\circ}$ sur
${\cal E}_0^{\vee\circ}$, on tensorise cette image inverse par
le faisceau d'Artin-Schreier obtenu en tirant $\sigma_n^*
\theta_n^*{\cal L}_\psi$ par le morphisme
${\cal E}_{n-1}^{\vee\circ}\times_{{\cal C}_{n-1}}
{\cal E}_{n-2}^{\vee\circ}\times_{{\cal C}_{n-2}}
\cdots {\cal E}_1^{\vee\circ}\times_{{\cal C}_1}
{\cal E}_0^{\vee\circ}\rightarrow {\cal D}_n$ d\'efini
ci-dessus et enfin on prend l'image directe \`a supports propres
par la projection de ${\cal E}_{n-1}^{\vee\circ}
\times_{{\cal C}_{n-1}}{\cal E}_{n-2}^{\vee\circ}
\times_{{\cal C}_{n-2}}\cdots {\cal E}_1^{\vee\circ}
\times_{{\cal C}_1}{\cal E}_0^{\vee\circ}$ sur
${\cal E}_{n-1}^{\vee\circ}\cong {\cal E}_n^\circ$.

Il ne reste plus qu'\`a remarquer que le foncteur
$$
R\rho_{n,!}\rho_n^*(\cdot )[{\rm dim}_kU_n]
$$
n'est autre que le foncteur
$$
(\cdot )[-{\rm dim}_kU_n]=(\cdot )[-(n-1)-{n(n-1)(2n-1)\over 6}
(g-1)].
$$
\hfill\hfill\cqfd

\vfill\eject
\centerline{\douzebf Bibliographie}
\vskip 20mm
%-------------------------------
\parindent=10mm

\newtoks\ref
\newtoks\auteur
\newtoks\titre
\newtoks\editeur
\newtoks\annee
\newtoks\revue
\newtoks\tome
\newtoks\pages
\newtoks\reste
\newtoks\autre

\def\livre{\leavevmode
\llap{[\the\ref]\enspace}%
\the\auteur\pointir
{\sl \the\titre},
\the\editeur,
{\the\annee}.
\smallskip
\filbreak}

\def\article{\leavevmode
\llap{[\the\ref]\enspace}%
\the\auteur\pointir
\the\titre,
{\sl\the\revue}
{\bf\the\tome},
({\the\annee}),
\the\pages.
\smallskip
\filbreak}

\def\autre{\leavevmode
\llap{[\the\ref]\enspace}%
\the\auteur\pointir
\the\reste.
\smallskip
\filbreak}

%-----------------------------------------------------

\ref={Br}
\auteur={J.-L. {\pc BRYLINSKI}}
\titre={Transformations canoniques, dualit\'e projective,
th\'eorie de Lefschetz, transformations de Fourier et sommes
trigonom\'etriques {\it dans} G\'eom\'etrie et Analyse
Microlocales}
\revue={Ast\'erisque}
\tome={140-141}
\annee={1986}
\pages={3-134}
\article

\ref={Dr}
\auteur={V.G. {\pc DRINFELD}}
\titre={Two dimensional $\ell$-adic representations of the
fundamental group of a curve over a finite field and
automorphic forms for $GL(2)$}
\revue={Amer. J. Math.}
\tome={105}
\annee={1983}
\pages={85-114}
\article

\ref={La}
\auteur={G. {\pc LAUMON}}
\titre={Correspondance de Langlands g\'eom\'etrique
pour les corps de fonctions}
\revue={Duke Math. J.}
\tome={54}
\annee={1987}
\pages={309-359}
\article

\ref={Ma}
\auteur={J.G. {\pc MACDONALD}}
\titre={Symmetric Functions and Hall Polynomials}
\editeur={Clarendon Press, Oxford}
\annee={1979}
\livre

\ref={Ne}
\auteur={P.E. {\pc NEWSTEAD}}
\titre={Introduction to Moduli Problems and Orbit Spaces}
\editeur={Tata Institute of Fundamental Research, Bombay}
\annee={1978}
\livre

\ref={Sh}
\auteur={J.A. {\pc SHALIKA}}
\titre={The multiplicity one theorem for $GL_n$}
\revue={Ann. of Math.}
\tome={100}
\annee={1974}
\pages={171-193}
\article

\vskip 10mm
\let\+=\tabalign

\line{\hbox{\kern 5mm{\vtop{\+ G\'erard LAUMON\cr
\+ URA 752 du CNRS\cr
\+ Universit\'e de Paris-Sud\cr
\+ Math\'ematiques, b\^at. 425\cr
\+ 91405 ORSAY C\'edex (France)\cr}}}}

\bye